\documentstyle[10pt]{article}
\input amssymb.sty
\input epsf
\newtheorem{theorem}{Theorem}[section]

\newtheorem{lemma}[theorem]{Lemma}
\newtheorem{proposition}[theorem]{Proposition}
\newtheorem{corollary}[theorem]{Corollary}
\newtheorem{conjecture}[theorem]{Conjecture}
\newtheorem{definition}[theorem]{Definition}

\begin{document}
\newcommand{\Z}{{\Bbb Z}}
\newcommand{\R}{{\Bbb R}}
\newcommand{\Q}{{\Bbb Q}}
\newcommand{\C}{{\Bbb C}}
\newcommand{\lms}{\longmapsto}
\newcommand{\lra}{\longrightarrow}
\newcommand{\ra}{\rightarrow}

\begin{titlepage}
\title{Multiple  polylogarithms and mixed Tate motives }
\author{A.B. Goncharov }

\date{}
\end{titlepage}
\maketitle

\vskip 3mm

{\it $\qquad$$\qquad$$\qquad$$\qquad$$\qquad$$\qquad$ To Don Zagier for his 50th birthday}

\vskip 5mm

\tableofcontents

\section  {Introduction}

{\bf An outline}. In this paper, which continues the series [G0-5],  we develop the theory of 
multiple polylogarithms 
\begin{equation} \label{zhe5}
{\rm Li}_{n_{1},...,n_{m}}(x_{1},...,x_{m})   
\quad = \quad 
\sum_{0 < k_{1} < k_{2} < ... < k_{m} } \frac{x_{1}^{k_{1}}x_{2}^{k_{2}}
... x_{m}^{k_{m}}}{k_{1}^{n_{1}}k_{2}^{n_{2}}...k_{m}^{n_{m}}},
\end{equation}
from 
 analytic, Hodge-theoretic and motivic points of view. 

Let $\mu_N$ be  the group of all $N$-th roots of unity. 
One of the reasons to develop such a theory was a mysterious correspondence (loc. cit.) 
between 
the 
 structure of the motivic fundamental group $\pi_1^{{\cal M}}({\Bbb G}_m - \mu_N)$ and 
the geometry and topology of modular varieties 
$$
Y_1(m;N):= \Gamma_1(m;N) \backslash GL_m(\R)/\R^*_+\cdot O_m 
$$  
Here $\Gamma_1(m;N) \subset GL_m(\Z)$ is the subgroup fixing  vector 
$(0, ..., 0, 1)$ modulo $N$. This paper  provides a background 
for the investigation of this correspondence. Indeed, as we show below, 
the study of 
motivic multiple polylogarithms at  $N$-th roots of unity is 
equivalent to the study of 
the motivic torsor ${\cal P}^{{\cal M}}({\Bbb G}_m - \mu_N; v_0, v_1)$ 
of path 
on ${\Bbb G}_m - \mu_N$ between the natural tangential base points at $0$ and $1$. 

The Hodge side of the story looks as follows. 
First of all we show that multiple polylogarithms are periods of $\Q$-rational framed 
Hodge-Tate structures, that is 
 mixed Hodge structures such that $h^{p,q}=0$ if $p \not = q$. 
The framing is an additional data which allows  to consider  
a specific period of a mixed Hodge structure. If we want to keep the information 
about this period only we come to an  equivalence relation on the set of all framed 
Hodge-Tate structures. 
The equivalence classes form a commutative graded Hopf algebra ${\cal H}_{\bullet}$ 
over $\Q$. The product structure is compatible with the product of the periods. The coproduct 
is something really new: it is invisible on the level of numbers. 
We show that the framed Hodge-Tate structures corresponding to 
multiple polylogarithms form a Hopf subalgebra ${\cal Z}^{\cal H}_{\bullet}(\C^*)$ of 
${\cal H}_{\bullet}$, and compute it very explicitly. 
Moreover, for any subgroup $G \subset \C^*$ the framed Hodge-Tate structures 
 related to multiple polylogarithms with the  arguments  in $G$ form a Hopf subalgebra  
${\cal Z}^{\cal H}_{\bullet}(G)$. In particular when $G = \mu_N$ we get the Hopf algera
${\cal Z}^{\cal H}_{\bullet}(\mu_N)$.  
We call it the {\it level $N$ cyclotomic Hopf algera} and  consider its study 
as the subject of higher cyclotomy theory. Indeed, let
$$
S_N:= {\rm Spec}\frac{\Z[x]}{(x^N-1)}[\frac{1}{N}]
$$
be the scheme of $N$-cyclotomic integers.
The   first  nontrivial component of 
 ${\cal Z}^{\cal H}_{\bullet}(\mu_N)$ is 
\begin{equation} \label{1.20.01.2}
{\cal Z}^{\cal H}_{1}(\mu_N) \quad = \quad (\mbox{the group of 
cyclotomic units in ${\cal O}^*(S_N) $}) \otimes \Q
\end{equation} 
The analytic incarnation of (\ref{1.20.01.2}) is provided by the 
formula ${\rm Li}_1(\zeta_N) = - \log (1-\zeta_N)$.

The Hopf algebra ${\cal Z}^{\cal H}_{\bullet}(\mu_N)$ describes the 
structure of 
the Hodge realization ${\cal P}^{{\cal H}}({\Bbb G}_m - G; v_0, v_1)$ of the torsor of path.
In particular we show that all the periods of this mixed Hodge structure 
are expressed via the special values of multiple polylogarithms at $N$-th roots of unity. 

The relationship of the 
structure of the  Hopf algebra ${\cal Z}^{\cal H}_{\bullet}(\mu_N)$ 
with the geometry of modular varieties $Y_1(m;N)$ 
will be discussed in detail [G10], see [G4] for the l-adic version. 

Now let us turn to the motivic aspect of the story. 
M. Levine [L],  V. Voevodsky [V], and also M.Hanamura [H], independently, 
 constructed triangulated categories of motives over an arbitrary  field.  
Unfortunately  so far a construction similar to [L] or [V] of the 
triangulated category of 
motives over a 
Noetherian scheme, or even over ${\rm Spec}\Z$,  with the desired ${\rm Ext}$-groups, 
is still far from being available. 

Nevertheless,  using the results of [L] or  [V] combined with some other results 
and Tannakian formalism, we construct in this paper 
  the abelian category ${\cal M}_T({\cal O}_{F,S})$ of mixed Tate motives over 
a ring ${\cal O}_{F,S}$ of $S$-integers in a number field $F$. 
 It has all the desired properties, including the basic formulas 
(\ref{1*1=}) - (\ref{1*1==}) 
expressing 
the ${\rm Ext}$-groups  via the K-groups of ${\cal O}_{F,S}$. 

We define the motivic torsor  
of path ${\cal P}^{\cal M}(X; x,y)$ on a regular variety $X$ over a field $F$. 
In particular if $F$ is a number field,  $X = {\Bbb A}^1 - \{z_1, ..., z_m\}$
 and  $ z_i, x, y \in F$, it is a pro-object in ${\cal M}_T(F)$. 

 In the continuation of this  paper we will show that the 
motivic torsor ${\cal P}^{{\cal M}}({\Bbb G}_m - \mu_N; v_0, v_1)$ 
of path 
on ${\Bbb G}_m - \mu_N$ between the tangential base points at $0$ and $1$ is a 
pro-object in the category ${\cal M}_T(S_N)$. 
This, combined with formulas (\ref{1*1=}) - (\ref{1*1==}),  
implies strong estimates on the dimensions of the $\Q$-vector 
spaces generated by the all periods of this motive, which are hypothetically exact 
when $N=1$. 
The torsor of path ${\cal P}({\Bbb P}^1 - \{0, 1, \infty\}; v_0, v_1)$ 
has been considered from motivic
point of view in the fundamental paper by P. Deligne [D]. Notice that in  
that paper the motives were 
 understood via their realizations, so the formulas for the Ext-groups 
were unaccessible.

We formulate conjectures describing the structure of the motivic Galois group of 
the category of mixed Tate motives over an arbitrary field $F$ 
via the motivic multiple polylogarithms.

{\bf 1. The multiple $\zeta$ - values}. They were invented by L.Euler in the
correspondence with Goldbach [Eu] as the  power series
\begin{equation} \label{lb*}
\zeta(n_{1}, n_2, ...,n_{m}) : = \sum_{0 < k_{1} < k_{2} < ... < k_{m} }
\frac{1}{k_{1}^{n_{1}}k_{2}^{n_{2}}...k_{m}^{n_{m}}} \qquad n_m >1
\end{equation}
  Here $w := n_1+...+n_m$ is the {\it  weight} of (\ref{lb*}) and $m$
is the {\it  depth}. 
Euler discovered some linear relations
over ${\Bbb Q}$ among them. Namely, he proved that   if $m+n$ is odd then 
$\zeta(m,n)$ is a linear combination with rational coefficients
of $\zeta(m+n)$ and  $\zeta(k)\zeta(m+n-k)$. 
Then these numbers were neglected, see however the work of 
M. Hoffman [Hof]  and references there.

In the very end of 80'th they were resurrected as coefficients of 
Drinfeld's associator [Dr]; 
couple of years later 
rediscovered by D.Zagier [Z] who, in particular,  found and studied 
the double shuffle relations for the multiple zeta numbers; 
 appeared in the works of 
M. Kontsevich [K1] on knot invariants,  and the author [G0-1]
 on mixed Tate motives. Later on, in the mid of 90'th 
D. Broadhurst and D. Kreimer [Kr], [BK] discovered them 
in quantum field theory. As coefficients of Drinfeld's associator appear
 in quantization problems ([K2]).  
In fact in [Dr] and [K1] appeared  not power
series (\ref{lb*}) but the so-called Drinfeld integrals related to 
the multiple $\zeta$-values 
by the following formula first noticed by Kontsevich ([K1]):
\begin{equation} \label{5*5}
{  \zeta}(n_{1},...,n_{m})  =  
 \int_{0}^{1} \underbrace {\frac{dt}{1-t}  \circ \frac{dt}{t} \circ
... \circ 
\frac{dt}{t} }
_{ n_{1} \quad  \mbox {times}} \circ \quad
... \quad \circ \underbrace 
{\frac{dt}{1-t} \circ \frac{dt}{t} \circ ...
\circ \frac{dt}{t}}_{ n_{m} \quad \mbox  {times}}
\end{equation}

The expression on the right is an iterated integral. The general definition of iterated integrals 
looks as follows. 
Let $\omega_{1},...,\omega_{n}$ be one-forms on a manifold $M$ and $\gamma : [0,1] \rightarrow M$ be a smooth path. 
The iterated integral $\int_{\gamma} \omega_{1} \circ ... \circ \omega_{n}$ 
is defined inductively:
\begin{equation} \label{h1}
\int_{\gamma} \omega_{1} \circ ... \circ \omega_{n} \quad : = \quad \int_{0}^{1}(\int
_{\gamma_{t}} \omega_{1} \circ ... \circ \omega_{n-1
}) \gamma^{*} \omega_{n}
\end{equation}
Here $ \gamma_{t}$ is the restriction of $\gamma$ to the interval $[0,t]$
and
$\int_{\gamma_{t}} \omega_{1} \circ ... \circ \omega_
{n-1}$ is considered as a function on $[0,1]$. We multiply this function by the
 one-form $\gamma^{*} \omega_{n}$  and integrate.

In the depth $1$ case (\ref{5*5}) is the famous Leibniz formula which shows that $\zeta(n)$ 
is an iterated integral of length $n$. For example
$$
\zeta(3) = \int_{0 \leq t_1 \leq t_2 \leq t_3 \leq 1}\frac{dt_1}{1-t_1}\wedge \frac{dt_2}{t_2}\wedge \frac{dt_3}{t_3}
$$

{\bf 2. Multiple polylogarithms}. They were defined in [G0-1] by the power series (\ref{zhe5}) 
which obviously generalize both
classical polylogarithms ${\rm Li}_n(x)$ (if $m=1$) 
and  multiple $\zeta$-values (when $x_1 =  ... = x_m =1$). 
We show that multiple polylogarithms are periods of 
 variations of mixed Tate motives, 
and hypothetically they provide  all such period functions, see conjecture \ref{1.4} below. 
This was the main reason for me to study them. 
 
The study of multiple polylogarithms at roots of unity seems to be especially interesting 
because of the 
 link to the  geometry 
of modular varieties $Y_1(m;N)$, as was outlined above.  

{\bf 3.  Conjectural description of the algebra of multiple $\zeta$-values}. 
Let  ${\cal Z}$   be the space
of 
  $\Bbb Q$ -  linear combinations of multiple $\zeta$'s. ${\cal Z}$   
is a  commutative algebra 
 over $\Bbb Q$.  For instance 
$$
  \zeta (m)\cdot  \zeta (n) =   \zeta (m,n) +   \zeta(n,m) +  \zeta (m+n)
$$
because
$$
\sum_{k_1,k_2 >0}\frac{1}{k_1^{ m}k_2^{n}} = \Bigl(\sum_{0 < k_1 < k_2 } + \sum_{0 < k_1 = k_2 } + \sum_{ k_1  > k_2 > }\Bigr)\frac{1}{k_1^{ m}k_2^{n}}
$$

 New  relations  were 
found by   
D.Zagier. 
Surprisingly
the dimension of the space of cusp forms for $SL_2(\Bbb Z)$ 
showed up in his study of  the depth 2 double shuffle relations [Z], [Z1]. 
Shortly after Kontsevich [K1]
showed that the new relations are coming from the product formula for 
 the iterated integrals (\ref{5*5}).

The main problem  is to describe explicitly the structure of 
the algebra ${\cal Z}$, i.e. to find 
all algebraic relations over $\Bbb Q$ between multiple $\zeta$ -values. 
Here is a hypothetical answer ([G1]). Denote by ${\cal Z}_k$ the $\Q$-vector space 
spanned over $\Q$ by the weight $k$ multiple zetas.

Let ${\cal F}(3,5,...)_{\bullet}$ be the free graded Lie algebra
generated by elements $e_{-(2n+1)}$ of degree $-(2n+1)$ ($n \geq 1$) and 
$$
U{\cal F}(3,5,...)_{\bullet}^{\vee}:= \oplus_{n \geq 1} 
\Bigl(U{\cal F}(3,5,...)_{-(2n+1)}\Bigr)^{\vee} 
$$ 
be the graded dual to its universal
enveloping algebra. It is graded by positive integers.   

Let $f_3, f_5, ...$ be the functionals on the vector space generated by the vectors $e_{-3}, e_{-5}, ... $ such that $<f_i,e_{-j}> =\delta_{i,j}$.  Then  $U{\cal
F}(3,5,...)_{\bullet}^{\vee}$ is isomorphic to the space of
noncommutative polynomials on $f_{2n+1}$ with the shuffle product.

\begin{conjecture} \label{Conjecture A}. 
a) The  weight provides a grading on the algebra ${\cal
Z}$.   

b) One has an isomorphism of graded algebras over $\Q$, where $\pi^2$ has the degree 2.
$$
{\cal Z}_{\bullet} \quad \stackrel{?}{=} \quad \Q[\pi^2] \otimes_{\Q} U{\cal F}(3,5,...
)_{\bullet}^{\vee} 
$$
\end{conjecture}

The part a) means that relations between $\zeta$'s of different weight, like $\zeta(5) = 
\lambda \cdot \zeta(7)$ where $\lambda \in \Q$, are impossible.

Let ${\cal F}(2,3)_{\bullet}$ be the free graded Lie algebra generated by two 
elements of degree $-2$ and $-3$.  The graded dual $U{\cal F}(2,3)_{\bullet}^{\vee}$ 
is isomorphic as a
graded vector space  to the space of noncommutative
polynomials in two variables $p$ and $g_3$ of degrees 2 and 3.
There is canonical isomorphism of {\it graded vector spaces}
$$
\Q[\pi^2] \otimes U{\cal F}(3,5,...
)_{\bullet}^{\vee}  \quad = \quad U{\cal F}(3,5)_{\bullet}^{\vee}
$$
The rule is clear from the pattern
$
(\pi^2)^3 f_3(f_7)^3(f_5)^2 \longrightarrow p^3 g_3(g_3p^2)^3(g_3p)^2
$. 
In particular if $d_k := {\rm dim} {\cal Z}_{k}$ then one should have $d_k =
d_{k-2} + d_{k-3}$. Computer calculations of D.Zagier [Z] confirmed this prediction 
for $k\leq 12$. Later on  much more extensive calculations  
were made by D. Broadhurst [Br].

One may reformulate  \ref{Conjecture A} as a statement about the space of 
irreducible multiple zeta values:
\begin{equation} \label{lbb}
{\overline {\cal Z}}_{ >0} :=  \frac{  {\cal Z}_{>0}}{ {\cal Z }_{>0} \cdot
   {\cal Z}_{>0}} \quad \stackrel{?}{=} \quad <\pi^2> \oplus {\cal F}(3,5,...)_{\bullet}^{\vee}
\end{equation}
Here $<\pi^2>$ is a one dimensional $\Q$-vector space generated by $\pi^2$ and 
the augmentation ideal ${\cal Z }_{>0}$  is generated by $\pi^2$ and $\zeta(n_1, ..., n_m)$.  
The left hand side in (\ref{lbb}) the space of algebra generators in
${\cal Z }_{\bullet}$. 
We get the following list of conjectural dimensions of weight $ w$
pieces of the left hand side of (\ref{lbb}):
$$
w \qquad : \quad 1 \quad 2 \quad   3 \quad 4 \quad 5\quad 6 \quad 7 \quad 8
\quad 9 \quad 10 \quad 11 \quad 12 \quad 13 \quad 14 \quad 15
$$
$$
{\rm dim} \quad: \quad 0 \quad 0 \quad 1 \quad 0 \quad 1 \quad 0 \quad 1 \quad 1
\quad 1 \quad 1 \qquad 2 \quad 2 \qquad 3 \quad 3 \qquad 4
$$

For example  one should have   ${\rm dim} \overline Z_{ 11} = 2$
     because ${\cal F}(3,5,...)_{11}$ is spanned 
over $\Q$ 
by  $e_{11}$ and $[[e_{3},e_{5}],e_3]$. Observe that there are $2^9=512$ convergent multiple 
$\zeta$-values of weight $11$. So there are a lot of  relations between them. The following 
theorem implies that ${\rm dim} \overline Z_{ 11} \leq  2$.

 {\bf 4. The Hodge version of the multiple zeta Hopf algebra}. Denote by 
$\widetilde {\cal Z}_w$ the $\Q$-span of the numbers 
$(2\pi i)^{-w}\zeta(n_1, ..., n_m)$ of weight $w$. 
We add $\Q$ as the weight zero component.  Then 
$\widetilde {\cal Z}:= \cup \widetilde {\cal Z}_w$ is a commutative algebra filtered 
by the weight filtration $W_{\bullet}$. Observe that there is no weight grading since $(2\pi i)^{-2}\zeta(2) = - 1/24$. 
Here is a version of conjecture \ref{Conjecture A}:

$$
{\rm Gr}_{\bullet}^W\widetilde {\cal Z} \quad \stackrel{?}{=} \quad 
 U{\cal F}(3,5,...
)_{\bullet}^{\vee} 
$$
The right hand side of this formula has a Hopf algebra structure. 
This suggests that  ${\rm Gr}_{\bullet}^W\widetilde {\cal Z}$ should have 
a natural structure of a commutative graded Hopf algebra. This is 
far away from being proved. 
However replacing the multiple zeta numbers by more sofisticated objects 
we start to see the Hopf algebra structure. Let me explain how to do this. 

By definition (see chapter 3 or [BGSV] for more details) an $n$-framing on a Hodge-Tate structure $H$ consists of a 
choice of non zero morphisms
$$
 v: \Q(0) \lra {\rm Gr}_{0}^WH; \qquad f: {\rm Gr}_{-2m}^WH \lra \Q(m)
$$
Consider the coarsest equivalence relation $\sim$ on the set of all 
$n$-framed Hodge-Tate structures such that $(H_1, v_1, f_1) \sim (H_2, v_2, f_2)$ if there 
is a morphism of Hodge-Tate structures $H_1 \to H_2$ respecting the frames. 
The equivalence classes 
of $n$-framed Hodge-Tate structures form an abelian group ${\cal H}_{n}$, and 
${\cal H}_{\bullet}:= \oplus_{n \geq 0}{\cal H}_{n}$ has a commutative 
graded Hopf algebra structure.

Now we can explain the precise relationship between the multiple zeta values and 
the Hodge realization 
\begin{equation} \label{3.6.01.7}
{\cal P}^{\cal H}(\C P^1 - \{0, 1, \infty\}; v_0, v_1)
\end{equation}
 of the torsor of path 
on $\C P^1-\{0, 1, \infty\}$ between the tangential base points 
at $0$ and $1$ provided by the  natural parameter $t$ on $\C P^1$. 
The mixed Hodge structure  (\ref{3.6.01.7}) is  a projective limit 
of the Hodge-Tate structures (see [D] or chapter 4 below). 
It carries a weight filtration $W_{\bullet}$ such that 
\begin{equation} \label{3.6.01.11}
{\rm Gr}_{-2m}^W{\cal P}^{\cal H}(\C P^1 - \{0, 1, \infty\}; v_0, v_1) \quad = \quad \otimes^m 
H_1(\C P^1 - \{0, 1, \infty\})
\end{equation}
is a direct sum of $2^m$ copies of $\Q(m)$. 
The forms $\omega_0:= d\log (t)$ and $\omega_1:= d\log (1-t)$ form a natural basis in 
$H_{DR}^1(\C P^1 - \{0, 1, \infty\})$. Thus the elements
$$
\omega(\varepsilon_1, ..., \varepsilon_m):= \omega_{\varepsilon_1} \otimes ... \otimes \omega_{\varepsilon_m}; \qquad \varepsilon_i 
\in \{0,1\}
$$
form a basis in the dual to (\ref{3.6.01.11}). We denote by 
$\gamma(\varepsilon_1, ...,\varepsilon_m)$ the corresponding dual basis in 
(\ref{3.6.01.11}).

\begin{definition} \label{3.6.01.2} We define 
\begin{equation} \label{3.6.01.12}
\zeta^{\cal H}(n_1, ..., n_m) \in {\cal H}_{w}
\end{equation}
as  the 
Hodge-Tate structure (\ref{3.6.01.7}) framed by 
\begin{equation} \label{3.6.01.5}
\underbrace {\frac{dt}{1-t}  \otimes \frac{dt}{t} \otimes
... \otimes
\frac{dt}{t} }
_{ n_{1} \quad  \mbox {times}} \otimes \quad
... \quad \otimes \underbrace 
{\frac{dt}{1-t} \otimes \frac{dt}{t} \otimes ...
\otimes \frac{dt}{t}}_{ n_{m} \quad \mbox  {times}} \quad \mbox{and} \quad \gamma_{\emptyset}
\end{equation}
\end{definition}

Formula (\ref{5*5}) just means that $\zeta(n_1, ..., n_m)$ is a period of 
$\zeta^{\cal H}(n_1, ..., n_m)$.

\begin{theorem} \label{3.6.01.3} The Hodge-Tate structure 
${\cal P}^{\cal H}(\C P^1 - \{0, 1, \infty\}; v_0, v_1)(-n)$ with an arbitrary framing given 
by 
$\gamma(\varepsilon_1, ...,\varepsilon_n)$ and $\omega(\delta_1, ..., \delta_m) $
is equivalent to a product of the multiple zeta 
structures (\ref{3.6.01.12}). 

In particular all periods of the mixed Hodge structure (\ref{3.6.01.7}) are  $\Q$-linear 
combinations 
of multiple zeta numbers multiplied by appropriate powers of $2\pi i$. 
\end{theorem}

It follows immediately from this that the elements (\ref{3.6.01.12}) form a 
graded Hopf algebra denoted ${\cal Z}^{\cal H}_{\bullet}$ and called the multiple zeta 
Hopf algebra.

\begin{conjecture} \label{1.3} There exists  an isomorphism of graded 
Hopf algebras over $\Q$
$$
{{\cal Z}}^{\cal H}_{\bullet} \quad \stackrel{?}{=}\quad  U{\cal F}(3,5, ... )_{\bullet}^{\vee}
$$
\end{conjecture}
It is the Hodge version of conjecture 17b) in [G1], as well as   
Conjecture \ref{Conjecture A}.

The  $l$-adic version of Conjecture \ref{Conjecture A} is related to the 
questions raised by P. Deligne, Y.Ihara
([Ih3]) and V.G.Drinfeld ([Dr]). 

Theorem (\ref{3.6.01.3}) implies the following result:

\begin{theorem} \label{3.6.01.1} There exists canonical surjective homomorphism of commutative graded algebras ${\cal Z}^{\cal H}_{\bullet} \lra {\rm Gr}_{\bullet}^W \widetilde 
{\cal Z}$. Thus  ${\rm dim} {\rm Gr}_k^W \widetilde {\cal Z}\leq 
{\rm dim} {\cal Z}^{\cal H}_k $. 
\end{theorem}

It follows that conjecture \ref{Conjecture A} implies Conjecture \ref{1.3}.

It is hard to estimates ${\rm dim}   {\cal
Z}_{k}$ from below: nobody can  prove that  
${ \zeta}(5) \not \in \Q$. The works of Apery on $\zeta(3)$ and Rivoal  
are amazing  breakthroughs in this direction. 
Our point is that 
it is much easier to deal with  framed 
Hodge-Tate structures then with real numbers! For example it is obvious that 
there are no linear relations between 
$\zeta^{\cal H}(2n+1)$'s.  So we eliminated the transcendental part of the problem by  working
with more  sophisticated   Hodge multiple zetas. I think that 
the  treatment of the transcendental aspects of multiple zetas is impossible without 
  investigation of the corresponding Hodge/motivic objects.

To explain where the right hand side  in the conjectures above is coming from 
we need  mixed motives over ${\rm Spec}(\Z)$.

   {\bf 5. Mixed Tate motives of the ring of $S$-integers in a number field}. 
Let $F$ be a number field and $S$ a finite set of prime ideals in  
the ring of integers in $F$. Denote by ${\cal O}_{F, S}$ 
the localization of the ring of integers in $S$. 
In section 3 we construct the abelian category ${\cal M}_{T}({\cal O}_{F, S})$ 
of mixed Tate motives over  
$Spec({\cal O}_{F, S})$. It contains the Tate motives $\Q(n)_{\cal M}$ which are mutually non isomorphic. The ${\rm Ext}$-groups between the 
Tate motives are calculated by the following basic formula conjectured by Beilinson: 
\begin{equation} \label{1*1=}
Ext^1_{{\cal M}_{T}({\cal O}_{F, S})}(\Q(0)_{\cal M} , \Q(n)_{\cal M}) \quad 
= \quad K_{2n-1}({\cal O}_{F, S})\otimes \Q
\end{equation} 
\begin{equation} \label{1*1==}
Ext^i_{{\cal M}_{T}({\cal O}_{F, S})}(\Q(0)_{\cal M} , \Q(n)_{\cal M}) \quad 
= 0 \quad \mbox{for $i>1$}
\end{equation}
The category ${\cal M}_{T}({\cal O}_{F, S})$ is equivalent to 
the category of finite dimensional modules over
 a pro-algebraic group $G_{{\cal M}_{T}}({\cal O}_{F, S})$ over $\Q$, 
called the motivic Galois group
of this category. This group is a semidirect product of the multiplicative group
${\Bbb G}_m$ and a pro-unipotent algebraic group $U({\cal O}_{F, S})$ over $\Q$. 
The Lie algebra of this prounipotent group is 
denoted by 
${   L}_{\bullet}({\cal O}_{F, S})$. The action of ${\Bbb G}_m$ provides it with a structure of 
graded Lie algebra. The category ${\cal M}_{T}({\cal O}_{F, S})$ 
is canonically equivalent to the category of finite dimensional graded 
${   L}_{\bullet}({\cal O}_{F, S})$ - modules. Denote by $H_n^i$ the degree $n$ part of $H^i$. 
It follows from this and formula (\ref{1*1==}) that 
$$
H_n^1({   L}_{\bullet}({\cal O}_{F, S})) = 
K_{2n-1}({\cal O}_{F, S})\otimes \Q; \qquad 
H_n^i({   L}_{\bullet}({\cal O}_{F, S})) = 0 \quad \mbox{for $i>1$}
$$
Thus ${   L}_{\bullet}({\cal O}_{F, S})$ is  isomorphic to a  free graded Lie algebra 
generated by the finite dimensional $\Q$-vector spaces 
$(K_{2n-1}({\cal O}_{F, S})\otimes \Q)^{\vee}$ in degree $-n$. 

{\bf Example 1}.  One has ([Bo])
\begin{equation} \label{kg}
{\rm dim} K_{2n-1}(\Z)\otimes \Q = 
\left\{ \begin{array}{ll}
1 &: \quad n \quad odd
\\ 
0 &: \quad n \quad even \end{array}\right.
\end{equation}
So the  Lie algebra ${   L}_{\bullet}(\Z)$ is 
isomorphic to ${\cal F}(3,5, ... )_{\bullet}$, and 
the category of mixed Tate motives over ${\rm Spec}(\Z)$ is equivalent to 
the category of 
graded ${\cal F}(3,5, ... )_{\bullet}$-modules.

{\bf Example 2}. Let $p(N)$  be the number of prime factors of $N$. One has ([Bo])
\begin{equation} \label{kgt}
{\rm dim} K_{2n-1}(S_N)\otimes \Q = 
\left\{ \begin{array}{ll}
\varphi(N)/2 &: \quad N>2, n >1
\\ 
\varphi(N)/2  +p(N)-1 &: \quad N>2, n =1 
\\
1 & : \quad N=2\end{array}\right.
\end{equation}

The category ${\cal M}_{T}({\cal O}_{F, S})$ is equipped with an array of 
realization functors including the Hodge and \'etale realizations.

There is a  notion of a 
framing on a mixed motive and  an equivalence relation on the set of all 
framed  mixed motives.  Informally  the equivalence classes are 
algebraic versions of  periods of mixed motives. These equivalence classes 
 form a Hopf algebra in the category of all pure motives (see ch. 3 of [G9]). 
In particular, the set of equivalence classes of  
framed mixed Tate motives over $\Z$ 
has a structure of a graded, commutative Hopf algebra over $\Q$, which is  
isomorphic to the graded dual of universal enveloping of ${   L}_{\bullet}(\Z)$. 
We will show in the continuation of this   paper that 
$\zeta^{\cal H}(n_1, ..., n_m)$ is a Hodge realization of 
the framed  mixed Tate motive $\zeta^{\cal M}(n_1, ..., n_m)$ over $\Z$. 

{\bf 6. Generalizations}. Let $\zeta_N$ be a primitive $N$-th root of unity. We define the framed Hodge-Tate structures
\begin{equation} \label{3.6.01.112}
{\rm Li}^{\cal H}_{n_1, ..., n_m}(\zeta_N^{\alpha_1}, ..., \zeta_N^{\alpha_m})
\end{equation}
as the  Hodge-Tate structure ${\cal P}^{\cal H}(\C P^1 - \{0, \mu_N, \infty\}; v_0, v_1)$ 
with the  framing suggested by formula (\ref{5*}) below, so that  
${\rm Li}_{n_1, ..., n_m}(\zeta_N^{\alpha_1}, ..., \zeta_N^{\alpha_m})$ is the main period 
of (\ref{3.6.01.112}).

 \begin{theorem} \label{3.6.01.33} The Hodge-Tate structure 
${\cal P}^{\cal H}(\C P^1 - \{0, \mu_N, \infty\}; v_0, v_1)$ with an arbitrary framing 
is equivalent to a product of the 
structures  (\ref{3.6.01.112}). 

In particular all periods of the mixed Hodge structure (\ref{3.6.01.112}) are  
$\Q$-linear combinations 
of multiple polylogarithms at $N$-th roots of unity multiplied 
by appropriate powers of $2\pi i$, twisted by $\Q(n)$. 
\end{theorem}
Similar results remain valid if $\mu_N$ is replaced by any subgroup 
$G \subset \C^*$. 

It follows that the elements (\ref{3.6.01.112}) 
form a commutative graded Hopf algebra denoted
${\cal Z}_{\bullet}^{\cal H}(\mu_N)$.

Let $\widetilde {\cal Z}_w(\mu_N)$ be the $\Q$-vector space spanned over $\Q$ 
by $(2 \pi i)^{-w}$ times the values of weight $w$ 
multiple 
polylogarithms at $N$-th roots of unity, where we consider all  brunches of the 
multivalued analytic functions provided by the power series (\ref{zhe5}). 
Then the union $\widetilde {\cal Z}(\mu_N) := \cup \widetilde {\cal Z}_w(\mu_N)$ 
is a commutative algebra filtered by the weight. 
Conjecture \ref{Conjecture A} has the following generalization.

\begin{conjecture} \label{th1.2} One has an isomorphism of the commutative algebras: 
$$  
{\rm Gr}^W_{\bullet}\widetilde {\cal Z}(\mu_N) \quad = \quad   {\cal
Z}^{\cal H}_{\bullet}(\mu_N)
$$
\end{conjecture}

 \begin{theorem} \label{th1.221} There exists canonical surjective homomorphism ${\cal Z}_{\bullet}^{\cal H}(\mu_N) \lra {\rm Gr}^W_{\bullet}\widetilde {\cal Z}(\mu_N)$ of commutative 
graded algebras. Thus 
$$  
{\rm dim} {\rm Gr}^W_{k}\widetilde {\cal Z}(\mu_N)  \quad \leq  \quad  {\rm dim} {\cal
Z}_{k}^{\cal H}(\mu_N)
$$
\end{theorem}

The space of indecomposables
$$
 {\cal C}^{{\cal H}}_{\bullet}(\mu_N):= \quad 
\frac{{\cal Z}^{{\cal H}}_{>0}(\mu_N)}{
{\cal Z}^{{\cal H}}_{>0}(\mu_N) \cdot {\cal Z}^{{\cal M}}_{>0}(\mu_N)}
$$
inherits a Lie coalgebra structure with the cobracket induced by the coproduct in 
${\cal Z}^{{\cal H}}_{>0}(\mu_N)$. Its dual is the Lie algebra 
${C}^{\cal H}_{\bullet}(\mu_N)$, called the cyclotomic Lie algebra. 
It follows from the result of Deligne [D2-3] 
that $C^{\cal H}_{\bullet}(\mu_N)$ is free 
for $N=2,3,4$. Nevertheless it follows from the results of this paper and [G4], [G10] 
that $C^{\cal H}_{\bullet}(\mu_N)$ is not free for sufficiently big $N$,  for instance for 
all prime $p>5$. The  obstruction is provided by the cusp 
forms on $Y_1(2;N)$.

{\bf 7. The universality conjecture}. 
\begin{conjecture} \label{1.4}  Periods of any variation of 
 mixed Tate motives can be expressed via multiple polylogarithms.
\end{conjecture}

For a  more refined version see conjecture 17a) in [G1] and conjecture \ref{1.10.01.3} 
 below.
The motivation comes from the following result, see theorem \ref{th2.7}: 

{\it Any function on a rational algebraic variety $Y$ given by an iterated integral of rational one-forms on $Y$ can be expressed by multiple polylogarithms}.

{\bf 8. The structure of the paper.} In chapter 2 we derive basic analytic properties 
of multiple polylogarithms. 
In the chapter  3 we introduce the abelian category 
of mixed Tate motives over ${\cal O}_{F,S}$, 
and show that it has all the expected properties. Let $X$ be  a regular variety over 
a field $F$ and $x, y \in X(F)$. In chapter 4 we define the motivic torsor of path 
${\cal P}^{\cal M}(X; x,y)$ as a pro-object in the triangulated category of motives 
${\cal D}{\cal M}_F$. If 
$F$ is a number field and the motive of $X$ belongs to the 
triangulated category of mixed Tate motives,  ${\cal P}^{\cal M}(X; x,y)$ is defined as a 
pro-object in the abelian category of mixed Tate motives over $F$. 
In chapter 5 we define a variation of  Hodge-Tate structures 
whose period functions are multiple polylogarithms. In chapter  6 we compute explicitly 
the Hopf algebra formed by  the multiple polylogarithm framed Hodge-Tate structures. 
Explicit formulas for the coproduct in small depths 
are presented in the end. 

One can show ([G10]) that the associated graded 
for the depth filtration of the 
cyclotomic Lie algebra $C^{\cal H}_{\bullet}(\mu_N)$ can be  naturally realized 
 as a Lie subalgebra of 
 the dihedral Lie algebra of the group $\mu_N$ which is defined in [G3], see also 
[G4]. 
In fact this result reverses the history: 
computation of the coproduct in the multiple polylogarithm Hopf algebra  was 
originally the source of our definition 
of the cobracket in the dihedral Lie algebra. 
Using this and  results in the[G3-4] one gets  the connection 
with the geometry of modular varieties. 

In section 7 we formulate some general 
conjectures relating the multiple polylogarithms and motivic Galois groups, 
and complete the proofs of 
the theorems from the introduction.

{\bf Acknowledgements}. A considerable part of this work was done during my stay in 
Max-Planck-Institute (Bonn) in 1992 and MSRI(Berkeley) in 1992-1993. It
was supported by   
NSF Grant DMS-9022140 in MSRI and some of the results appeared in the MSRI preprint [G0]. The work on this project continued during my several fruitful 
visits to  the 
Max-Planck-Institute (Bonn) in 1996-1999. 
I would
like to these institutions for hospitality and support.  This work was 
supported and NSF grants DMS-9500010 and DMS-9800998. 

It is my pleasure to thank A.A. Beilinson, M.L. Kontsevich,  and D. Zagier for
many very helpful discussions.

\begin{center}
\hspace{4.0cm}
\epsffile{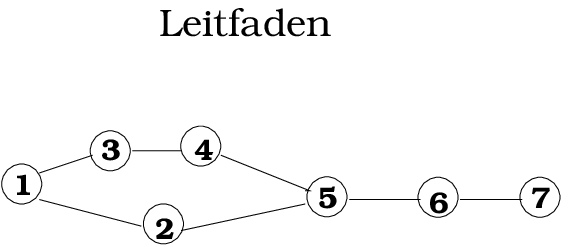}
\end{center}

\section{Analytic properties}

In this chapter we investigate the analytic properties of the 
multiple polylogarithms. Some of these results will get another 
interpretation in s. 5.3. 

{\bf 1. Definitions}. 
 Hyperlogarithms are defined as the following iterated integrals:
$$
{\rm I}(a_1: ... :a_m:a_{m+1}) \quad : = \quad \int_{0}^{a_{m+1}} \frac{dt}{t-a_{1}} \circ \frac{dt}{t-a_{2}} \circ ... \circ 
\frac{dt}{t-a_{m}} 
$$
\begin{equation} \label{2}
\int_{0 \leq t_1 \leq t_2 \leq ... \leq t_m \leq a_{m+1}} \frac{dt_1}{t_1-a_{1}} \wedge \frac{dt_2}{t_2-a_{2}} \wedge ... \wedge 
\frac{dt_m}{t_m-a_{m}} 
\end{equation}
They were considered by E.Kummer [Ku], H.Poincare, 
Lappo-Danilevsky [LD], and many others , but treated as analytic 
functions of one variable, the upper limit of integration. 
We study them as multivalued analytic functions on  $a_{1}, ... , a_{n}, a_{m+1}$.

The point $0$  is a special point for $a_{1}, ... , a_{m+1}$ in (\ref{2}).
The properties of integral (\ref{2}) change drastically if one of the
variables $a_{i}$ becomes $0$. Therefore we will study the integrals
\begin{equation} \label{3??}
{\rm I}_{n_{1},...,n_{m}}(a_{1}:...:a_{m}:a_{m+1}) : =
\end{equation}
$$
\int_{0}^{a_{m+1}} \underbrace {\frac{dt}{t-a_{1}} \circ \frac{dt}{t} \circ ... \circ 
\frac{dt}{t}}_{n_{1} \quad  \mbox {times}} \circ \underbrace {\frac{dt}{t-a_{2}} \circ \frac{dt}{t} \circ... \circ 
\frac{dt}{t}}_{n_{2} \quad  \mbox {times}} \circ \quad ... \quad \circ 
\underbrace {\frac{dt}{t-a_{m}} \circ \frac{dt}{t} \circ ... \frac{dt}{t}}_
{n_{m} \quad  \mbox {times}}
$$
as a function of integers $n_{1}, ... , n_{m}$ and complex variables 
$a_{1}, ... , a_{m}, a_{m+1}$,  where $a_{i} \neq 0$. 
More deep reasons for doing so will be apparent later on. 
Set $w:= n_1 +...+n_m$. We call integral (\ref{3??}) 
 multiple polylogarithm of weight $w$ and  depth $m$. 
Clearly 
\begin{equation} \label{3}
{\rm I}_{n_{1},...,n_{m}}(a_{1}:...:a_{m}:a_{m+1}) \quad = \quad 
{\rm I}_{n_{1},...,n_{m}}(\frac{a_{1}}{a_{m+1}}:...:\frac{a_{m}}{a_{m+1}}:1)
\end{equation}
so one can always suppose $a_{m+1} =1$.  We will use a notation
\begin{equation} \label{ur8-12}
{\rm I}_{n_{1},...,n_{m}}(a_{1},...,a_{m}) : 
={\rm I}_{n_{1},...,n_{m}}(a_{1}:...:a_{m}:1)
\end{equation}

Consider also the iterated integral
\begin{equation} \label{2**}
{\rm I}(a_0; a_1, ... ,a_m; a_{m+1}) \quad : = \quad \int_{a_0}^{a_{m+1}} \frac{dt}{t-a_{1}} \circ \frac{dt}{t-a_{2}} \circ ... \circ 
\frac{dt}{t-a_{m}}
\end{equation}
It is invariant under the affine transformations $a_i \lms \alpha a_i + \beta$, and in particular  equals to ${\rm I}(0; a_1-a_0, ... ,a_m-a_0; a_{m+1}-a_0)$.  

{\bf 2.The differential equation}. 
 \begin{theorem} \label{2.2} 
$$
d{\rm I}(a_0; a_{1},...,a_{m}; a_{m+1}) = 
$$
\begin{equation} \label{300}
\sum_{i=1}^m {\rm I}(a_0; a_{1},...,\hat a_{i},...,a_m; a_{m+1})\cdot \Bigl(
d\log(a_{i+1} - a_i) - d\log(a_{i-1} - a_i) \Bigr)
\end{equation}
\end{theorem}

{\bf Proof}. It is easy to check that one can assume, without loss of generality, that $a_0 =0, a_{m+1} =1$. 
Differentiating integral (\ref{2})  by $a_i$ we get
$$
\int_{0 \leq t_1 \leq  ... \leq  t_m \leq 1} 
\frac{dt_1}{t_1-a_1} \wedge ... \wedge \frac{dt_{i-1}}{t_{i-1} - a_{i-1}} \wedge \frac{dt_i}{(t_i-a_i)^2} 
\wedge \frac{dt_{i+1}}{t_{i+1} - a_{i+1}} \wedge ... \wedge \frac{dt_m}{t_m-a_m} 
$$
Let us assume $1 < i < m$ (the cases $i=1$ and $i=m$ are treated similarly). 
 Integrating by $t_i$ and using  
$$
\int_{t_{i-1}}^{t_{i+1}}\frac{dt_i}{(t_i-a_i)^2} = \frac{1}{t_{i-1}-a_i} - \frac{1}{t_{i+1}-a_i}
$$
and then applying the formulas 
$$
\frac{dt_{i-1}}{(t_{i-1}-a_{i-1})(t_{i-1}-a_{i})} = 
 \frac{1}{a_{i-1}-a_i}\left( \frac{dt_{i-1}}{t_{i-1}-a_{i-1}}  - 
\frac{dt_{i-1}}{t_{i-1}-a_{i}} \right)
$$
$$
-\frac{dt_{i+1}}{(t_{i+1}-a_{i+1})(t_{i+1}-a_{i})} =
 -\frac{1}{a_{i+1}-a_i}\left( \frac{dt_{i+1}}{t_{i+1}-a_{i+1}}  - 
\frac{dt_{i+1}}{t_{i+1}-a_{i}} \right)
$$
and integrating with respect to the rest of the variables we come to formula (\ref{300}). 

A formula for $d{\rm I}_{n_{1},...,n_{m}}(a_{1},...,a_{m})$ can be obtained by similar considerations. The answer can be presented as follows. 
Take the $1$-form (\ref{300}) written for 
$d{\rm I}(0; a_1, ..., a_{w}; a_{w+1})$ where $w = n_1+...+n_m$ and specialize it formally to 
the submanifold $a_2 = ... = a_{n_1} = 0$, 
$a_{n_{1}+2} = ... = a_{n_1+n_2} = 0$, ...  
using the receipt that if $a_i = a_{i+1} =0$ 
then the specialization of the $1$-form $d\log(a_i-a_{i+1})$ is zero. 
(The reasons for this rule will be apparent below when we 
 talk about canonical regularization).  
For example if all $n_i >1$ then
\begin{equation} \label{3000}
d{\rm I}_{n_1,...,n_m}(a_{1},...,a_{m}) \quad = \quad \sum_{i=1}^m {\rm I}_{n_1,...,n_i-1,...,n_m}(a_{1},...,a_m)\cdot d\log(\frac{a_{i+1}}{a_i}).
\end{equation}

{\bf 3. Power series expansion}.  
Recall that multiple polylogarithms has been defined by power series expansion 
\begin{equation} \label{5}
{\rm Li}_{n_{1},...,n_{m}}(x_{1},...,x_{m})   
\quad = \quad 
\sum_{0 < k_{1} < k_{2} < ... < k_{m} } \frac{x_{1}^{k_{1}}x_{2}^{k_{2}}
... x_{m}^{k_{m}}}{k_{1}^{n_{1}}k_{2}^{n_{2}}...k_{m}^{n_{m}}}
\end{equation}
which are convergent if $|x_{i}| < 1$. On the other hand we used the same 
name for the iterated integrals (\ref{3??}). 
Here is the presize relation justifying this.

  \begin{theorem}  \label{2.3}  
Suppose $|x_{i}| < 1$. Then 
$$
{\rm Li}_{n_{1},...,n_{m}}(x_{1},...,x_{m})  = 
$$
\begin{equation} \label{5*}
(-1)^m \cdot \int_{0}^{1} \underbrace {\frac{dt}{t-(x_1...x_m)^{-1}}  \circ 
... \circ 
\frac{dt}{t} }
_{ n_{1} \quad  \mbox {\rm times}} \circ \quad
... \quad \circ \underbrace 
{\frac{dt}{t-x_m^{-1}} \circ ... 
\circ \frac{dt}{t}}_{ n_{m} \quad  \mbox {\rm times}}
\end{equation}
\end{theorem}

This theorem provides an analytic continuation of multiple polylogarithms. 

Formula (\ref{5*}) just means that  
\begin{equation} \label{5*??}
 {\rm I}_{n_{1},...,n_{m}}(a_{1}:...:a_{m}:a_{m+1})\quad  = \quad (-1)^m \cdot {\rm Li}_{n_{1},...,n_{m}}
(\frac{a_2}{a_1},\frac{a_3}{a_2},...,\frac{a_m}{a_{m-1}},\frac{a_{m+1}}{a_m})  
\end{equation}
\begin{equation} \label{5*??**??}
 {\rm Li}_{n_{1},...,n_{m}}(x_{1},...,x_{m} )\quad  = \quad (-1)^m \cdot 
{\rm I}_{n_{1},...,n_{m}}
( (x_1...x_m)^{-1},(x_2...x_m)^{-1},...,x_m^{-1})  
\end{equation}

\vskip 3mm \noindent
{\bf Proof of Theorem  \ref{2.3}}. Suppose for simplicity that $m = 2$. Then integral (\ref{5*}) is equal to
$$
(-1)^2\cdot\int \int_{0 \leq t_1 \leq ... \leq t_{n_1 + n_2} \leq 1}\sum_{k_1 = 1}^{\infty} (x_1x_2t_1)^{k_1}\frac{dt_1}{t_1}  ...   \frac{dt_{n_1}}{t_{n_1}}\cdot
\sum_{k_2 = 1}^{\infty} (x_1t_{n_1+1})^{k_2}\frac{dt_{n_1 +1}}{t_{n_1 +1}}  ...   \frac{dt_{n_1 +n_2}}{t_{n_1 + n_2}} \quad =
$$
$$
\sum_{k_1,k_2 \geq 1}^{\infty} \frac{x_1^{k_1}x_2^{ k_1 + k_2}}
{ k_1^{n_1}(k_1 + k_2)^{n_2}}
$$
The general case is completely similar. 



{\it Generating series for multiple polylogarithms}. Set
$$
{\rm Li} (x_{1},...,x_{m}|t_1,...,t_m ) := \sum_{n_i \geq 1}{\rm Li}_{n_{1},...,n_{m}}(x_{1},...,x_{m} )t_1^{n_1-1}... t_m^{n_m-1}
$$
 \begin{lemma}  \label{2.new}  
One has the following identity of formal power series
\begin{equation} \label{5.new}
{\rm Li} (x_{1},...,x_{m}|t_1,...,t_m ) 
=   
\sum_{0 < k_{1} < k_{2} < ... < k_{m} } \frac{x_{1}^{k_{1}}x_{2}^{k_{2}}
... x_{m}^{k_{m}}}{(k_{1} -t_1)(k_{2} -t_2) ... (k_{m} -t_m)}
\end{equation}
\end{lemma}

{\bf Proof}. Follows immediately from the identity
$$
\frac{1}{k -t} = \frac{1}{k(1 -t/k)} = \sum_{n=1}^{\infty}\frac{t^{n-1}}{k^n}
$$

{\bf 4. Monodromy of multiple logarithms}. We will deduce it from the following 
general statement. Let $X$ be a complex variety, $f_1, ..., f_m$ 
rational functions on $X$ and $D$ an irreducible divisor in $X$. We will assume that 

 {\it $D$ is not a component of 
both ${\rm div}f_i$ and ${\rm div}f_{i+1}$ for any $1 \leq i\leq m-1 $} 

We will study the monodromy of the iterated integral 
$\int_{\gamma } 
d\log f_1  \circ  ...\circ   d\log f_m$ when the path $\gamma$ with fixed endpoints $a,b$ 
crosses the generic point $p$ of the divisor $D$. 
Let $\delta$ be a little path around the divisor $D$ near $p$. 
Define an index $v_{\delta}(f) \in \Z$ by 
$$
v_{\delta}(f):= \frac{1}{2\pi i} \int_{\delta} d\log f 
$$

Let $\alpha_1$ be a path from $a$ to a point $q \in \delta$, and $\alpha_2$ a 
path from $q$ to $b$. 
Let $\gamma:= \alpha_2\alpha_1$, i.e. we go first by 
$\alpha_1$, then by $\alpha_2$, and
 $\gamma':= \alpha_2\delta\alpha_1$. 
Denote by $\overline \alpha_i$ the limiting path obtained from $\alpha_i$ when the cycle $\delta$ shrinks to the point $p$ on the divisor $D$. 

\begin{proposition} \label{hov2} Assuming that $D$ and the functions $f_i$ satisfy the condition above, one has 
$$
\frac{1}{2\pi i}\Bigl(\int_{\gamma'} - 
\int_{\gamma }\Bigr) 
(d\log f_1  \circ  ...\circ   d\log f_m) = 
$$
\begin{equation} \label{UPssala}
\sum_{k=1}^mv_{\delta}(f_k)\Bigl(\int_{\overline \alpha_1} 
d\log f_1  \circ   ...\circ   d\log f_{k-1}\Bigr) \cdot 
\Bigl(\int_{\overline \alpha_2} 
d\log f_{k+1}  \circ   ...\circ   d\log f_m\Bigr)  
\end{equation}
\end{proposition}

{\bf Remark}. The condition on  $f_1, ... ,  f_m$ implies that the  integrals 
$$
\int_{\overline \alpha_1} 
d\log f_1  \circ    ...\circ   d\log f_{k-1}\quad \mbox{and} \quad 
\int_{\overline \alpha_2} 
d\log f_{k+1}   \circ   ...\circ   d\log f_{m} 
$$
are convergent provided that $v_{\delta}(f_k) \not = 0$: indeed, then neither $f_{k-1}$ nor $f_{k+1}$  can have $D$ as a component of their divisor.  
Therefore the right hand side of the formula makes sense. 

{\bf Proof}. 
The key tool is  the formula
\begin{equation} \label{BT}
\int_{\alpha\beta  }\omega_1 \circ ... \circ \omega_k  = \sum_{i=0}^k 
(\int_{\alpha }\omega_1 \circ ... \circ \omega_i) \cdot 
(\int_{\beta}\omega_{i+1} \circ ... \circ \omega_k)  
\end{equation}

Let $\delta(\varepsilon)$ be a loop around the divisor $D$ which shrinks to  $p \in D$ as $\varepsilon \lra 0$. 

\begin{lemma} \label{hov3} Assume that at least one of the functions $g_1, ..., g_k$ is regular at the generic point of the divisor $D$. Then 
$$
\lim_{\varepsilon \to 0}\int_{\delta(\varepsilon)} 
d\log g_{1}  \circ   ...\circ   d\log g_k  = 0
$$
\end{lemma}

{\bf Proof}. Taking a one dimensional slice transversal to $D$ at the point $p$  we reduce the problem to the case when $X$ is isomorphic to a neighborhood of zero in $\C$ and $\delta(\varepsilon)$ is radius $\varepsilon$ a loop around zero  starting at the point $\varepsilon$. Then
$$
\int_{\delta(\varepsilon)} 
d\log g_{1}  \circ   ...\circ   d\log g_k = 
\int_{0\leq \varphi_1 \leq ... \leq  \varphi_k \leq 2\pi} 
d\log g_{1}(\varepsilon e^{i\varphi_1}) \wedge  ...\wedge   
d\log g_k(\varepsilon e^{i\varphi_k})
$$
Write $g_j(t) = t^{n_j} \widetilde g_j(t)$ 
where $\widetilde g_j(0) \not = 0$. Then 
$$
\vert \int_{\varphi_{j-1}}^{\varphi_{j+1}} d \log \widetilde g_j(\varepsilon 
e^{i\varphi_j})\vert \leq C\cdot \varepsilon 
$$
On the other hand 
\begin{equation} \label{TALM}
\int_{\delta(\varepsilon)} \underbrace{d \log t \circ ... \circ d \log t}_{p\quad \mbox{times}} = \frac{(2 \pi i)^p}{p!}
\end{equation}
Since $g_j(t) = \widetilde g_j(t)$ for some $1 \leq j \leq k$, we are done. The proposition follows immediately from this lemma and 
formula (\ref{BT}).

\begin{corollary} \label{MCH}
Let $\gamma'$ be a deformation of the path $\gamma$ with fixed endpoints 
in $0$ and $1$  which crosses 
only once a singular point, $a_i$, $ 1 \leq i \leq m$, as shown on the 
picture. 

\begin{center}
\hspace{4.0cm}
\epsffile{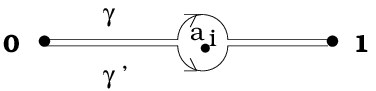}
\end{center}

Then
$$
{\rm I}_{\gamma'}(0; a_1, ..., a_m; 1)  - 
{\rm I}_{\gamma}(0; a_1, ..., a_m; 1)  =  
$$
\begin{equation} \label{ur3}
2 \pi i \cdot 
{\rm I}(0; a_1, ..., a_{i-1}; a_i) \cdot {\rm I}(a_i; a_{i+1}, ..., a_{m}; 1) 
\end{equation} 
\end{corollary}

{\bf Proof}. Follows immediately from the proposition \ref{hov2}. 

Since $\pi_1(\C - \{\infty, a_1, ..., a_m\})$ is generated by loops around the points 
$a_1, ..., a_m$ 
this corollary establishes the monodromy properties of multiple logarithms. The monodromy of 
multiple polylogarithms will be studied  in chapter 5.

{\bf 5.  The first shuffle relations: the product relations for power series (\ref{5})}. 
The product of two multiple polylogarithms is equal to a sum of multiple polylogarithms. 
More precisely, one has 
\begin{equation} \label{11.21.0.1}
{\rm Li}(x_1,...,x_m| t_1,...,t_m)\cdot {\rm Li}(x_{m+1},...,x_{m+n}|t_{m+1},...,t_{m+n}) \quad =
\end{equation}
$$
\sum_{\sigma \in \Sigma_{m,n}}{\rm Li}(x_{\sigma (1)},...,x_{\sigma (m+n) }|t_{\sigma ( 1) }, 
...,t_{\sigma (m+n)}) \quad + \quad \sum^{m+n-1}_{p= 1}\Bigl(\mbox{the depth $=m+n-p$ terms}\Bigr)
$$
where $\Sigma_{m,n}$ is the set of 
all shuffles $\sigma$ of $\{1,...,m\}$ and $\{m+1,...,m+n\}$, i.e. 
permutations $\sigma $ such that $\sigma(1) < ... < \sigma(m)$ 
and $\sigma(m+1) < ... < \sigma(m+n)$.

For example
$$
{\rm Li}_{m}(x) \cdot {\rm Li}_{n}(y) = {\rm Li}_{m,n}(x,y) + {\rm Li}_{n,m}(y,x) + {\rm Li}_{m+n}(xy)
$$
Indeed 
$$
{\rm Li}_{m}(x) \cdot {\rm Li}_{n}(y) = \sum_{0< k  , l  }^{\infty} \frac{x^{k}y^{l}}
{ k^{m}l^{n}} = 
\sum_{0< k < l  }^{\infty} \frac{x^{k}y^{l}}
{ k^{m}l^{n}} + \sum_{ k > l >0 }^{\infty} \frac{x^{k}y^{l}}
{ k^{m}l^{n}} + \sum_{ k  >0 }^{\infty} \frac{(xy)^{k} }
{ k^{m+n} }
$$
We can rewrite these formulas as
$$
{\rm Li}(x_1|t_1) \cdot {\rm Li}(x_2|t_2) =  {\rm Li}(x_1, x_2|t_1,  t_2) + {\rm Li}(x_2, x_1|t_2,  t_1) + 
$$
$$
\frac{1}{t_1-t_2}
\Bigl( {\rm Li}_1(x_1x_2|t_1 ) - {\rm Li}_1(x_1x_2|t_2 )\Bigr) 
$$
The proof can be also obtained from lemma \ref{2.new} using the identity
$$
\frac{1}{(k - t_1)(k - t_2)} = \frac{1}{t_1 - t_2}\Bigl(\frac{1}{k - t_1} - \frac{1}{k - t_2}\Bigl)
$$

One more example:
$$
{\rm Li} (x_1,x_2|t_1,t_2) \cdot {\rm Li} ( x_3|t_3) = 
$$
$$
{\rm Li} (x_1, x_2, x_3|t_1, t_2, t_3) + 
{\rm Li} (x_1, x_3, x_2|t_1, t_3, t_2) + {\rm Li} (x_3, x_1, x_2, | t_3, t_1, t_2 ) + 
$$
$$
\frac{1}{t_2-t_3}\Bigl({\rm Li} (x_1, x_2x_3|t_1, t_2 ) - {\rm Li} (x_1, x_2x_3|t_1, t_3 )\Bigr) +
$$
$$
 \frac{1}{t_1-t_3}\Bigl({\rm Li} (x_1x_3, x_2|t_1 , t_2) - {\rm Li} (x_1x_3, x_2|t_3, t_2)\Bigr) 
 $$
which is equivalent to 
$$
 {\rm Li}_{n_1,n_2}(x_1,x_2) \cdot {\rm Li}_{n_3}(x_3) =
$$
$$
 {\rm Li}_{n_1,n_2,n_3 }(x_1,x_2,x_3) + {\rm Li}_{n_1,n_3,n_2 }(x_1,x_3,x_2) + {\rm Li}_{n_3,n_1,n_2 }(x_3,x_1,x_2) + 
$$
$$
{\rm Li}_{n_1+n_3,n_2 }(x_1x_2, x_3) + {\rm Li}_{n_1,n_2+n_3 }(x_1, x_2x_3)
$$

The proof of the general statement is completely similar.
The reader can find in 
G.Racinet's theses [R] a beautiful interpretation of the 
first shuffle relations in terms of a 
group like element of a certain Hopf algebra. 

{\bf 6. The inversion formula}. 
The series
$$
B(\varphi_1, ..., \varphi_m | t_1, ..., t_m):= \sum_{-\infty <k_1< ...< k_m < \infty}
\frac{e^{2\pi i (\varphi_1 k_1 + ... + \varphi_m k_m)}}{(k_1-t_1) \cdot ... \cdot (k_m-t_m)}
$$
define a distribution on the $m$-dimensional torus with coordinates 
$0 \leq \varphi_i \leq 1$  depending on parameters $t_i$. For $m=1$ there is 
 a classical identity which goes back to Kronecker ([We], ch. 7.4) 
and provides the generating series for the Bernoulli polynomials:
$$
B(\varphi| t) \quad = \quad  -2\pi i
\frac{e^{2\pi i \{\varphi\}t}}{e^{2\pi i t}-1}\quad = \quad -2\pi i \cdot \sum_{n \geq 0} B_n(\{\varphi\}) 
\frac{(2 \pi i t)^{n-1}}{n!}
$$
Here $\{\varphi\}$ is the fractional part of $\varphi $. Indeed, 
$(\frac{1}{2\pi i} \frac{\partial}{\partial \varphi} - t)B(\varphi|t) = \delta(\varphi)$, 
and the distribution in the middle satisfies the same equation, 
so the difference between these two 
distribution  is a constant. Since clearly $B(0|t) = -\frac{1}{t}$ this constant is zero. 

Consider  the generating series 
$$
{\cal L}i(\varphi_{1}, ..., \varphi_m| t_{1}, ..., t_{m}):= \quad 
{\rm Li}(e^{2\pi i \varphi_{1}}, ..., e^{2\pi i \varphi_{m}}| t_{1}, ..., t_{m})
$$

Using the decomposition 
$$
\{-\infty <k_1< ...< k_m < \infty\} = 
$$
$$
\cup_{j}\Bigl(  \{k_1< ...< k_j < 0 < k_{j+1}... < k_m\} 
\cup  \{k_1< ...< k_{j-1} < k_j = 0 < k_{j+1}... < k_m\}\Bigr) 
$$ in the sum  (\ref{5.new}) we immediately get
\begin{equation} \label{12.2.00.4}
B(\varphi_1, ... ,\varphi_m | t_1, ..., t_m) = 
\end{equation}
$$
\sum_{j=1}^m (-1)^j {\cal L}i(-\varphi_j, ..., -\varphi_1 | -t_j, ..., -t_1) 
{\cal L}i(\varphi_{j+1}, ..., \varphi_m | t_{j+1}, ..., t_m) + 
$$
$$
\sum_{j=1}^m \frac{(-1)^{j}}{t_j}
{\cal L}i(-\varphi_{j-1}, ..., -\varphi_1| -t_{j-1}, ...,
- t_{1}) \cdot 
{\cal L}i(\varphi_{j+1}, ..., \varphi_m | t_{j+1}, ..., t_m) =
$$
Therefore the following proposition provides a formula expressing 
$$
{\cal L}i(\varphi_{1}, ..., \varphi_m| t_{1}, ..., t_{m}) + (-1)^m 
{\cal L}i(-\varphi_{m}, ..., -\varphi_1| -t_{m}, ..., -t_{1})
$$
via the products of lower depth multiple polylogarithms and Bernoulli polynomials. 
For example
$$
{\cal L}i(\varphi_{1}| t_{1}) -
{\cal L}i(-\varphi_1| - t_{1})  \quad = \quad -2\pi i 
\cdot \sum_{n \geq 1} B_n(\{\varphi\})\frac{(2 \pi i t)^{n-1}}{n!}
$$
and for $m=2$ one has 
$$
{\cal L}i(-\varphi_{2}, -\varphi_{1}| -t_{2}, -t_{1}) \quad  + \quad {\cal L}i(-\varphi_{2}, -\varphi_{1}| -t_{2}, -t_{1})  
$$
$$- \quad 
{\cal L}i(-\varphi_{2}| -t_{2}){\cal L}i(\varphi_{1}| t_{1}) \quad - \quad \frac{{\cal L}i(\varphi_{2}| t_{2})}{t_{1}} \quad + \quad 
\frac{{\cal L}i(-\varphi_{1}| t_{1})}{t_{2}} \quad =
$$
$$
{\cal L}i(\varphi_{2}| t_{2}-t_{1}) B(\varphi_{1}+ \varphi_{2}|t_1 )\quad  - \quad 
{\cal L}i(-\varphi_{1}| t_{1}-t_{2}) B(\varphi_{1}+ \varphi_{2}|t_2)
$$

\begin{proposition} \label{12.1.00.2} 
One has
$$
B(\varphi_1, ..., \varphi_m | t_1, ..., t_m) \quad = \quad \sum_{j=1}^m (-1)^{j-1}
{\cal L}i(-\varphi_{j-1}, ..., -\varphi_1| t_{j}-t_{j-1}, ...,
t_{j} - t_{1}) \cdot 
$$
\begin{equation} \label{12.1.00.4}
B(\varphi_1 + ... + \varphi_m | t_j) \cdot 
{\cal L}i(\varphi_{j+1}, ...,  \varphi_m| t_{j+1}- t_{j}, ...,
t_{m} -t_{j}) 
\end{equation}
\end{proposition}

{\bf Proof}. 
We need the following lemma

\begin{lemma} \label{12.1.00.2} Let $k_{i,j}:= k_i - k_j$, $t_{i,j}:= t_i - t_j$. 
Then 
\begin{equation} \label{12.1.00.1}
\frac{1}{(k_1-t_1) ...
 (k_m-t_m) } \quad = \quad \sum_{j=1}^m \frac{1}{(k_j-t_j) \prod_{i \not = j}
(k_{i,j}-t_{i,j}) }
\end{equation}
\end{lemma}

{\bf Proof}. 
We use the induction on $m$ to prove this identity simultaneously with the one
\begin{equation} \label{12.1.00.3}
\sum_{j=1}^m \frac{1}{\prod_{i \not = j}(k_{i,j}-t_{i,j}) } = 0
\end{equation}
For $m=2$ both identities are trivial. Suppose we proved (\ref{12.1.00.1})
 for $m=n$. Then it implies (\ref{12.1.00.3}) for $m=n+1$. Indeed, since 
$k_{j,1} - k_{s,1} = k_{j,s}$ and similarly for $t$'s, using  
(\ref{12.1.00.1}) we write the term in (\ref{12.1.00.3})
 corresponding to $j=1$ as:
$$
\frac{1}{\prod_{i = 2}^{n+1}(k_{i,1}-t_{i,1}) } = \sum_{j \not = 1} 
\frac{1}{(k_{j,1} - t_{j,1})\prod_{s \not = j}(k_{s,j}-t_{s,j}) }= 
-\sum_{j=2}^{n+1} \frac{1}{\prod_{i \not = j}(k_{i,j}-t_{i,j}) } 
$$
which is just equivalent to (\ref{12.1.00.3}). (The product in the denominator 
of the last sum is over $i \in \{2, 3, ..., \widehat j, ...,  n+1\}$). To prove 
(\ref{12.1.00.1}) for $m=n+1$ we write it, using the induction assumption, as
$$
\Bigl(\sum_{j=1}^n \frac{1}{(k_j-t_j) \prod_{i \not = j}(k_{i,j}-t_{i,j}) }\Bigr)\frac{1}{k_{n+1}-t_{n+1}}
$$
and then use 
$$
\frac{1}{(k_j-t_j)(k_{n+1}-t_{n+1})} \quad = \quad \frac{1}{k_{ n+1, j} - t_{n+1, j}}\cdot 
\Bigl( \frac{1}{k_j-t_j} - \frac{1}{k_{n+1}-t_{n+1}}\Bigr)
$$
and the identity (\ref{12.1.00.3}). The lemma is proved. 

To deduce the proposition from the lemma we write 
$$
B(\varphi_1, ... , \varphi_m| t_1, ..., 
t_m) \quad = \quad \sum_{0< k_1 < ... < k_m} \sum_{j=1}^{m}\frac{e^{2\pi i 
(\varphi_1 k_1 + ... + \varphi_m k_m)}}{(k_{j} - t_{j})\prod_{i \not = j}(k_{i,j} - t_{i,j})}
$$
and use (\ref{5.new}) and the formula
$$
\varphi_1 \cdot k_1 + ... + \varphi_m \cdot k_m = \sum_{i=1}^{j-1} -\varphi_i \cdot (k_j-k_i)
+  (\varphi_1 + ... + \varphi_m) \cdot k_j + \sum_{i=j+1}^m \varphi_i \cdot (k_i-k_j)
$$
to evaluate the sum related to the  $j$-th term of this formula. 
The proposition is proved.

{\bf 7. The second shuffle relations}. Let us set
\begin{equation} \label{shuffle2}
I^![a_1,...,a_m|t_1,...,t_m]:= I[a_1,...,a_m|t_1, t_1+t_2,t_1+t_2+t_3,...,t_1+...+t_m]
\end{equation}

\begin{theorem} \label{shuffle}
\begin{equation} \label{shuffle1}
I^![a_1,...,a_{k}|t_1,...,t_{k}]\cdot
I^![a_{k+1},...,a_{k+l}|t_{k+1},...,t_{k +l}] =
\end{equation}
$$
\sum_{\sigma \in \Sigma_{k,l}}I^![a_{\sigma (1)},...,a_{\sigma (k+l) }|
t_{\sigma ( 1) },...,t_{\sigma (k+l)}]
$$
\end{theorem}

We will give three different proofs of this theorem. Each of them sheds a new 
light on the result. 
 One needs the following preliminary result. 
Suppose we are given certain letters
$$
a_i, y_1^{(i)}, ... , y_{n_i-1}^{(i)} \qquad 1 \leq i \leq L
$$

Let us say that a string formed by these letters is
admissible if the following conditions are satisfied:
$$
a_p \quad \mbox{is before} \quad a_q \quad <=> \quad p<q
$$
$$
x^{(i)}_p \quad \mbox{is before} \quad x^{(i)}_q \quad <=> \quad p<q
$$

{\bf Combinatorial lemma}. {\it The number of admissible strings
 is equal to the coefficient with which $t_1^{n_1-1}
 ... t_L^{n_L-1}$ appears in the sum}
$$
\sum_{k_i \geq 0} t_1^{k_1}(t_1+t_2)^{k_2-1}...(t_1+ ... + t_L)^{k_L-1} 
$$

{\bf Proof}. The bijective correspondence between the monomials 
$t_1^{n_1-1} ... t_L^{n_L-1}$ appearing in the sum and admissible
strings is clear from the following pattern. It is convenient to treat for a
moment $t_i$ as non commuting variables. Then the monomial
$$
t_1^2\cdot t_2 t_1 t_1 \cdot t_3 t_2 t_2 t_1 t_1 \quad \mbox{of}
\quad t_1^2(t_1+t_2)^3(t_1+t_2+t_3)^5
$$
corresponds to the string
$$
a_1 y_1^{(1)} y_2^{(1)}\cdot a_2 y_1^{(2)} y_3^{(1)} y_4^{(1)} \cdot a_3
        y_1^{(3)} y_2^{(2)} y_3^{(2)} y_5^{(1)}  y_6^{(1)}
$$
So the variables $a_i$ separate the parts of a monomial in $t$'s corresponding to different factors 
$(t_1+...+t_i)^{k_i-1}$; to produce a string corresponding to a monomial in $t$'s we assume that 
in the letter $y_j^{(i)}$ the lower index $j$ shows how many times $t_i$ already appeared in the monomial. 

{\bf First proof of the theorem}: {\it via iterated integrals}. 
Recall the  well known  shuffle product formula for the product of iterated
integrals  
\begin{equation} \label{schprfo}
\int_{\gamma} \omega_{1} \circ ... \circ \omega_{k_1} \cdot \int_{\gamma}
\omega_{k_1+1} \circ ... \circ \omega_{k_l+k_2} = \sum_{\sigma \in \Sigma_{k_1,k_2}} \int_{\gamma} \omega_{\sigma(i_1)} \circ ... \circ \omega_{\sigma(i_{k_1 + k_2})}
\end{equation}

The  theorem follows immediately from this and the
combinatorial 
lemma.

For example to find the product of iterated integrals
$$
 \int_{0}^{1} \underbrace {\frac{dt}{t-a_{1}} \circ \frac{dt}{t} \circ ...
\circ \frac{dt}{t}}_{k_1 \quad  \mbox {times}} \quad \cdot \quad \int_{0}^{1} \underbrace {\frac{dt}{t-a_{2}} \circ \frac{dt}{t} \circ ...
\circ \frac{dt}{t}}_{k_2 \quad  \mbox {times}}
$$
we need to apply the combinatorial lemma twice: for the shuffles of the differentials where $dt/(t-a_1)$ goes before $dt/(t-a_2)$, and for the shuffles where it goes after.

{\bf Second proof of the theorem}:  {\it via the iterated
  integral presentation  
for the generating function}.

\begin{lemma} \label{iterfpse}
\begin{equation} \label{43212}
I^![x_1,...,x_m|t_1,...,t_m] = \int_0^1\frac{ s^{-t_1}}{x_1 -
  s}ds \circ ... \circ \frac{s^{-t_m}}{x_m - s}ds =
\end{equation}
$$
\int_{0 < s_1 < ... < s_m < 1}\frac{s_1^{-t_1} \cdot ... \cdot
  s_m^{-t_m}}{(x_1 - s_1) \cdot ... \cdot (x_m - s_m)}ds_1 ... ds_m
$$
\end{lemma}

The theorem follows immediately  from lemma \ref{iterfpse} and
the product formula for the iterated integrals.

{\bf Remark}. Compare formulas (\ref{43212}) and (\ref{5.new}). 

{\bf First proof of the lemma: {\it via iterated integrals}}. We need the following identity:
\begin{equation} \label{ur8-4}
\int_a^b \underbrace{\frac{dt}{t} \circ \frac{dt}{t}\circ  ... \circ \frac{dt}{t}}_{\mbox{$n$ times}} \quad = \quad \log^n\frac{b}{a} 
\end{equation}
One can
write the statement of the lemma  as follows
$$
I^![x_1,...,x_m|t_1,...,t_m] \stackrel{?}{= }
$$
$$
\sum_{n_1,...,n_m
  \geq 0}\frac{t_1^{n_1} \cdot ... \cdot t_1^{n_m}}{(n_1)!
  \cdot ... \cdot (n_m)!} \cdot \int ... \int_{0 < s_1 < ... <
  s_m < 1}\frac{(-\log 
  s_1)^{n_1} \cdot ... \cdot (-\log s_m)^{n_m} }{(x_1 -
  s_1) \cdot ... \cdot (x_m - s_m)}ds_1 ... ds_m
$$

Rewriting identity (\ref{ur8-4}) in the case $b=1, a=s$ as 
$$
\frac{(-\log s_i)^{n_i}}{n_i!} = \int_{\Delta}\frac{d y^{(i)}_1}{y^{(i)}_1} 
  ...  \frac{d y^{(i)}_{n_i}}{y^{(i)}_{n_i}} \qquad \Delta:= \{0 < s_i < y^{(i)}_1
  ... < y^{(i)}_{n_i} < 1\}
$$ 
and using the combinatorial lemma we get the needed formula.

{\bf Second proof of the lemma: {\it  via power series
  expansions}}. According to theorem \ref{2.3} 
formula (\ref{43212}) is equivalent to  the following formula
$$
{\rm Li}^![x_1,...,x_m|t_1,...,t_m] = \int_0^1\frac{(x_1
  ...  x_m)  s^{-t_1}}{1-(x_1 
  ...  x_m)  s}ds \circ ... \circ \frac{x_m 
  s^{-t_m}}{1-x_m  s}ds
$$
One has
$$
\int_0^{u}\frac{x s^{-t}}{1-xs}ds = \sum_{k=1}^{\infty} x^k
\int_0^u s^{k-t-1}ds = \sum_{k=1}^{\infty} \frac{x^k}{k-t}u^{k-t}
$$
Therefore
$$
\int_0^1 \Bigl(\int_0^{s_2}\frac{x_1x_2 s_1^{-t_1}}{1-x_1x_2s_1}ds_1\Bigr)
\frac{x_2 s_2^{-t_2}}{1-x_2s_2}ds_2 = 
$$
$$
\int_0^1
\sum_{k_1=1}^{\infty}
\frac{(x_1x_2)^{k_1}}{k_1
  -t_1}\sum_{k_2=1}^{\infty}x_2^{k_2}\int_0^1s_2^{k_1 + k_2 - t_1
  - t_2 -1} ds_2 \quad = \quad 
\sum_{k_1,k_2=1}^{\infty}
\frac{x_1^{k_1}x_2^{k_1+k_2}}{(k_1-t_1)(k_1+k_2 - t_1 -t_2)}
$$
and a similar proof in general. 
The lemma is proved.

More generally, let 
$$
{\rm I}_{n_1, ..., n_m}(a_0; a_1, ..., a_m; a_{m+1}):= 
$$
\begin{equation} \label{ur4-8,1}
\int_{a_0}^{a_{m+1}} 
\frac{dt}{t-a_1}\underbrace{\frac{dt}{t} \circ \frac{dt}{t}\circ  ... \circ \frac{dt}{t}}_{\mbox{$n_1-1$ times}} \circ  ... \circ \frac{dt}{t-a_m}\underbrace{\frac{dt}{t} \circ \frac{dt}{t}\circ  ... \circ \frac{dt}{t}}_{\mbox{$n_m-1$ times}} 
\end{equation}
Consider the generating series
$$
{\rm I}(a_0; a_1, ..., a_m; a_{m+1}| t_1; ...; t_m) := \sum_{n_i \geq 0}{\rm I}_{n_1, ..., n_m}(a_0; a_1, ..., a_m; a_{m+1})t_1^{n_1-1} ... t_m^{n_m-1}
$$
$$
I^!(a_0; a_1, ..., a_m; a_{m+1}| t_1; ...; t_m) := {\rm I}(a_0; a_1, ..., a_m; a_{m+1}| t_1;
t_1 + t_2;  ...; t_1+...+t_m)
$$
Then we have the following result similar to the lemma above.
\begin{lemma} \label{iterfpse*}
\begin{equation} \label{ur4-8,2}
I^!(a_0; a_1, ..., a_m; a_{m+1}| t_1; ...; t_m) = a_{m+1}^{t_1+ ... + t_m}\cdot 
\int_{a_0}^{a_{m+1}}\frac{ s^{-t_1}}{x_1 -
  s}ds \circ ... \circ \frac{s^{-t_m}}{x_m - s}ds
\end{equation}
\end{lemma}

{\bf Proof}.  One can proceed similarly to the first proof of lemma \ref{iterfpse}. 

{\bf   Third proof of the theorem}:  {\it via power series expansions}.  We 
will use generating functions given in lemma \ref{2.new}.  In the simplest case $ k=l=1$  we need to prove that
$$
{\rm Li}(x_1|t_1) \cdot {\rm Li}(x_2|t_2) = {\rm Li}(\frac{x_1}{x_2},x_2|t_1,t_1+t_2) + {\rm Li}(\frac{x_2}{x_1},x_2|t_2,t_1+t_2) 
$$
The right hand side is equal to
$$
\sum_{m_1,m_2 >0}\frac{(x_1/x_2)^{m_1}x_2^{m_1 + m_2}}{(m_1 -t_1)(m_1 + m_2 - t_1 -t_2)} + 
\sum_{m_2,m_1 >0}\frac{(x_2/x_1)^{m_2}x_1^{m_1 + m_2}}{(m_2 -t_2)(m_1 + m_2 - t_1 -t_2)} =
$$
$$
\sum_{m_1,m_2 >0}x_1^{m_1}x_2^{m_2} \cdot \Bigl(\frac{1}{(m_1 -t_1)(m_1 + m_2 - t_1 -t_2)} + \frac{1}{(m_2 -t_2 )(m_1 + m_2 - t_1 -t_2)}\Bigr)
$$
and the statement follows from the identity
$$
\frac{1}{p_1p_2} = \frac{1}{p_1( p_1 + p_2)} + \frac{1}{p_2( p_1 + p_2) }
$$
In general the proof follows similarly from the following identity:

\begin{lemma}
\begin{equation} \label{A3}
\frac{1}{p_1( p_1 + p_2) ... ( p_1 + ... + p_k)} \cdot \frac{1}{p_{k+1}( p_{k+1}  + p_{k+2}) ... ( p_{k+1}   + ... + p_{k+l})} = 
\end{equation} 
$$
\sum_{\sigma \in \Sigma_{k,l}} \frac{1}{p_{\sigma(1)}( p_{\sigma(1)} + p_{\sigma(2)} ) ... ( p_{\sigma(1)} + ... + p_{\sigma(k+l)} )}
$$
\end{lemma}

{\bf Proof}. We will use induction. The shuffles $\Sigma_{k,l}$ are divided on two parts 
$\Sigma^I_{k,l}$ and $\Sigma^{II}_{k,l}$ as follows. A shuffle is in the first (resp. second) 
part if it is ended by $p_k$ (resp. $p_{k+l}$). By induction
$$
\sum_{\sigma \in \Sigma^I_{k,l}} \frac{1}{p_{\sigma(1)}( p_{\sigma(1)} + p_{\sigma(2)} ) ... ( p_{\sigma(1)} + ... + p_{\sigma(k+l)} )} = \frac{1}{p_1+ ... + p_{k+l}}\cdot 
$$
\begin{equation} \label{A1}
\frac{1}{p_1 (p_1+p_2) ... (p_1+...+p_{k-1}) \cdot p_{k+1} (p_{k+1} + p_{k+2}) ... (p_{k+1} + ... + p_{k+l}) }
\end{equation} 
A similar sum over $\Sigma^{II}_{k,l}$ is equal to 
\begin{equation} \label{A2}
\frac{1}{p_1  ... (p_1+...+p_{k}) \cdot p_{k+1}  ... (p_{k+1} + ... + p_{k+l-1}) \cdot (p_1+ ... + p_{k+l})}
\end{equation} 
using the identity 
$$
\frac{1}{(p_1+...+p_{k})(p_{k+1} + ... + p_{k+l})} = 
$$
$$
\frac{1}{(p_1+...+p_{k})(p_{1} + ... + p_{k+l})} + 
 \frac{1}{(p_{k+1} + ... + p_{k+l})(p_{1} + ... + p_{k+l})} 
$$
we see that the sum of (\ref{A1}) and (\ref{A2}) equals  to  (\ref{A3}). The lemma is proved.

{\bf Remark}.  Let us rewrite the product relation (\ref{shuffle1}) in terms of the generating series for 
the iterated integrals. Set
$$
{\rm I}_!(a_1,...,a_m|t_1,...,t_m):= {\rm I}(a_1\cdot  ...\cdot a_m,a_2\cdot ...\cdot a_m,...,a_m|t_1, ..., t_m)
$$
Then 
$$
{\rm I}_!(a_1,...,a_k|t_1,...,t_k)\cdot {\rm I}_!(a_{k+1},...,a_{k+l}|t_{k+1},...,t_{k+l}) =
 $$
$$
\sum_{\sigma \Sigma_{k,l}}{\rm I}_! (a_{\sigma (1)},...,a_{\sigma (k+l) }|t_{\sigma ( 1) },...,t_{\sigma (k+l)}) + \mbox{lower depth terms}
$$
This makes the first shuffle relations surprisingly similar (I would say dual) to the second 
shuffle relations. I do not know what stays behind this duality. 
 
{\bf 8. Distribution relations}. Recall that $w:= n_1 + ... + n_m$. 
\begin{proposition} \label{2.5}
   Let us suppose that  $|x_i| \leq 1$. Then the following identities hold 
for $w\geq 2$ or $w=1$ and $x_1^l \not = 1$:
\begin{equation} \label{5n}
{\rm Li}_{n_{1},...,n_{m}}(x_{1}^l,...,x_{m}^l)   = l^{w-m} \sum_{0 \leq a_i \leq l-1} {\rm Li}_{n_{1},...,n_{m}}(\zeta^{a_1}x_{1},...,\zeta^{a_m}x_{m})   
\end{equation}
If $x_i <1$ one can rewrite these identities using the generating series:
\begin{equation} \label{5n.new}
 {\rm Li} (x^l_{1},...,x^l_{m}|t_{1},...,t_{m})   = 
\sum_{ y^l_i  =  x_i^l} {\rm Li}( y_{1},..., y_{m}|lt_{1},...,lt_{m})   
\end{equation}
\end{proposition}

{\bf Proof}. The conditions $w\geq 2$ or $w=1$ and $x_1^l \not = 1$ are the necessary and sufficient 
conditions for convergence of the series on the left hand side of (\ref{5n}). 
Notice that $w=1$ implies that $m=1$. 
Let $\zeta$ be a primitive $l$-th root of unity:
$\zeta^l =1$. 
The right hand side of (\ref{5n}) is equal to 
$$
\sum_{0 \leq a_i \leq l-1} \sum_{k_1 < ... < k_m} \frac{x_{1}^{k_{1}}
... x_{m}^{k_{m}}}{k_{1}^{n_{1}}...k_{m}^{n_{m}}}
\zeta^{a_1 k_1 + ... + a_m k_m}
$$
If $l \not|k_i $ for some $i$ then the sum over $a_i$ is zero. So the
right hand side  of (\ref{5n}) is equal to
$$
l^{w-m} \cdot l^m \sum_{k_1 < ... < k_m}\frac{x_{1}^{lk_{1}}
... x_{m}^{lk_{m}}}{(lk_{1})^{n_{1}}...(lk_{m})
^{n_{m}}} = 
{\rm Li}_{n_{1},...,n_{m}}(x_{1}^l,...,x_{m}^l) 
$$

{\bf 9. Canonical regularization}. (Compare with the end of s. 5.3 below). 
Consider the integral 
\begin{equation} \label{2**}
\int_{0}^{1} \underbrace{\frac{dt}{t} \circ ... \circ \frac{dt}{t}}_{p} 
\circ
\frac{dt}{t-a_1} \circ ... \circ 
\frac{dt}{t-a_m} \circ 
\underbrace{\frac{dt}{t-1} \circ ... \circ \frac{dt}{t-1}}_{q}
\end{equation}
where $a_1 \not =  0, a_m \not = 1$.
It is divergent if $p$ or $q$ is different from zero.
To  regularize it we will  prove that the integral
\begin{equation} \label{2*12321}
\int_{\varepsilon}^{1-\varepsilon} \underbrace{\frac{dt}{t} \circ ... \circ \frac{dt}{t}}_{p} 
\circ
\frac{dt}{t-a_1} \circ ... \circ 
\frac{dt}{t-a_m} \circ 
\underbrace{\frac{dt}{t-1} \circ ... \circ \frac{dt}{t-1}}_{q}
\end{equation}
is a polynomial in $\log\varepsilon$ and take the constant term as the regularized value.   
For example the regularization of $\log(1-x)$ at $x=1$ is zero.
Notice that  $\varepsilon$
is the canonical coordinate on ${{\Bbb P}}^1\backslash \{0,1,\infty\}$.

\begin{proposition} \label{p2.5} a) Integral (\ref{2*12321}) equals to 
$$
(-1)^{q}\int_{0}^{1} \frac{\log^p(t/\varepsilon)}{p!} \frac{dt}{t-a_1} \circ...\circ  \frac{dt}{t-a_m}  \frac{\log^q((1-t)/\varepsilon)}{q!}
$$
It is a polynomial in $\log\varepsilon$ whose coefficients are multiple polylogarithms with special points at $0,1, a_1,...,a_m$.

  b) The regularized value of integral (\ref{2*12321}) 
is obtained as follows: take the integral
$$
\int_0^1\frac{dt}{t-a_1}\circ ... \circ \frac{dt}{t-a_m} 
$$
and insert to it in all possible ways ordered $q$ (resp. $p$) - tuple of 1-forms $\frac{dt}{t-1}$ (resp. $\frac{dt}{t}$) before (resp. after)
$\frac{dt}{t-a_m}$ (resp. $\frac{dt}{t-a_1}$). 
Then take the sum of all multiple polylogarithms obtained this  way and multiply it by $(-1)^{p+q}$.
\end{proposition}

One can reformulate the part b) as follows:

Let $\omega_i = f_i(t)dt$. Suppose that $f_{m}(t)$ (resp.
$f_{1}(t)$) does not have a pole at $t=1$ (resp. $t = 0$). Then the regularized value of the iterated integral
\begin{equation} \label{r1}
\int_0^1\frac{dt}{t}\circ...\circ\frac{dt}{t}\circ\omega_1\circ...\circ\omega_{m}\circ \frac{dt}{t-1}\circ...\circ\frac{dt}{t-1}
\end{equation}
is equal to 
\begin{equation} \label{r2}
(-1)^{p+q}\cdot \int_{\Delta}\frac{dy_1}{y_1}...\frac{dy_p}{y_p}f_1(t_1)dt_1 ...f_{m}(t_{m})dt_m
\frac{dx_1}{x_1-1}...\frac{dx_q}{x_q-1}
\end{equation}
where
$$
\Delta = \{t_1 < y_1 < ... < y_p < 1,   0 < t_1 < ... < t_{m} < 1; 0 < x_1 < ... < x_q < t_{m} \}
$$

{\bf Proof}. Integrating (\ref{2*12321}) $p$ times we get
$$
 \int_{\varepsilon}^{1-\varepsilon}
\frac{\log^p(t/\varepsilon)}{p!}\frac{dt}{t-a_1} \circ ... \circ \frac{dt}{t-a_m} \circ \frac{dt}{t-1}\circ ... \circ \frac{dt}{t-1}
$$
Further, integrating the last $q$ differentials in (\ref{2*12321})  we get 
we get
$$
(-1)^{q}\int_{\varepsilon}^{1-\varepsilon} \frac{\log^p(t/\varepsilon)}{p!} \frac{dt}{t-a_1}  \circ...\circ
\frac{dt}{t-a_m}\frac{\log^q((1-t)/\varepsilon)}{q!}
$$
This is a polynomial in $\log \varepsilon$ whose coefficients have limit when $\varepsilon \lra 0$. 
In particular the regularized value
of (\ref{r1}) is
\begin{equation} \label{r21}
(-1)^{q}\int_{0}^1\frac{\log^p(t)}{p!}\frac{dt}{t-a_1} \circ...\circ\frac{dt}{t-a_m} \frac{\log^q(1-t)}{q!}
\end{equation}
This coincides with (\ref{r2}). Notice that we are picking up the sign $(-1)^p$ when go from the integral 
in (\ref{r21}) to the integral in  (\ref{r2}). The proposition is proved.

Let us choose once forever a coordinate $t$ in $\C P^1 \backslash \{\infty\}$ which 
identifies it with $\C$.  Then there is a tangent vector $\partial_t $ such that 
$<dt, \partial_t > = 1$.
Choose a tangent vector $v_{a} = \lambda_a \partial_t$  at every point $a \in \C$.  
Let $a_0, ..., a_{m+1}$ be {\it arbitrary} configuration of points 
in $\C$. 
Then for any path $\gamma$ between the tangent vectors $v_{a_0}$ and 
 $v_{a_{m+1}}$ we have the regularized integral
\begin{equation} \label{12.21.00.1q}
{\rm I}_{\{\gamma, v_a\}}(a_0; a_1, ..., a_m; a_{m+1})
\end{equation} 
It is defined follows.  
The integral 
\begin{equation} \label{12.21.00.112q}
\int_{a_0 + \varepsilon \lambda_{a_0}}^{a_{m+1} + \varepsilon \lambda_{a_{m+1}}} 
\frac{dt}{t-a_1} \circ ... \circ \frac{dt}{t-a_m}
\end{equation} 
admits an asymptotic expansion as $\varepsilon \to 0$, 
which is a polynomial in $\log \varepsilon$ 
whose coefficients have limit when  $\varepsilon \to 0$. To prove we proceed just 
like in the proof of 
proposition \ref{p2.5}. Another, more conceptual,  
proof see in the end of s. 5.3 below. See also [D] about  the torsor 
of path between the tangential base points.
When  
$\varepsilon \to 0$ the limit value of the 
free term of this polynomial is, by definition, the regularized value 
 (\ref{12.21.00.1q}). The regularization  depends only  
on the choice of tangent vectors  $v_{a_0}$ and 
 $v_{a_{m+1}}$, but not on the choice of the coordinate $t$. 
If $a_0 \not = a_1$ and $a_{m} \not = a_{m+1}$ 
the integral is convergent and its value is 
the regularized value.

Let
\begin{equation} \label{GE1}
{\rm I}_{n_0; n_1, ..., n_m}(a_0; a_1, ..., a_m;a_{m+1}):=
\end{equation}
$$
\int_{a_0}^{a_{m+1}} \underbrace{\frac{dt}{t} \circ ... \circ \frac{dt}{t}}_{n_0} 
\circ
\frac{dt}{t-a_1} \circ \underbrace{\frac{dt}{t} \circ ... \circ \frac{dt}{t}}_{n_1-1} \circ ... \circ 
\frac{dt}{t-a_m} \circ 
\underbrace{\frac{dt}{t} \circ ... \circ \frac{dt}{t}}_{n_m-1}
$$
where the integral is regularized using the tangent vector $\partial_t$ at $a_0$ and 
$-\partial_t$ at $a_{m+1}$.
Let us define the generating series 
\begin{equation} \label{GE2}
{\rm I}(a_0; a_1, ..., a_m; a_{m+1}| t_0; t_1; ... ; t_m):= 
\end{equation}
$$
\sum_{n_0 \geq 0, n_i \geq 1} {\rm I}_{n_0; n_1, ..., n_m}(a_0; a_1, ..., a_m;a_{m+1})
t_0^{n_0} t_1^{n_1-1} ... t_m^{n_m-1}  
$$

{\bf Remark}. It would be more natural to use a notation 
$I_{n_0, n_1-1, ..., n_m-1}(...)$ for the iterated integral 
(\ref{GE1}), so the indices correspond to the numbers of $dt/t$ differentials 
between two consequitive differentials of type $dt/(t-a_i)$. 
However this would eventually contradict to the classical notation ${\rm Li}_n(z)$

When $a_0=0$ the coefficients of the generating series (\ref{GE2}) 
with $n_0 >0$ are divergent. The regularized value of  
this generating series is given by the following proposition.

\begin{proposition} \label{GENP}
$$
{\rm I}(0; a_1, ..., a_m; a_{m+1}| t_0; t_1; ... ; t_m) = 
$$
\begin{equation} \label{3.10.01.1}
e^{t_0 \cdot \log a_{m+1}  } \cdot {\rm I}(a_1: ...: a_m: a_{m+1}| t_1 - t_0, ... , t_m  - t_0) 
\end{equation}
\end{proposition}

{\bf Proof}. Let $\omega_i = f_i(x)dx$. Suppose that $f_{1}(x)$ 
does not have a pole at $x=0$, and $f_{k}(x)$ 
does not have a pole at $x=a_{m+1}$. Then similar to proposition \ref{p2.5} the regularized value of 
$$
\sum_{n\geq 0} t_0^n \cdot \int_0^{a_{m+1}} 
\underbrace{ \frac{dx}{x}\circ ... \circ \frac{dx}{x}}_{n} \omega_1 \circ ... \circ \omega_k
$$
is equal to
$$
\sum_{n\geq 0} t_0^n \cdot \int_0^{a_{m+1}} \frac{(\log x)^n}{n!}\omega_1 \circ ... \circ \omega_k = 
\int_0^{a_{m+1}} e^{t_0 \cdot \log x  } \cdot \omega_1 \circ ... \circ \omega_k = 
$$
\begin{equation}\label{GE3}
e^{t_0 \cdot \log a_{m+1}  } \cdot \int_0^{a_{m+1}} e^{t_0 \cdot \log x/a_{m+1}  } 
\omega_1 \circ ... \circ \omega_k 
\end{equation}
Using the formula
$$
\frac{(\log x/a_{m+1})^n}{n!} = \int_x^{a_{m+1}} \underbrace{ \frac{dt}{t}\circ ... \circ \frac{dt}{t}}_{n} 
$$
in order to write integral (\ref{GE3}) in the form similar to (\ref{r2}) we prove the proposition. 

{\bf Second proof}. Formula (\ref{3.10.01.1}) 
follows formally from the shuffle product formula for the regularized iterated 
integrals. To prove this we use the identity
$$
\Bigr(\int_{\varepsilon}^{a_{m+1}}\frac{dt}{t}\Bigl) \cdot 
\int_{\varepsilon}^{a_{m+1}} \underbrace { \frac{dt}{t} \circ ... \circ 
\frac{dt}{t}}_{k \quad  \mbox {times}} \circ 
\omega_1 \circ ... \omega_m \quad = \quad (k+1) \int_{\varepsilon}^{a_{m+1}} 
\underbrace { \frac{dt}{t} \circ ... \circ 
\frac{dt}{t}}_{k+1 \quad  \mbox {times}} \omega_1 \circ ... \omega_m +
$$
$$
\int_{\varepsilon}^{a_{m+1}} 
\underbrace { \frac{dt}{t} \circ ... \circ 
\frac{dt}{t}}_{k \quad  \mbox {times}} \circ 
\left\{ \frac{dt}{t} \circ_{Sh} (\omega_1 \circ ... \omega_m) - 
\frac{dt}{t} \circ \omega_1 \circ ... \omega_m \right\} 
$$
where 
$$
(\omega_1 \circ ... \circ\omega_ k) \circ_{Sh}(\omega_{k+1} \circ ... \omega_{k+l})
:= \sum_{\sigma \in \Sigma_{k,l}} \omega_{\sigma(1)} \circ ... \circ \omega_{\sigma(k+l)}
$$
and proceed by induction on $k$. We will use the second proof 
in the Hodge version of the story.

\begin{corollary} \label{incor}
The regularized value of a divergent multiple polylogarithm at $N$-th roots of unity is a linear combination with rational coefficients of convergent multiple polylogarithms at $N$-th roots of unity.
\end{corollary}

\begin{lemma} \label{inlem}
The second shuffle relations hold for the canonical regularization of multiple polylogarithms. 
\end{lemma}

{\bf Proof}. The proof of the second shuffle relations works without any changes if we replace $\int_0^1$ by $\int_{\varepsilon}^{1-\varepsilon}$.   So the second shuffle relations hold for the canonical 
asymptotic expansions in $\log \varepsilon$, and in particular for the constant terms. 

{\bf Remark}. In chapter 5 we  use the tangent vector $\partial_t$ 
instead of $-\partial_t$ at $a_{m+1}$.  
The difference is computed as follows.
If $$
I(\varepsilon)= I_0(\varepsilon)  + \log(\varepsilon)I_1(\varepsilon) + ... 
+ \log^p(\varepsilon)I_p(\varepsilon)
$$
 is the asymptotic expansion for the 
integral of type $\int_{x}^{a_{m+1} + \varepsilon }$ then 
the asymptotic expansion for $\int_{x}^{a_{m+1} - \varepsilon }$ 
is given  by $I(-\varepsilon)$. To compute the constant term of $I(-\varepsilon)$ 
observe that $I_k(\varepsilon)$ has a limit when $\varepsilon \to 0$, and 
 $\log (-\varepsilon ) = \pm \pi i +
\log (\varepsilon ) $ with the sign depending whether the path between the tangential 
base points  $v_{a_{m+1}}$ and $-v_{a_{m+1}}$ is going clockwise or counterclockwise.

{\bf 10. The power series regularization}.  Let us suppose that $|x_i| \leq 1$. 
The especially interesting for us case is when all the variables 
$x_i$ are $N$-th roots of unity. Then    
the power series ${\rm Li}_{n_1,...,n_m}(x_1,...,x_m)$ are divergent if and only if 
$n_m=1, x_m=1$.  
  
 In this subsection we present an approach to the   regularization 
of the first shuffle relations for multiple polylogarithms which has been 
developed in the  end of [G3]. 

D. Zagier was  first to discover (around 1993, unpublished) 
that one must regularize 
the power series shufle relations for the multiple zeta values 
in a different way 
and worked out several regularization procedures. 
Another approach was recently  developed Boutet de Monvil, see a detailed 
account in 
the theses of G. Racinet [R]. In particular it 
provides an explicit formula relating this regularization of multiple zeta values with 
the canonical regularization. 
I left to the reader to work out 
the relationship between the regularization used below and the one used 
in [R]. 

\begin{lemma} \label{olmn}
Let $|x_i | \leq 1 $. Then the power series
\begin{equation}  \label{ol1}
{\rm Li}_{n_1,...,n_l,1, ...,1}(x_1,...,x_l, 1- \varepsilon, ... , 1- \varepsilon)
\end{equation} 
have  an asymptotic expansion which is a polynomial in $  \log
\varepsilon $ whose coefficients are explicitly computable  $\Q$-linear combinations of lower depth multiple polylogarithms. This polynomial has weight $w$ provided that the weight of 
$  \log \varepsilon $ is $2$. 
\end{lemma}

 {\bf Proof}.   We will prove it by  induction on $m-l$. Using   
$$
(- \log \varepsilon)^m  = \Bigl(\sum_{k>0}\frac{(1- \varepsilon)^k}{k}\Bigr)^m  
$$  
together with the induction assumption we get      
$$
{\rm Li}_{ 1, ...,1}(  1- \varepsilon, ... , 1- \varepsilon) \quad = \quad
\frac{(- \log \varepsilon)^m }{m!} \quad + 
$$
\begin{equation} \label{req212}
\sum_{0<i<m}\mbox{(lower depth multiple $\zeta$'s)} (\log \varepsilon)^i
\end{equation} 
   
 Applying   the 
 power series product formula to
 $$ 
 {\rm Li}_{n_1,...,n_l }(x_1,...,x_l) \cdot  
  {\rm Li}_{ 1, ...,1}(  1- \varepsilon, ... , 1- \varepsilon) 
 $$
 and using the induction  on $m-l$, we get the lemma. 
Here is an example. 
$$
{\rm Li}_1(x_1) \cdot {\rm Li}_{ 1, 1}(  1- \varepsilon, 1- \varepsilon) = 
$$
$$
{\rm Li}_{ 1, 1, 1}(  x_1, 1- \varepsilon, 1- \varepsilon) + {\rm Li}_{ 1, 1, 1}(  1- \varepsilon, x_1, 1- \varepsilon) + {\rm Li}_{ 1, 1, 1}(  1- \varepsilon, 1- \varepsilon,  x_1) + 
$$
$$
{\rm Li}_{ 2, 1}(  x_1(1- \varepsilon), 1- \varepsilon) + {\rm Li}_{ 1,2}(1- \varepsilon,   x_1(1- \varepsilon)) 
$$
Since $x_1 \not = 1$ the asymptotic expansion of this expression coincides with the asymptotic expansion of 
$$
{\rm Li}_{ 1, 1, 1}(  x_1, 1- \varepsilon, 1- \varepsilon) + {\rm Li}_{ 1, 1, 1}(  1, x_1, 1- \varepsilon) + 
{\rm Li}_{ 1, 1, 1}(  1, 1,  x_1) + {\rm Li}_{ 2, 1}(  x_1, 1- \varepsilon) + {\rm Li}_{ 1,2}(1,   x_1) 
$$

Let us assume that $n_l \not = 1$ or $x_l \not = 1$. 
Denote by ${\widetilde {\rm {\rm Li}}}_{n_1,...,n_l,1, ...,1}(x_1,...,x_l, 1, ... , 1)$ 
the constant term of the expansion 
of ${\rm Li}_{n_1,...,n_l,1, ...,1}(x_1,...,x_l, 1- \varepsilon, ... , 1- \varepsilon)$ 
in $\log \varepsilon$ and call it 
the {\it power series regularization}. 

The asymptotic expansion of 
${\rm Li}_{n_1,...,n_l,1, ...,1}(x_1,...,x_l, 1 - \varepsilon, ... , 1- \varepsilon)$ is called
{\it the power series asymptotic expansion of}
${\rm Li}_{n_1,...,n_l,1, ...,1}(x_1,...,x_l, 1, ... , 1)$.

\begin{proposition} \label{olmnq}
The power series asymptotic expansions of multiple polylogarithms satisfy the first shuffle relations. 
In particular the power series regularization also satisfy the first shuffle relations.
\end{proposition}

{\bf Proof}. 
Say that in ${\rm Li}_{n_1,...,n_m}(x_1,...,x_m)$ the variable $x_i$ enters with the 
multiplicity  $n_i$. 
In ${\rm Li}_{n_1,...,n_l,1, ...,1}(x_1,...,x_l, 1- \varepsilon, ... , 1- \varepsilon)$ 
each of the variables $1- \varepsilon$ enters with the multitude  $1$, and also $x_l \not = 1$. 
We call such power series {\it admissible}. 

Taking the shuffle product of two such power series we usually get terms which are not admissible 
for one of the following two reasons: there are variables 
$(1- \varepsilon)^p$ with $ p>1$, or $x_i \cdot (1- \varepsilon)^q$, $q\geq 1$. 
However each of these terms can be replaced by an admissible term with the same 
asymptotic expansion. Indeed, 
the variables  $(1- \varepsilon)^p$ enter with the multitude  $p$, and if $p>1$ then we can put $\varepsilon = 0$ in this variables, and also in all variables to the left of this one, keeping the same asymptotic expansion. A similar statement is valid for terms $x_l(1- \varepsilon)^q$, and for the terms staying on the left from their term with $x_l$. Now the proposition is obvious.

The canonical and power series regularization {\it do not} coincide in general.
For example,
$$
(-\log \varepsilon)^2 = 2 {\rm Li}_{1,1}(1-\varepsilon, 1- \varepsilon) + \zeta (2)
$$
so 
$$
\widetilde \zeta (1,1) = - \frac{1}{2}\zeta (2), \qquad \mbox{but} \quad \zeta (1,1) = 0
$$
Notice that it is definitely false that the multiple $\zeta$-values obtained 
 by {\it canonical} regularization satisfy the first shuffle  relations. For example
$$
0 = \zeta(1)^2 \not = 2\cdot \zeta(1,1) + \zeta(2) 
$$

\begin{proposition} \label{olmn3}
 The difference
$$
{\overline {\rm {\rm Li}}}_{n_1,...,n_l, 1 , ... , 1}(x_1, ... , x_l,1, ... , 1) - 
{\rm Li}_{n_1,...,n_l, 1 , ... , 1}(x_1, ... , x_l,1, ... , 1) 
$$
is an explicitly computable  $\Q$-linear combination of lower depth multiple polylogarithms. In particular if $x_i^N=1$ it is a $\Q$-linear combination of lower depth multiple polylogarithms at $N$-th roots of unity. 
\end{proposition}

{\bf Proof}. 
According to formula (\ref{2.3}) one has 
\begin{equation} \label{ol2}
{\rm Li}_{n_1,...,n_l, 1 , ... , 1}(\frac{a_2}{a_1}, ... , \frac{1}{a_{l}},1, ... , 1,1- \varepsilon) = 
\end{equation}
$$
{\rm I}_{n_1,...,n_l, 1 , ... , 1}(a_1: ... : a_l:1: ... : 1:1- \varepsilon)
$$
 So to compare the power series and canonical regularizations we 
need to compare the asymptotic expansions of  
$$
{\rm Li}_{n_1,...,n_l, 1 , ... , 1}(x_1, ... , x_l,1, ... , 1,1- \varepsilon) \quad \mbox{and} \quad 
{\rm Li}_{n_1,...,n_l, 1 , ... , 1}(x_1, ... , x_l,1- \varepsilon , ... , 1- \varepsilon)
 $$ 
Using (\ref{2.3}) and direct integration we get     
$$
{\rm Li}_{1,...,1}(  1, ..., 1, 1- \varepsilon) = {\rm I}_{ 1, ...,1}(  1 : ... :1:1- \varepsilon) = 
$$
$$  
     \int_{0 < t_1 < ... < t_m < 1 - \varepsilon } \frac{dt_1}{t_1-1}
      \wedge ... \wedge \frac{dt_m}{t_m-1} \quad = \quad
     \frac{(- \log \varepsilon)^m }{m!} 
$$  
Therefore by (\ref{req212}) 
\begin{equation} \label{[!]}
 {\rm Li}_{1,...,1}(  1, ..., 1, 1- \varepsilon) - {\rm Li}_{1,...,1}(  1- \varepsilon, ..., 1- \varepsilon) 
\end{equation}
is a $\Q$-linear combination of lower depth multiple $\zeta$'s. Multiplying 
(\ref{[!]}) by ${\rm Li}_{n_1,...,n_l, 1 , ... , 1}(x_1, ... , x_l,1, ... , 1)$ and using induction on $m-l$ 
we get the proposition.

\begin{lemma} \label{DRc}
The distribution relations holds for the regularized multiple polylogarithms
\end{lemma}

{\bf Proof}. For $l>0$ this is clear since the distribution relations obviously hold for the power series asymptotic expansions. It is easy to see that they also hold for $l= -1$.

As we see the regularization of multiple
  polylogarithms does not lead us to new functions. However 
adding the regularized double shuffle relations we get new relations between 
them, even in the case of
  the multiple $\zeta$-values. D. Zagier 
 conjectured around 1993,  on the basis of his numerical calculations, 
that these relations together with Euler's formula $\zeta(2) = \pi^2/6$ 
provides all the relations between the multiple zeta values.  
Here are some examples.




{\bf Examples}. 1. {\it The weight 3, depth 2 multiple $\zeta$-values}. 
Recall that the canonical regularization gives $\zeta(1) =0$. 
The regularized  double shuffle relations are
$$
0 = \widetilde \zeta(1) \cdot \widetilde \zeta(2) = \widetilde \zeta(2,1) + \zeta(1,2) + \zeta(3)
$$
and 
$$
0 = \zeta(1) \cdot \zeta(2) = 2 \cdot \zeta(1,2) + \zeta(2,1)
$$
The comparison of canonical and power series regularizations leads to one more relation:
$
\widetilde \zeta(2,1) = \zeta(1,2)
$.  Solving them we get $\zeta(1,2) = \zeta(3)$.

Let me remark that the canonical regularization gives   
$\zeta(2,1)
= -2 \zeta(1,2)$, and we also got this from  the second shuffle relation. 

Notice that in this case none of the 
double shuffle relations makes sense without regularization, 
so we must add the
regularized  relations in order to express
$\zeta(1,2)$ via $\zeta(3)$!

2. {\it The weight 4, depth 2 multiple $\zeta$-values}.  
The  double shuffle relations which make sense without regularization  are
$$
\zeta(2) \cdot \zeta(2) = 2 \cdot \zeta(2,2) +  \zeta(4)
$$
and
$$
\zeta(2) \cdot \zeta(2) = 2 \cdot \zeta(2,2) + 4 \cdot \zeta(1,3)
$$
The regularization adds two  more:
$$
0 = \zeta(1) \cdot \zeta(3) =  \zeta(1,3) +  \zeta(3,1) + \zeta(4)
$$
$$
0 = \zeta(1) \cdot \zeta(3) = 2 \cdot \zeta(1,3) + \zeta(2,2) + \widetilde \zeta(3,1)
$$
Finally, the comparison between the canonical and power series regularizations leads to 
$
\widetilde \zeta(3,1) = \zeta(3,1) 
$. 
Solving these equations we get 
$$
\zeta(3,1) = -\frac{5}{4}\zeta(4),\quad \zeta(1,3) = \frac{1}{4}\zeta(4),\quad \zeta(2,2) = \frac{3}{4}\zeta(4),\quad
\zeta(4) = \frac{2}{5}\zeta(2)^2
$$

{\bf 11. Universality of multiple polylogarithms: analytic version}. Suppose that $\omega_{1}, ... , \omega_{n}$
   are such 1-forms that $\int_{\gamma(0)}^{\gamma(1)}\omega_{1}
\circ ... \circ \omega_{n}$ does not depend of a deformation of path $\gamma$ preserving the ends. Then this integral is a multivalued function of variable $y = \gamma(1)$ on $Y$.

\begin{theorem} \label{th2.7}
Suppose that a multivalued analytical function $F(y)$ 
on a rational variety $Y$ is given by a  iterated integral of rational 
1-forms. Then it can be expressed by multiple polylogarithms.
More precisely, there exist rational functions $f^{(i)}_j(y) \in \C(Y)$ such that  
$$
F(y) = \sum_i{\rm I}_{n^{(i)}_{1},...,n^{(i)}_{l^{(i)}}}(f^{(i)}_{1}(y),...,
f^{(i)}_{l}(y)) \quad + \quad C
$$
  where $C$ is a constant 
\end{theorem}

 {\bf Proof}. One can suppose that $Y$ is a (Zariski) open subset in $\C^{N}$. Choose a point $y$ and for any other point $x$ consider line through $y$ and 
$x$. Let us prove that the restriction of our function to this line can be expressed by multiple polylogarithms.

\begin{proposition} \label{p2.8}
  Suppose $f_{i}(t), g_{i}(t)$  are rational functions. Then 
$$
\int_{a_1}^{a_{2}} \frac{dg_{1}(t)}{f_{1}(t)} \circ 
\frac{dg_{2}(t)}{f_{2}(t)} \circ ... \circ \frac{dg_{n}(t)}{f
_{n}(t)} 
$$
is a sum of multiple polylogarithms.
\end{proposition}

{\bf Proof}. Using factorization  and decomposition into 
partial fractions one can reduce the proposition to study of the the following integral
$$
\int_{x_1 \leq 
t_1 \leq ... \leq t_n \leq x_2} (t-a_{1})^{k_1} dt \wedge  (t-a_{2})^{k_2} dt \wedge ... \wedge  
(t-a_{n})^{k_n} dt  
$$
If $k_i \not = -1$ then 
$$
\int_{t_{i-1}}^{t_{i+1}}(t-a_{i})^{k_i} dt = 
\frac{(t_{i+1}-a_{i})^{k_i+1}}{k_i+1} - \frac{(t_{i-1}-a_{i})^{k_i+1}}{k_i+1}
$$ 
 \begin{lemma}  \label{2.9}  The function ${\rm I}(x_1; a_1, ..., a_n; x_2)$ (see (\ref{2**})) 
is a sum of multiple polylogarithms. 
\end{lemma}

\vskip 3mm \noindent
{\bf Proof}. 
Presenting a path from $x_1$ to $x_2$ as a composition of paths from $x_1$ to $0$ and from $0$ to $x_2$ and using the product formula (\ref{BT}) we get 
\begin{equation} \label{2$}
{\rm I}(x_1;a_{1},...,a_n;x_2) \quad= \quad  \sum_{k=0}^{n} {\rm I}(x_1; a_{1},...,a_k;0)   
{\rm I}(0;a_{k+1},...,a_n;x_2) 
\end{equation} 
where we assume that ${\rm I}(x_1; \emptyset ;0) := 1$. This together with 
\begin{equation} \label{QN1}
\int_{\alpha} \omega_1 \circ ... \circ  \omega_n \quad = \quad(-1)^n\int_{\alpha^{-1}} \omega_n \circ ... \circ  \omega_1
\end{equation}
prove the lemma. 
Theorem (\ref{th2.7}) is proved.

{ \bf 12. Multiple $\zeta$ and $L$ functions and their analytic properties. }
We define the multiple zeta functions by 
$$
\zeta(s_1,...,s_m):= \sum_{0< n_1 < ... < n_m}\frac{1}{n_1^{s_1} ...  n_m^{s_m}}
$$
Let us
see what happens if we try to calculate the value
of $\zeta(s_1,s_2)$ at $(s_1,s_2)=(0,0)$. The series defining $\zeta(s_1,s_2)$
are absolutely convergent if $Re (s_1) >1$ and $ Re (s_1+s_2)
>1$. So
$$
\zeta( 0 ,s_2) =  \sum_{0 < k_1 < k_2} \frac{1}{k_2^{s_2}} =
\sum_{0  < k_2} \frac{k_2-1}{k_2^{s_2}} = \zeta( s_2-1) + \zeta( s_2)
$$
On the other hand in the domain of absolute convergence
$$
\zeta(s_1,s_2) + \zeta(s_2,s_1) + \zeta(s_1+s_2)= \zeta(s_1) \zeta(s_2) 
$$
Thus if we assume analytic continuation of the double zeta we
should have
$$
\zeta(s_1,0)= -\zeta(0,s_1) - \zeta(s_1) + \zeta(s_1) \zeta(0)
$$
Now from the first formula we get
$$
\zeta(0,0) = \zeta(-1) - \zeta(0) = \frac{5}{12}
$$
and from the second
$$
2\zeta(0,0) = - \zeta(0) + \zeta(0)^2 = \frac{3}{4}
$$
So we get two different values for $\zeta(0,0)$!.
This argument led to claims in the literature  that 
 $\zeta(s_1,...,s_m)$ ``does not seem to be very nice: in
 particular, it does not have a unique analytic continuation to
 values outside the domain of absolute convergence''.

We show below that multiple $\zeta$, and more
generally multiple Dirichlet $L$-series, do have an analytic
continuation to a meromorphic function in  $(s_1,...,s_m)$
variables with simple 
poles at cetain families of hyperplanes. We compute the residues.

In particular $\zeta(s_1,s_2) =  
s_1/(s_1+s_2) \quad + \mbox{regular terms}$ near $(0,0)$. So 
$\zeta(s_1,s_2) + \zeta(s_2,s_1)$ has no singularity at $(0,0)$,
but indeed
$$
\lim_{s_1 \to 0} \zeta(s_1,0) \not = \lim_{s_2 \to
  0}\zeta(0,s_2) 
$$
and moreover the restriction of   $\zeta(s_1,s_2)$  to line through zero
has a limit at $(0,0)$, but these limits depend on the direction of
the line. 

We view multiple zeta functions as an analytic interpolation of the family of framed 
mixed Tate 
motives over $\Z$ given by their special values at positive integers. 
The multiple zeta functions do not have an Euler product.

{\it Multiple Dirichlet series}. 
Let $\chi_1,...,\chi_m$  be periodic functions on $\Z$ with periods  $N_1,...,N_m$.
In this section we study analytic properties of multiple Dirichlet L-functions:
$$
L_{\chi_1,...,\chi_m}(s_1,...,s_m) = 
\sum_{0< n_1 < ... < n_m}\frac{\chi_1(n_1)\cdot...\cdot\chi_m(n_m)}{n_1^{s_1} ...  n_m^{s_m}}
$$
In particular we have multiple $\zeta$-functions
$$
\zeta(s_1,...,s_m):= \sum_{0< n_1 < ... < n_m}\frac{1}{n_1^{s_1} ...  n_m^{s_m}}
$$

 Product of two multiple Dirichlet $L$-series  is a  sum
of  multiple Dirichlet $L$-series. For example  
$$
L_{\chi_1}(s_1 )L_{\chi_2}( s_2 ) = L_{\chi_1,\chi_2}( s_1, s_2) + L_{\chi_2,\chi_1}(s_ 2,s_1) + L_{\chi_1\cdot\chi_2}( s_1 + s_2)
$$

{\it Analytic continuation of  multiple $\zeta$-functions}. The results below were 
obtained by the author and by M.Kontsevich (unpublished) in 1995. 
    For the  Riemann $\zeta$-function it follows  from the integral formula  
\begin{equation}  \label{disrt}
\Gamma(s) \zeta(s) = \int_0^{\infty} \frac{x^{s-2}\cdot xdx}{e^{x}-1}
\end{equation} 
Namely,  let $x_+^{s}$ be the distribution whose value on a smooth test function $\varphi(x)$ is given for $Re s >1$ by    $<x_+^{s},\varphi(x)>: = \int_0^{\infty}x^{s}\varphi(x) dx$.   
It is well known (see [GS]) that
the generalized function
$x_+^{s-1} /\Gamma(s)$
admits an analytic continuation to a holomorphic function in $s$ such that     
$$
\frac{x_+^{s-1}}{\Gamma(s)}\vert_{s=-n} = \delta^{(n)}(x), \qquad n \geq 0
$$
where $<\delta^{(n)}(x), \varphi(x)>:= (-1)^n\varphi^{(n)}(0)$.
Considering the right hand side of  (\ref{disrt}) as the value of the distribution
$x_+^{s-2}$ on the smooth test function $x/e^{x}-1$ we immediately get analytic continuation and special values at  $s = 0,-1,-2, ...$ of $\zeta(s)$.

For the multiple $\zeta$-functions  we  proceed in a similar way.

\begin{theorem} \label{ancont}
$\zeta(s_1,...,s_m)$ extends to a meromorphic function in   $\C^m$ with the following singularities:

i) a simple pole along the hyperplane $s_m =1$ with the residue $\zeta(s_1,...,s_{m-1})$  

ii) simple poles along the hyperplanes of the following families: 
$$
s_1 + ... +s_m= m-k_1, \quad s_2 + ... +s_m = m-1-k_2, \quad ..., \quad s_{m-1} +   s_m= 2-k_{m-1}
\quad  k_i \geq 0
$$
\end{theorem}

{\bf Proof}. The series defining $\zeta(s_1,...,s_m)$  are absolutely convergent
 if Re $s_m$ $>1$ and Re $(s_1 + ... + s_m) >m$.

\begin{lemma} \label{l1mzf1}
\begin{equation} \label{mzf1 }
\Gamma(s_1)\cdot ... \cdot \Gamma(s_m)\zeta(s_1,...,s_m) = \int_0^{\infty}...\int_0^{\infty} \frac{t_1^{s_1-1}\cdot ... \cdot t_m^{s_m-1} dt_1...dt_m}{(e^{t_1+...+t_m}-1)  (e^{t_2+...+t_m}-1)\cdot ... \cdot(e^{t_m}-1)}
\end{equation}
\end{lemma}

{\bf Proof}. Notice that
$$
\prod_{k=1}^{m}\frac{1}{e^{t_k+...+t_m}-1} \quad = \quad \sum_{k_1  >0}e^{-k_1(t_1+t_2+...+t_m)}\cdot \sum_{k_2>0}e^{-k_2(t_2 +...+t_m)} \cdot ... \cdot  
\sum_{k_m>0}e^{-k_mt_m} = 
$$
$$
\sum_{k_i  >0}e^{-k_1t_1 - (k_1+k_2)t_2 - (k_1+...+k_m)t_m}
$$
 So after the integration of the   right hand side   we get
$$
\Gamma(s_1) \cdot ... \cdot\Gamma(s_m)\sum_{k_1,...,k_m >0} \frac{1}{k_1^{s_1}\cdot ... \cdot (k_1+...+k_m)^{s_m}} 
$$
The lemma is proved.

Unfortunately the expression
$$
\prod_{k=1}^{m}\frac{1}{e^{t_k+...+t_m}-1}
$$ 
has a singularity of type 
$
\prod_{k=1}^{m}( t_k+...+t_m )^{-1} 
$ 
so the right hand side  of (\ref{mzf1 })  is {\it not} the value of  the  distribution 
\begin{equation} \label{maindistrib1}
{t_{1+}^{s_1-1}} \cdot ... \cdot {t_{m+}^{s_m-1}}  dt_1...dt_m
\end{equation} 
 on a smooth test function.
 The situation is saved by   substitution 
$$
x_1 = t_1 + ... + t_m; \quad x_1x_2 = t_2 + ... + t_m; \quad ...; \quad x_1 ... x_m =  t_m \quad <=>
$$
$$
t_m = x_1 ... x_m; \quad t_{l-1} = x_1 ... x_{m-1}(1-x_m); \quad ... ;\quad  t_1 = x_1(1-x_2)
$$
Notice that
$t_i \geq 0 $  if and only if  $ x_1 \geq 0$ and  $ 0 \leq x_i \leq  1$  for $i>1$. 
Using it we rewrite the integral as
$$
\Gamma(s_1)\cdot ... \cdot \Gamma(s_m)\zeta(s_1,...,s_m) = 
$$
$$
\int_0^{\infty}\int_0^{1}...\int_0^{1} 
\frac{     \prod_{p=1}^m x_p^{s_p+ ... + s_m -1} \prod_{q=2}^m (1-x_q)^{s_q-1} dx}
{\prod_{p=1}^m (e^{x_1 ... x_p  }-1)}
$$
 Now
$$
\varphi(x_1,...,x_m):= \frac{ x_1^mx_2^{m-1}...x_m}{\prod_{k=1}^{m} (e^{x_1 ... x_k }-1)}
$$
is a perfectly smooth function and 
$$
  \frac{(s_m-1)\cdot\zeta(s_1,...,s_m)}{\prod_{k=1}^{m-1}\Gamma(s_k + ... + s_m-(m-k+1)) }  
$$
is the value of the distribution
\begin{equation} \label{maindistrib}
 \prod_{k=1}^{m}\frac{x_{k+}^{s_k+...+s_m-(m-k+2)}}{\Gamma(s_k + ... +s_m-(m-k+1)) }    \cdot
\prod_{k=1}^{m-1}\frac{(1-x_{k+1})_+^{s_{k}-1}  }{\Gamma(s_k)} \cdot dx 
\end{equation} 
on the smooth  function $\varphi(x_1,...,x_m)$. 
The distributions appearing as the factors in this product 
are holomorphic functions in the $s$-variables. Their wave fronts are transversal.  
So their product makes sense and it is  an entire function in $s$. The theorem is proved.

{\bf  Example}. One has 
$$
\frac{  (s_2 -1)}{\Gamma(s_1 + s_2 -2) }\cdot\zeta(s_1, s_2)    \quad = \quad \int_0^{\infty}\int_0^{1}\frac{x_1^{s_1+s_2-3}x_2^{s_2-2}  (1-x_2)^{s_1-1} \cdot x_1^2x_2 
\cdot dx_1dx_2}{\Gamma(s_1 + s_2 -2) \Gamma(  s_2 -1)  \Gamma(  s_1 ) \cdot 
(e^{x_1 } -1)(e^{ x_1x_2} -1)}  
$$
Thus the double $\zeta$-function $\zeta(s_1, s_2)$ is the value of the distribution
$$
 \frac{x_{1+}^{s_1+s_2-3}x_{2+}^{s_2-2}(1-x_2)_+^{s_1-1}    \cdot dx_1dx_2}
{\Gamma(  s_1 )\Gamma(  s_2 )   }  
\quad \mbox{on} \quad 
\varphi(x_1,x_2) = \frac{x_1^2x_2}{(e^{x_1} -1)(e^{x_1x_2} -1)}
$$

{\it Residues}. 
The multiple residue of $\zeta(s_1, ... , s_m)$ at $s_1 = ... = s_m =1$ is $1$.
This is the most singular point of the multiple zeta function, 
the only point where the $m$-th multiple residue is non zero.

To find all the residues of the 
multiple zeta function  it is sufficient to find the residues on the hyperplanes.  
In general the residue on a hyperplane $s_l + ... + s_m = (m-l+1) -k$ has the following structure. It is a linear combination with rational coefficients  of   $\zeta(s_1,...,s_{l-1})$ and rational  functions (products of "binomial coefficients" ) on variables $s_{l+1},...,s_m$. The  rational coefficients are products of values of Riemann $\zeta$'s at negative points.  The precise statements see in the theorem \ref{rezid}.

Recall the definition of Bernoulli numbers and their relationship 
with the special values of the Riemann $\zeta$-function:
$$
\frac{x}{e^x-1} \quad =\quad  \sum_{n\geq 0}\frac{B_n}{n!}x^n, 
\qquad
\zeta(1-n) = -\frac{B_n}{n}, \quad n>1; \qquad  \zeta( 0) =  B_1 = -1/2 
$$
For any integer $k \geq -1$ we will use notations
$$
{ s \choose k }: = \frac{\Gamma(s+1)}{k!\Gamma(s-k+1)},\quad k \geq 0; \qquad  { s \choose -1 }: = \frac{1}{s+1} 
$$
  Set $\beta_p:= B_p/p$ for $p>0$ and $\beta_0= 1$, and 
$$
Q_{p_l,...,p_m}(s_l,  ... , s_m):= \prod_{i=m}^{l} {\sum_{\alpha =i}^m(s_{\alpha}
 + p_{\alpha} -1) -1 \choose p_i-1} 
$$
 \begin{theorem} \label{rezid}
 The residue of $\zeta(s_1,...,s_m)$ on the hyperplane 
\begin{equation} \label{hypk}
s_l + ... + s_m = (m-l+1) -k 
\end{equation}
 where $k \geq 0$ is an integer, equals
$$
\sum_{p_{l+1}+...+p_m =k} \beta_{p_{l+1}}\cdot ... \cdot \beta_{p_m} \cdot\zeta(s_1,...,s_{l-1}) \cdot Q_{p_{l+1},...,p_m}(s_{l+1},  ... , s_m)  
$$
\end{theorem}
 
{\bf Examples}. 0. If $l=m$ and $k=0$ we recover part i) of theorem \ref{ancont}.

1.  One has
$$
Res_{s_1+s_2 = 2-k}\zeta(s_1,s_2) = \beta_k \cdot \frac{ \Gamma(s_2+k-1)}{\Gamma(s_2) } = \beta_k \cdot { s_2 + k-2 \choose k-1}
$$

For instance the residues of $\zeta(s_1,s_2)$ at the lines $s_1+s_2 = 2-k$ for $k = 0,1,2,3,4$ are
$$
 \frac{ 1}{s_2-1}, \quad -\frac{ 1}{ 2}, \quad \frac{ s_2}{12 }, \quad
 0, \quad -\frac{ s_2(s_2+1)(s_2+2)}{ 720}, \quad
$$

2. The residue of $\zeta(s_1, s_2, s_3)$ at $s_1 + s_2 +s_3 = 3-k$ is
$$
\sum_{p_2+ p_3 =k} \beta_{p_2} \beta_{p_3} {
   s_2 + s_3 + p_2+ p_3  -3 \choose p_2-1} {  s_3 + p_3   - 2 \choose p_3-1}
$$
  
3. The residue of $\zeta(s_1, s_2, s_3)$ at $  s_2 +s_3 = 2-k$ is
$$
    \beta_{ k} \cdot \zeta(s_1) \cdot  { s_3 +  k   - 2 \choose  k-1}
$$

{\bf Proof}.  It  consists of two steps. 
First we derive a formal expression 
$$
\zeta(s_1,  ... , s_m) "=" \sum_{p_l,...,p_m \geq 0} \beta_{p_l} \cdot ... \cdot 
\beta_{p_m} Q_{p_l,...,p_m}(s_l,  ... , s_m) \times
$$
\begin{equation} \label{formaldec}
\zeta\Bigl(s_1,...,s_{l-2}, s_{l-1} + \sum_{\alpha =l}^m (s_{\alpha}  +  p_{\alpha}-1) 
\Bigr)
\end{equation}
paying no attention to the problems of convergence of the right hand side.

Then we show that only finite number of the terms on the right 
have non zero residues on the hyperplane (\ref{hypk}). Going back to the first step 
it  is  easy  to justify that the residue of the left hand side 
indeed equals to the residue of 
the right hand side. One  could  take residues from the beginning, 
thus avoiding suspicious infinite sums;
however (\ref{formaldec})  seems to be a quite  useful  expression itself: for instance 
 the right hand side is also a finite sum when $s_l,...,s_m$ are non 
positive integers, thus giving a formula for the values of multiple zetas at 
non positive points, see theorem \ref{decf} below. 

{\bf Examples of formula (\ref{formaldec})}.  1. 
$$
\zeta(s_1,  s_2) =  \sum_{p \geq 0} \beta_p {  s_2+p-2 \choose p-1} \zeta(s_1+s_2+p-1)
$$


$$
\zeta(s_1,  s_2, s_3) = 
$$
\begin{equation} \label{adr11}
\sum_{ p_2, p_3 \geq 0}   \beta_{ p_2}  \beta_{ p_3} {s_2+s_3 + p_2 +  p_3 -3 \choose  p_2-1} {  s_3 + p_3   - 2 \choose  p_3-1} \zeta(s_1+s_2+ s_3 +  p_2 +  p_3   - 2) 
\end{equation}

{\bf Proof}. Let us spell  the details  of the 
 proof of formula (\ref{adr11}). Recall that
$$ 
 \Gamma(s_1)\Gamma(s_2)\Gamma(s_3) \cdot \zeta(s_1,s_2,s_3) 
$$
is the value of the distribution
 \begin{equation} \label{dr111}
x^{s_1 + s_2 + s_3 -4}_{1+} \cdot  x_{2+}^{s_2 + s_3-3}  (1-x_2)_+^{s_1-1} \cdot x_{3+}^{s_3-2}(1-x_3)_+^{s_2-1}dx 
 \end{equation}
on $\varphi(x_1,x_2,x_3)$. Consider the following expansion of $\varphi(x_1,x_2,x_3)$ into formal   series:
\begin{equation} \label{dr333} 
 \frac{x_1}{e^{x_1}-1}\cdot \sum_{p_2,p_3 \geq 0}\frac{B_{p_2}}{p_2!}\frac{B_{p_3}}{p_3!} (x_1x_2)^{p_2 }(x_1x_2x_3)^{p_3}
\end{equation} 
Notice that for given $x_1$ this series converge  only in a small domain near $(x_2,x_3)=(0,0)$. So if we simply evaluate the distribution (\ref{dr111}) 
on each term of this series and then take a sum it  {\it a priory} 
does not converge. 
However  let us compute for given $p_2,p_3$ the corresponding triple integral.  It is clearly a product of the constant  
\begin{equation} \label{for1}
\frac{B_{p_2}}{p_2!}\frac{B_{p_3}}{p_3!}
\end{equation}
 and  three one dimensional integrals.  The integral along $x_1$ equals
\begin{equation} \label{for2}
\Gamma(s_1 + s_2 + s_3 + p_2 + p_3 -  2) \zeta(s_1 + s_2 + s_3 + p_2 + p_3 -  2)
\end{equation}
The integrals along $x_2$  and $x_3$ are Euler's $B$-integrals, so the integral along $x_2$ gives
\begin{equation} \label{for3}
\frac{\Gamma(s_2 + s_3 + p_2 + p_3 -2)\Gamma(s_1)}{\Gamma(s_1 + s_2 + s_3 + p_2 + p_3 -2)}
\end{equation}
and the integral along $x_3$ equals
\begin{equation} \label{for4}
\frac{\Gamma(s_3+ p_3 -1)\Gamma(s_2)}{\Gamma(s_2+ s_3 +p_3 -1)}
\end{equation} 
 Multiplying the  factors  (\ref{for1}) - (\ref{for4}), dividing them by $\Gamma(s_1)\Gamma(s_2)\Gamma(s_3)$ and   we  obtain just the term   of the formula (\ref{adr11}) corresponding to given $p_2,p_3$. 
 
In general we proceed the same way. First of all we  write $\varphi(x_1,...,x_m)$  as series 
$$
   \prod_{i=1}^{l-1} \frac{x_1   ...   x_i}{e^{x_1   ...   x_i}-1} \cdot \sum_{p_{l},...,p_m \geq 0} \frac{B_{p_l}}{p_l!}\cdot ... \cdot \frac{B_{p_m}}{p_m!}(x_1   ...   x_l)^{p_l} \cdot ... \cdot (x_1 ... x_m)^{p_m}
$$
For given $p_l,...,p_m$ the $m$-dimensional integral splits to a product of
a constant  and the following  integrals: the integral along $x_1,...,x_{l-1}$, and 
integrals along $x_k, k \geq l$,  which are $B$-functions.
Taking their product we get  the needed summand of the  formula (\ref{formaldec}).

Now we calculate similarly the residue of $\zeta(s_1,  ... , s_m)$ on the hyperplane (\ref{hypk}).  Let us show that   we  get a (finite) sum of the residues of the terms of (\ref{formaldec})   with  $p_l + ... + p_m \leq (m-l+1) - k$. 

For simplicity  let us
 spell the details in the case $l=2, m=3$.  The distribution (\ref{dr111}) is a direct product of the three distributions sitting on the lines with coordinates $x_1$, $x_2$ and $x_3$.  Notice that
$$
\frac{x_{2+}^{s_2 + s_3-3}  (1-x_2)_+^{s_1-1}}{\Gamma(s_2 +s_3 -2)}|_{s_2 +s_3 -2 = -k} 
\quad = \quad \delta^{(k)}(x_2) \cdot (1-x_2)_+^{s_1-1}
$$
Therefore the residue of the  distribution (\ref{dr111}) at  $s_2 + s_3 = 2-k$ 
evaluated on a test function $\psi(x_1,x_2,x_3)$  will be zero provided 
 \begin{equation} \label{condicia}
\psi^{(i)}_{x_2} (x_1, 0,x_3) = 0 \quad \mbox{ for $i >k$}
\end{equation} 
Now if we subtract from the function $\varphi(x_1,x_2,x_3)$ the terms of the sum (\ref{dr333}) with $p_2 + p_3 \leq k$ we get a function satisfying the condition (\ref{condicia}) above. So calculating the residue we can pay attention only to the terms with $p_2 + p_3 \leq k$.
 
Now return to the general case and look more carefully  at the    
expression
\begin{equation} \label{expresso} 
 \frac{\Gamma(s_{l} + ... +s_m +  p_l + ... + p_m  - (m-l+1))}{(p_l-1)!\Gamma(s_{l} + ... +s_m +  p_{l+1} + ... + p_m  - (m-l ))}
 \end{equation}
which is the first factor in the product defining the rational function \linebreak $Q_{p_l,...,p_m}(s_l,  ... , s_m)$.    The denominator in (\ref{expresso}) is non zero on the hyperplane
(\ref{hypk}) only if  $p_{l+1}+...+p_m  \geq k$. This means that (\ref{expresso}) has a nonzero residue on the hyperplane (\ref{hypk}) only if $p_l =0$ and $p_{l+1}+...+p_m   = k$. In this case the residue of (\ref{expresso}) is equal to $1$. The theorem is proved.

{\it Special values  when some of $s_i$ are   not positive
  integers}. 
\begin{theorem} \label{decf}
If $s_l,...,s_m$ are  non positive   integers then the right hand side of formula (\ref{formaldec})    gives an expression of  $\zeta(s_1,  ... , s_m)$   as a finite linear combination of  multiple zeta functions of depth $l-1$.
 \end{theorem}
 
{\bf Examples}. 1. If  $s_3$   is a   non positive integer then 
$$
\zeta(s_1, s_2, s_3) =
\sum_{ p_3 \geq 0}     \beta_{ p_3} \zeta( s_1,s_2 + s_3 + p_3 -1)
{  s_3 + p_3   - 2 \choose  p_3-1}
$$
 where the sum on the right is finite.
   
2. If  $s_3, s_4$   are negative integers then 
$$
\zeta(s_1, s_2, s_3, s_4) =
$$
$$
\sum_{ p_3, p_4 \geq 0} \beta_{ p_3} \beta_{ p_4}
\zeta( s_1,s_2 + s_3 + s_4 + p_3 + p_4-2){s_3+s_4 + p_3 +  p_4 -3 \choose  p_3-1}{ s_4+  p_4  - 2 \choose  p_4-1}
$$

{\bf Proof}. Suppose that $s_l,...,s_m$ are  non positive   integers.
We claim that the terms of (\ref{formaldec}) could be non zero only if the following conditions are satisfied:
\begin{equation} \label{kkk} 
 \sum_{i=l}^m (s_i + p_i -1) \leq 0 \quad \mbox{for $l \leq m$} 
\end{equation}
 Indeed,
$$
{s_m + p_m -2\choose p_m-1} = \frac{\Gamma(s_m + p_m -1)}{(p_m-1)!\cdot \Gamma(s_m)}
$$
and $1/\Gamma(s_m)=0$ for non positive   integer $s_m$. So this binomial coefficient could be nonzero only if $s_m + p_m -1\leq  0$. 
And so on. 
Notice that only a finite number of non positive integers $p_l,...,p_m$ satisfy the  conditions (\ref{kkk}) provided $s_l,...,s_m$ are  non positive   integers.

Similarly to the proof of theorem \ref{rezid} one can justify formula (\ref{formaldec}) in our case.

{\it Analytic continuation of multiple Dirichlet series}. Set 
$$
\varphi_{\chi_1,...,\chi_m}(x_1,...,x_m):= \sum_{0< m_i \leq  N_i } \prod_{i=1}^m\frac{e^{-m_ix_1...x_i}\cdot\chi_i(m_1+...+m_i) x_1...x_i}{1-e^{-N_ix_1...x_i}}
$$
For example
$$
\varphi_{\chi_1 ,\chi_2}(x_1,x_2) = \sum_{0< m_i \leq  N_i }  \frac{e^{-m_1x_1  - m_2x_1x_2}\cdot\chi_1(m_1 ) \chi_2(m_1 +m_2)x^2_1 x_2}{(1-e^{-N_1x_1})(1-e^{-N_2x_1 x_2})}
$$

\begin{theorem} \label{anconti}Let $\chi_1,..., \chi_m$ are  
periodic functions with periods $N_1,..., N_m$. Then

a) $
\Gamma(s_1)\cdot ... \cdot\Gamma(s_m)\cdot L_{ \chi_1,..., \chi_m}(s_1,...,s_m)
$ 
equals to the value of the distribution (\ref{maindistrib})
on $\varphi_{\chi_1,...,\chi_m}(x_1,...,x_m)$.

b) $L_{\chi_1,..., \chi_m}(s_1,...,s_m)$ extends to a meromorphic 
function in   $\C^m$ which    has  simple poles along the hyperplanes 
of the following families: 
$$
s_1 + ... +s_m = m-k_1, \quad s_2 + ... +s_m= m-1-k_2, \quad ..., \quad s_{m-1} +   s_m = 
2-k_{m-1}\qquad  k_i \geq 0
$$
  and   a simple pole along the hyperplane $s_m =1$  with the residue 
$$ 
\sum_{a=1}^{N_m}\chi_m(a)\cdot  L_{\chi_1,..., \chi_{m-1}}(s_1,...,s_{m-1})
$$  
\end{theorem} 
 
{\bf Proof}. Part b) follows immediately from the part a). We  spell the details of the proof of part a) in the case $m=2$. The general case is   similar. By definition
$$
\Gamma(s_1)\Gamma(s_2)\cdot L_{ \chi_1, \chi_2}(s_1,s_2) = \Gamma(s_1)\Gamma(s_2)\cdot \sum_{n_i >0}\frac{\chi_1(n_1)\chi_2(n_1 + n_2)}{n_1  ^{s_1}(n_1 + n_2)^{s_2}} = 
$$
$$
\sum_{n_i >0}\int_0^{\infty}\int_0^{\infty} e^{-n_1 t_1} e^{-(n_1 + n_2) t_2}
\chi_1(n_1)\chi_2(n_1 + n_2)  t_1^{s_1-1} t_2^{s_2-1} d t_1d t_2
$$ 
Writing $n_i = k_iN_i + m_i$ where $k_i \geq 0$ and $0< m_i \leq N_i$ 
and setting 
$$
C_{\chi_1,\chi_2}( t_1,  t_2) := \sum_{0< m_i \leq  N_i } e^{-m_1( t_1 + t_2)}e^{-m_2  t_2 } \chi_1( m_1) \chi_2( m_1 +  m_2)
$$
we  rewrite the last expression as
$$
\int_0^{\infty}\int_0^{\infty} \Bigl(\sum_{k_i \geq 0}e^{-k_1N_1( t_1 + t_2)}  e^{-k_2N_2  t_2 }\Bigr)C_{\chi_1,\chi_2}( t_1,t_2) t_1^{s_1-1} t_2^{s_2-1} d t_1d t_2
=
$$
$$
\int_0^{\infty}\int_0^{\infty}  \frac{C_{\chi_1,\chi_2}( t_1,  t_2)  t_1^{s_1-1} t
_2^{s_2-1}}{(1 - e^{N_1( t_1+ t_2)}  )(1 - e^{N_2  t_2 } )}d t_1d t_2
$$

Making the usual substitution 
$$
t_1 + t_2 = x_1, \quad t_2 = x_1x_2  \qquad <=> \qquad t_1 = x_1(1-x_2), \quad  t_2 = x_1 x_2
$$
we get the integral representation from the part a).

\section{ The mixed Tate motives over a ring of 
$S$-integers in a number field}

The triangulated category of mixed motives over an arbitrary 
 filed $F$ 
 has been constructed, independently, in the fundamental works by M. Levine [L] and 
V. Voevodsky [V]. Their theories lead to  equivalent categories ([L]). 

The Beilinson-Soul\'e vanishing conjecture blocks any attempts 
to construct an abelian category of mixed motives over  an arbitrary 
 filed $F$. 
However if $F$ is a number field 
the vanishing conjecture follows from  the results of Borel and Beilinson [Bo], [Be1]. 
Using this one can deduce 
the existence of 
the abelian category ${\cal M}_T(F)$ of mixed Tate motives 
over a number field $F$ satisfying all the desired properties, see chapter 5 of [G8] or [L1]. 

Unfortunately  so far a construction similar to [L] or [V] of the triangulated category of 
motives over a 
Noetherian scheme, or even over ${\rm Spec}\Z$, with the expected ${\rm Ext}$'s,  
is still far from being available. 

Let ${\rm Spec}{\cal O}_{F, S}$ be the scheme obtained by deleting 
a finite set $S$  of closed points from the spectrum of the 
ring of integers of a number field $F$. We  show in this chapter that, 
using the Tannakian formalism,  one can nevertheless construct an abelian category 
${\cal M}_T({\cal O}_{F, S})$ 
of mixed motives over the scheme ${\rm Spec}{\cal O}_{F, S}$ 
with all the desired properties, 
including the formula expressing the Ext groups with $K$-theory 
of the ring ${\cal O}_{F, S}$.

We start  by  recalling
 the Tannakian formalism for mixed Tate categories addressed on the language 
of framed objects (see [BMS], [BGSV]).  
For the  proofs (in a more  general setting) see chapter 3 of [G9].

{\bf 1. A review of the Tannakian formalism for mixed Tate categories}. Let $K$ be a characteristic zero field. 
Let ${\cal M}$ be a Tannakian  $K$-category with an invertible object $K(1)$. So 
in particular ${\cal M}$ is an abelian tensor  category. 
Set $K(n):= K(1)^{\otimes n}$. 
Recall ([BD2]) that the pair $({\cal M}, K(1))$ is called a {\it mixed Tate category} 
if the objects $K(n)$
are mutually nonisomorphic, any simple object is isomorphic to one of them and
$Ext_{{\cal M}}^1(K(0),K(n)) = 0$ if $n \leq 0$. A pure functor 
between two mixed Tate categories $({\cal M}_1, K(1)_1)$ and $({\cal M}_2, K(1)_2)$ is a 
tensor functor $\varphi: {\cal M}_1 \lra {\cal M}_2$  equipped 
with an isomorphism $\varphi(K(1)_1) = K(1)_2$. 

It is easy to show that any object $M$ of a mixed Tate category has a canonical 
weight filtration $W_{\bullet}M$ indexed by $2\Z$ 
such that $gr^W_{2n}M$ is a direct sum of copies of $K(-n)$, and 
morphisms are strictly compatible with the weight filtration. The functor
$$
\Psi = \Psi_{{\cal M}}:  {\cal M} \longrightarrow Vect_{\bullet}, \quad M \longmapsto 
\oplus_{n} Hom_{{\cal M}}(K(-n), gr^W_{2n}M)
$$
  to the category of graded $K$-vector spaces is a fiber functor. Let 
$$
L_{\bullet}({\cal M}):= Der ( \Psi) : = \{F \in End \Psi | F_{X \otimes Y} = F_{X} 
\otimes id_{Y} +  id_{Y} \otimes F_{Y}\}
$$  
be the space of its derivations. 
It is a negatively graded pro-Lie algebra over $K$ called 
 the {\it fundamental Lie algebra} of the mixed Tate category ${\cal M}$. 

Let $\widetilde \Psi$ be the fiber functor to the category of finite dimensional 
$K$-vector spaces obtained from $\Psi$ by forgetting the grading. Then the automorphisms of the fiber functor 
respecting the tensor structure provide 
a pro-algebraic group scheme over $K$, denoted 
$ Aut^{\otimes} \widetilde \Psi$. It is a semidirect product of ${\Bbb G}_m$ 
and a pro-unipotent group scheme 
$U({\cal M})$. The pro-Lie algebra $L_{\bullet}({\cal M})$ is the 
Lie algebra of  $U({\cal M})$. The action of ${\Bbb G}_m$ provides 
a grading on $L_{\bullet}({\cal M})$. 

According to the Tannakian formalism the functor $\Psi$ provides an equivalence 
between the category $\cal M$ and 
 the category of finite dimensional modules over the
pro-group scheme $ Aut^{\otimes} \widetilde \Psi$.  This category is naturally equivalent 
to the category of graded finite dimensional modules over the group scheme $U({\cal M})$.
 Since $U({\cal M})$ is  pro-unipotent,  the last category is equivalent  to the category 
graded finite dimensional modules over the
graded pro-Lie algebra $L_{\bullet}({\cal M})$.  

Pure functors  $({\cal M}_1, K(1)_1)  \lra ({\cal M}_2, K(1)_2)$ between the mixed Tate categories 
 are in $1-1$ correspondence with the graded Lie algebra 
morphisms  $L_{\bullet}({\cal M}_2) \lra L_{\bullet}({\cal M}_1)$. 

 Let  ${\cal U}_{\bullet}({\cal M}):= End(\Psi)$ be the space of all endomorphisms of the 
fiber functor $\Psi$. It is a graded 
Hopf algebra isomorphic to the universal enveloping 
of the Lie algebra $L_{\bullet}({\cal M})$. 

Recall that a Lie coalgebra is a vector 
space ${\cal D}$ equipped with a linear map 
$
\delta:  {\cal D} \lra \Lambda^2{\cal D} 
$
such that the composition 
$
{\cal D} \stackrel{\delta}{\lra}\Lambda^2{\cal D} \stackrel{\delta \otimes 1 - 1 \otimes 
\delta}{\lra}
\Lambda^3{\cal D} 
$ 
is zero. If ${\cal D} $ is finite dimensional then it is a Lie coalgebra if 
and only if its dual is a Lie algebra.

 Let $G$ be a 
unipotent algebraic group over $\Q$. Then the ring of regular functions 
$\Q[G]$ is a Hopf algebra with the coproduct induced by the multiplication in $G$.
The (continuous) dual of its completion at the group unit $e$ is isomorphic 
to the universal enveloping 
of the Lie algebra ${\rm Lie}(G)$ of $G$. Let $\Q[G]_0$ be the ideal of functions
equal to
zero at $e$. Then the coproduct induces on $\Q[G]_0/\Q[G]_0 \cdot \Q[G]_0$ the 
structure of a Lie coalgebra  dual to ${\rm Lie}(G)$.

In this paper we employ the duality $V \lms V^{\vee}$ between the Ind- and pro-objects 
in the category of finite dimensional $K$-vector spaces. The graded dual Hopf algebra ${\cal U}_{\bullet}({\cal M})^{\vee}:= 
\oplus_{k \geq 0} {\cal U}_{-k}({\cal M})^{\vee}$ can be identified with the Hopf algebra  
of regular functions on the pro-group scheme $U({\cal M})$. 
Therefore its  quotient by the square of the augmentation ideal 
$$
{\cal L}_{\bullet}({\cal M}):= \frac{{\cal U}_{>0}({\cal M})^{\vee}}
{{\cal U}_{>0}({\cal M})^{\vee} \cdot {\cal U}_{>0}({\cal M})^{\vee}}
$$
is a Lie coalgebra. 
The cobracket $\delta$ on ${\cal L}_{\bullet}({\cal M})$ 
is induced by the restricted coproduct 
\begin{equation} \label{8-23.2/99}
 \Delta'(X):= \quad \Delta (X)-  (X \otimes 1 + 1 \otimes X)
\end{equation}
on ${\cal U}_{\bullet}({\cal M})^{\vee}$. 
The  graded dual to the Lie coalgebra ${\cal L}_{\bullet}({\cal M})$ 
is identified with the fundamental Lie algebra 
$L_{\bullet}({\cal M})$. Below we  recall a more efficient 
  way to think about  it. 

 {\bf 2. The Hopf algebra of framed objects in a mixed Tate category 
(cf. [BMS], [BGSV], [G9])}. Recall (cf. [ES], chapter 8) that 
a Hopf algebra over a field $k$ is a $k$-vector space $A$ with 
the following additional structures:

1. Multiplication $\mu: A \otimes A \lra A$.

2. Comultiplication $\Delta: A \lra A \otimes  A$.

3. Unit $i:A \lra k$. 4. Counit $\varepsilon: k \lra A$. 5. Antipode $S: A \lra A$.

They must obey the following properties:

1. The map $\mu$ defines a structure of an associative algebra on $A$ with unit $i(1)$.

2. The maps $\Delta$ and $\varepsilon$ define 
a structure of a coassociative algebra on $A$.

3. The maps $\Delta: A \lra A \otimes A$ and $\varepsilon: A \to k$ 
are homomorphisms of algebras.

4. The map $S$ is a linear isomorphism satisfying the relations
\begin{equation} \label{4.16.01.1}
\mu \circ (S \otimes {\rm id}) \circ \Delta  = \mu \circ ({\rm id} \otimes S) 
\circ \Delta  = i \circ \varepsilon
\end{equation}

A bilagebra is a $k$-vector space equipped with the structures $1-4$ satsfying the axioms
$1-3$. 

\begin{lemma} \label{4.15.01.1}
Suppose that $A = A_{\bullet}$ is a commutative bialgebra 
graded by non negative integers $n \geq 0$ such that 
$A_0 =k$. Then there exists unique 
antipode map $S$ on $A_{\bullet}$. 
\end{lemma}
 
{\bf Proof}. Indeed, the condition (\ref{4.16.01.1}) determines uniquely the restriction of the map 
$S$ to $A_n$ by induction. One has $S|A_0 = {\rm id}$. For example 
$S|A_1 = -{\rm id}$ and  $S|A_2 = -{\rm id} + \mu \Delta_{1,1}$ 
where $\Delta_{1,1}$ is the 
$A_1 \otimes A_1$ component of the coproduct. 

Let $n \geq 0$. Say
that   $M$ is an $n$-framed  object of ${\cal M}$  if it is supplied with a 
nonzero morphisms $v_{0} : K(0) \longrightarrow 
gr^W_{0}M$ and   $f_{n}: gr^W_{-2n} M \longrightarrow K(n)$.  

Consider the  coarsest equivalence
relation on the set of all $n$-framed  objects for which $M_1 \sim M_2$
if there is a map $ M_1 \to  M_2$ respecting the frames.  For example 
replacing $M$ by 
$W_0M/W_{-2p-2}M$ we see that any $n$-framed
 object is equivalent to a one $ M$ with $W_{-2-2n} M = 0$, $W_{0} M= M$.
Let ${\cal A}_n({\cal M})$ be the set of equivalence classes.
It  has a structure of an abelian group  
with the composition law defined as follows:
$$
[ M, v_{0},f_{n}] + [ M', v_{0}',f_{n}'] \quad := \quad [M \oplus M', (v_{0}, 
v_{0}'), f_{n} + f_{n}' ]
$$
It is straitforward to check that the composition law is well 
defined on equivalence classes of framed objects. Indeed, if 
$\varphi: \widetilde M \lra M$ 
is a morphism providing an equivalence of the framed objects 
$[ \widetilde M, \widetilde v_{0}, \widetilde f_{n}] \sim 
[ M, v_{0}, f_{n}]$ then 
$ \varphi \oplus {\rm id}: \widetilde M \oplus M' \lra M \oplus M'$ 
provides an equivalence. 

The neutral element is $K(0) \oplus K(n)$ with the obvious frame. 
The inversion is given by 
$$
-[M, v_{0}, f_{n}]:= \quad [M, -v_{0}, f_{n}] \quad = \quad [M, v_{0}, -f_{n}]
$$
See ch. 2 of [G9] for a proof of the fact that $K(0) \oplus K(n)$ is 
indeed a neutral elelment, and  this formula defines an 
inverse.

The composition $f_0\circ v_0: K(0) \to K(0)$ provides an 
isomorphism ${\cal A}_0({\cal M}) = K$. 

The tensor product   induces the commutative and associative multiplication
$$
\mu: {\cal  A}_{k}({\cal M}) \otimes {\cal  A}_{\ell}({\cal M})\to {\cal A}_{k+ \ell}({\cal M})
$$
One varifies that it is well defined on equivalence classes using an argument similar 
to the one used for the additive structure on ${\cal  A}_{k}({\cal M})$. 
The init is given by $1 \in K = {\cal A}_0({\cal M})$. 
The counit is the projection of ${\cal  A}_{\bullet}({\cal M})$ 
onto its $0$-th component. 

Let us define the comultiplication
$$
\Delta = \bigoplus_{0 \leq k \leq n} \Delta_{k, n-k}: \quad {\cal  A}_{n}({\cal M})
\to \bigoplus_{0 \leq k \leq n } 
{\cal  A}_{k}({\cal M})\otimes
{\cal  A}_{n-k}({\cal M})
$$

 Choose a basis $\{b_i\}$, where $1 \leq i \leq m$,  of ${\rm Hom}_{\cal M} (K(p),gr_{-2p}^W  M)$   and  the dual basis $\{b'_i\}$ of $({\rm Hom}_{\cal M} (gr_{-2p}^W M, K(p) )$.  
Then
$$
\Delta_{p, n-p } [M,v_{0},f_n]:=  \quad \sum_{i=1}^m [M,v_{0},b'_i] \otimes [M, b_i, f_n](-p)
$$
In particular $\Delta_{0,n} = id \otimes 1$ and $\Delta_{n,0} = 1 \otimes id$.  
Set ${\cal A}_{\bullet}({\cal M}):=\oplus {\cal A}_n({\cal M})$.

\begin{theorem} \label{8-23.3/99}
a) ${\cal  A}_{\bullet}({\cal M})$ has a structure of a 
graded Hopf algebra over $K$ with the commutative multiplication
$\mu$ and the comultiplication $\Delta$.

 b) The Hopf algebra ${\cal  A}_{\bullet}({\cal M})$ is 
canonically isomorphic to the Hopf algebra  ${\cal U}_{\bullet}({\cal M})^{\vee}$.
\end{theorem}

{\bf Proof}. a) Let us show that the coproduct is well defined on equivalence 
classes of framed objects. It is sufficient to prove this for equivalences given by 
injective and surjective morphisms in ${\cal M}$. Indeed, if 
$\varphi: M \lra M'$ respects the frames in $M$ and $M'$ then the projection 
$M \lra M/{\rm Ker} (\varphi)$ and injection $M/{\rm Ker} (\varphi) \hookrightarrow  M'$ 
also respect the frames. Let us suppose that 
$$
\varphi : [M, v_0, f_n] \hookrightarrow [M', v'_0, f'_n]
$$
is an equivalence. Choose a basis $\{b_i\}$ of 
${\rm Hom}_{\cal M} (K(p),gr_{-2p}^W  M')$ 
such that $b_1, ..., b_s$ is a basis in 
${\rm Hom}_{\cal M} (K(p),gr_{-2p}^W  M')$. Then 
$[M', v_0, b_{i}']=0$ for $i>s$. Indeed, we may assume that 
$W_{-2p-2}M'=0$, and also ${\rm Gr}^W_0M'= K(0)$. 
Then there is a natural injective morphism 
$K(p) \oplus M \cap M' \hookrightarrow M'$ respecting the frames, and the projection 
$$
K(p) \oplus (M\cap M') \lra K(p) \oplus (M\cap M')/W_{-2}(M\cap M') = 
K(p) \oplus K(0)
$$ The statement is proved. 
The arguments in the case of the projection are completly similar 
(and can be obtained by dualization).

It is straitforward to show 
that $\Delta$ and $\mu$ satisfy the Hopf axiom, i.e. the condition 2 above. 
It is also straitforward to check that $\Delta$ is coassociative. 
The part a) of the theorem is proved.

b)  We follow the proof of theorem 3.3 in [G9] making some necessary 
corrections. The canonical isomorphism $\varphi: {\cal  A}_{\bullet}({\cal M}) \to 
{\cal U}_{\bullet}({\cal M})^{\vee}$
is constructed as follows.
Let $F \in End(\Psi)_n$ and $ [M,v_{0}, f_n]\in {\cal  A}_{n}({\cal M})$. 
Denote by $F_M$ the endomorphism of $\Psi(M)$ provided by $\Psi$. Then 
$$
<\varphi([M,v_{0}, f_n]) , F>:= \quad <f_n, F_M(v_{0})>
$$
It is obviously well defined on equivalence classes of framed objects. 
One easyly checks that $\varphi$ is a morphism of graded Hopf algebras. 

Let us show that $\varphi$ is surjective. Recall that for a commutative non negatively graded 
Hopf algebra ${\cal A}_{\bullet}$ the space of primitives 
$
{ CoLie} ({\cal A}_{\bullet}):= 
{\cal A}_{>0}/{\cal A}_{>0}^2
$ 
has a Lie coalgebra structure, and the dual to the 
universal enveloping of the dual Lie algebra is 
canonically isomorphic to  ${\cal A}_{\bullet}$. 
Since $\varphi$ is a map of Hopf algebras it provides a morphism of the Lie coalgebras 
$$
\overline \varphi: {CoLie} ({\cal A}_{\bullet}({\cal M})) \lra 
{CoLie} (U_{\bullet}({\cal M})^{\vee})
$$
Let us show that this map is surjective. 
Let $f$ be a functional on ${CoLie} (U_{\bullet}({\cal M})^{\vee})$ which is zero 
on the graded components of the degree $\not = n$. Consider $f$ as a functional 
on $U_{\bullet}({\cal M})^{\vee}$, and denote by ${\rm Ker}(f)$ its kernel. 
Then $U_{\bullet}({\cal M})^{\vee}/{\rm Ker}(f) = K_{(-n)}$ is a one dimensional 
space sitting in degree $-n$. The graded $k$-vector space 
$K_{(0)} \oplus U_{\bullet}({\cal M})^{\vee}/{\rm Ker}(f) $ has a 
$U_{\bullet}({\cal M})^{\vee}$-module structure: an element 
$X \in U_{>0}({\cal M})^{\vee}$ sends 
$1 \in K_{(0)} \lms f(X) \in K_{(-n)}$ and annihilates $K_{(-n)}$. Since by the 
definintion of $f$ the square of the augmentation ideal acts by zero, 
we get a well defined action of $U_{\bullet}({\cal M})^{\vee}$. 
Therefore the map $\overline \varphi$ is surjective. 

Dualizing $\overline \varphi$ we get an injective map of 
Lie algebras, and hence an injective map of the corresponding 
universal enveloping algebras. Dualizing it we prove the surjectivity of $\varphi$.

Now let us show that $\varphi$ is injective. Using the Tannaka theory we may assume that 
${\cal M}$ is the category of finite dimensional representations of 
the Hopf algebra $U_{\bullet}({\cal M})$. Suppose that 
$\overline \varphi\Bigl([M, v_0, f_n] \Bigr) =0$. Consider the cyclic submodule 
$U_{\bullet}({\cal M}) \cdot K_{(0)}$. It has no non zero 
components in the degree $-n$, since otherwise $\varphi[M, v_0, f_n] \not = 0$. 
Thus there are maps 
$$
K_{(0)} \oplus K_{(-n)} \longleftarrow U_{\bullet}({\cal M}) \cdot K_{(0)} \oplus K_{(-n)} \lra M / W_{<-n}M
$$
respecting the frames. Thus $[M, v_0, f_n]=0$. 

The part b), and hence the theorem are proved.

Under the  equivalence between the  tensor category ${\cal M}$ and the category of 
graded finite dimensional comodules 
over the Hopf algebra 
${\cal  A}_{\bullet}({\cal M})$ an object $M$ of ${\cal M}$ corresponds to 
the graded comodule 
$\Psi (M)$ with 
${\cal  A}_{\bullet}({\cal M})$-coaction $\Psi (M )
\otimes \Psi (M^{\ast}) \longrightarrow {\cal  A}_{\bullet}({\cal M})$ given by the formula 
\begin{equation} \label{8-24.5/99}
x_{m} \otimes y_{n} \longrightarrow \quad \mbox{the class of $M$ 
framed by $x_{m}, y_{n}$}
\end{equation}
 We call the right hand side the {\it matrix coefficient} of $M$ corresponding to $x_{m}, y_{n}$.

The restricted coproduct
$\Delta'$  provides the quotient 
$
{\cal A}_{\bullet}({\cal M})/({\cal A}_{>0}({\cal M}) )^2
$ 
with the structure of a graded Lie coalgebra  with cobracket 
$\delta$. It is canonically isomorphic to ${\cal L}_{\bullet}({\cal M})$.

\begin{lemma} \label{2.3.01.2}
Each equivalence class of $n$-framed objects contains a unique 
minimal representative, which appears as a subquotient in any 
$n$-framed object from the given equivalence class. 
\end{lemma}

{\bf Proof}. Suppose we have morphisms $M_1 \stackrel{f}{\longleftarrow} N 
\stackrel{g}{\longrightarrow} M_2$ between the $n$-framed objects respecting the frames. For any $n$-framed object $X$ we can assume that ${\rm gr}^W_{2m}X =0$ unless $ -n \leq m \leq 0$ and 
${\rm gr}^W_{-2n}X = \Q(n)$; ${\rm gr}^W_{0}X = \Q(0)$. Taking the subquotient ${\rm Im}(f)$ 
of $M_1$ we may suppose that $f$ is surjective. Then ${\rm gr}^W_{0}{\rm Ker}(f) = 
{\rm gr}^W_{-2n}{\rm Ker}(f) = 0$. Therefore $g({\rm Ker}(f))$ has the same property, and hence $M_2$ is equivalent to $N_2':= M_2/{\rm Ker}(f)$. The lemma follows from these remarks.

{\bf 3. The triangulated category of mixed Tate motives over $F$ (cf. [V])}. 
Let ${\cal D}{\cal M}_F$ be  the triangulated category  of 
 motives over a field $F$ as constructed by V. Voevodsky in   [V].  
Recall that an object in ${\cal D}{\cal M}_F$
is a complex of regular (but not necessarily projective) varieties $X_1
\longrightarrow X_2 
\longrightarrow ... \longrightarrow X_n$ where the arrows are given by
finite 
correspondences and the composition of any two successive arrows is
zero. 
The Tate motive $\Q(1)$ is defined as the object $P^1 \lra \ast$ of ${\cal D}{\cal M}_F$
 sitting in degrees $[-2,-3]$, where $\ast$ is a point. 
One has the  basic  formula conjectured by Beilinson: 
\begin{equation} \label{vvvvv}
Ext^i_{{\cal D}{\cal M}_F}(\Q(0), \Q(n)) = 
gr^{\gamma}_nK_{2n-i}(F)\otimes \Q 
\end{equation} 
We define the triangulated category of mixed 
Tate motives ${\cal D}_T(F)$ as 
the triangulated subcategory of ${\cal D}{\cal M}_F$ generated by the objects $\Q(i)$, 
$i \in \Z$, see s. 5.4 of [G8] for more details. As we mentioned above, similar results are 
provided by M.Levine's theory [L]. 

Here is a useful specific way to get objects in the category ${\cal D}{\cal M}_F$ 
which will be used below. 
For a set $X$ denote by $\Z[X]$ the free abelian group generated by the elements of  $X$. 
Let ${\cal C}$ be a category. We define a new category $\Z[{\cal C}]$ whose 
objects are the objects of ${\cal C}$, and 
$$
{\rm Mor}_{\Z[{\cal C}]}(X, Y):= \Z[{\rm Mor}_{{\cal C}}(X, Y)]
$$  
A complex 
$X_1 \stackrel{d}{\lra} X_2 \stackrel{d}{\lra} ... \stackrel{d}{\lra} X_n$ of objects 
of the category ${\cal C}$ is defined by the morphisms 
 $d: X_k \lra X_{k+1}$ in $\Z[{\cal C}]$ such that $d \circ d =0$. 

Let $Y_{\bullet}$ be a complex of topological spaces. 
Applying the singular chain complex functor 
$Y \lra S_{*}(Y)$ to each 
member of the complex we get a bicomplex. The Betti homology  
$H_n^B(Y_{\bullet})$ are 
the homology of the total complex of this bicomplex. 
One can  generalize  all this 
to the case of  algebraic varieties, and 
define  the De Rham, \'etale,  or other  
cohomology of a complex of varieties.

{\bf 4. The abelian category of mixed Tate
  motives over a number field}. 
Let  $F$ be a  number field. Then it follows from the results [Bo], [Be1] that 
 the Beilinson and Borel regulators 
coincides up to a non zero rational factor. This implies 
the Beilinson-Soul\'e vanishing conjecture for number fields. One has 
\begin{equation} \label{vvvvvtt}
gr^{\gamma}_nK_{2n-1}(F)\otimes \Q  = K_{2n-1}(F)\otimes \Q; 
\quad K_{2n}(F)\otimes \Q = 0 \quad n>0, \quad K_0(F) = \Z
 \end{equation} 
Therefore, as explained in ch. 5 of [G8], or in [L1], 
one can define 
  the abelian category  ${\cal M}_T(F)$ 
of mixed Tate
  motives over  $F$ with  all the needed
  properties. In particular it is a Tate category and its derived category is equivalent to the 
triangulated subcategory of ${\cal D}_T(F)$. 
Combining this with (\ref{vvvvv})  and (\ref{vvvvvtt}) we get  
\begin{equation} \label{8-23/99}
Ext^1_{{\cal M}_T(F)}(\Q(0),\Q(n)) = K_{2n-1}(F) \otimes \Q, 
\end{equation}
\begin{equation}\label{8-23.1/99}
Ext^i_{{\cal M}_T(F)}(\Q(0),\Q(n)) = 0 \quad \mbox{for $ i>1$}
\end{equation} 
This means that the fundamental Lie algebra of the category ${\cal M}_T(F)$, which we denote by 
 ${L}_{\bullet}(F)$, is isomorphic to 
a free graded pro-Lie algebra over $\Q$ generated by the duals to 
finite dimensional $\Q$-vector spaces (computed by Borel) 
\begin{equation}\label{5-28.00}
K_{2n-1}(F) \otimes \Q \quad = \quad [\Z^{{\rm Hom}(F,\C)}\otimes\R(n-1)]^+, \quad n >1
  \end{equation} 
sitting in the degrees $-n$, 
 and the infinite dimensional $\Q$-vector space $(F^* \otimes \Q)^{\vee}$ for $n=-1$. 
The isomorphism (\ref{5-28.00}) is given by the regulator map. 

The category  ${\cal M}_T(F)$  is equipped with an array of realization functors. 
We construct in the next subsection  the Hodge and \'etale realizations. 


{\bf 5. Other examples of mixed Tate categories and realization functors}.  
Here is a general method of getting 
mixed Tate categories (cf. [BD2]). Let $F$ be a field of characteristic zero, 
${\cal C}$ be any Tannakian $F$-category and $K(1)$ be a rank one object of ${\cal C}$ such 
that the objects $K(i):= K(1)^{\otimes i}$ $i \in \Z$, are mutually nonisomorphic. 
An object $M$ of ${\cal C}$ is called a {\it mixed Tate object} if it admits a finite 
increasing filtration $W_{\bullet}$ indexed by $2\Z$, such that ${\rm Gr}^W_{2k}M$ is a 
direct sum of copies of $K(k)$. Denote by ${\cal T}{\cal C}$ the full subcategory of 
mixed Tate objects in ${\cal C}$. Then 
${\cal T}{\cal C}$ is a Tannakian subcategory in ${\cal C}$ and $({\cal T}{\cal C}, K(1))$ 
is a mixed Tate category. 

Further, if $({\cal C}_i, K(1)_i)$ are as above and $\varphi: {\cal C}_1 \lra {\cal C}_2$ is a 
tensor functor equipped with an isomorphism $\varphi(K(1)_1) = K(1)_2$, then 
$\varphi({\cal T}{\cal C}_1) \subset  {\cal T}{\cal C}_2$ and the restriction of $\varphi$
on ${\cal T}{\cal C}_1$ is a pure functor. 

We use this set up below to construct the realization functors. 

1. {\it The category of $\Q$-rational Hodge-Tate structures}. Applying the above construction 
to the category $MHS/\Q$ of mixed Hodge structures over $\Q$ we get the mixed Tate category 
${\cal H}_T$ of $ \Bbb Q$-rational Hodge-Tate structures. Namely, 
a $ \Bbb Q$-rational Hodge-Tate structure is a mixed Hodge structure over $\Bbb Q$ 
with $h^{p,q} =0$ if $p \not = q$.  
Equivalently, a Hodge-Tate structure is a $ \Bbb Q$-rational mixed Hodge structure with 
weight $-2k$ quotients isomorphic to 
a direct sum of copies of $\Q(k)$. 
Set ${\cal H}_{\bullet}:= {\cal A}_{\bullet}({\cal H}_{T})$. Then 
\begin{equation} \label{8-24/99}
{\cal H}_1 = Ext^1_{{\cal H}_{T}}(\Q(0),\Q(1))  = Ext^1_{MHS/\Q}(\Q(0),\Q(1)) = \frac{\C}{2 \pi i \Q} = \C^*\otimes {\Q}
\end{equation}

{\bf Example}. The extension  corresponding 
via isomorphism  (\ref{8-24/99}) to a given $z \in \C^*$ is provided by the 
$\Q$-rational Hodge-Tate structure 
$H(z)$, also denoted $\widetilde \log (z)$,  defined by the period matrix 
$$
 \left (\matrix{1&0 \cr \log( z)& 2 \pi i  \cr }\right ) 
$$ 
Let me recall its construction. 
Denote by  $H_{\C}$ the two dimensional $\C$-vector space with a basis 
$e_0, e_{-1}$. 
The $\Q$-vector space $H(z)_{\Q}$ is the subspace generated by 
the columns of the matrix, i.e. the vectors 
$e_0 + \log z \cdot e_{-1}$ and $2 \pi i e_{-1}$. 
The weight filtration on $H(z)_{\Q}$ is given by 
$$
W_0H(z)_{\Q} = H(z)_{\Q}, 
\quad W_{-1}H(z)_{\Q} = W_{-2}H(z)_{\Q} = 2 \pi i \cdot e_{-1}, \quad 
W_{-3}H(z)_{\Q} = 0
$$
The Hodge filtration 
is defined by 
$
F^{1}H_{\C} = 0, \quad F^0H_{\C} = <e_0>, 
\quad F^{-1}H_{\C} = H_{\C}
$.

It is easy to see that
$$
Ext^1_{{\cal H}_{T}}(\Q(0),\Q(n)) = \frac{\C}{(2 \pi i)^n \Q} = \C^*(n-1)\otimes {\Q}
\quad \mbox{for $n > 0$}
$$
It was proved by Beilinson  that the higher Ext groups are zero. 
So the fundamental Lie algebra of the category of 
Hodge-Tate structures is isomorphic to a  free graded pro-Lie 
algebra over $\Q$ generated by the 
$\Q$-vector spaces $(\C^*(n-1)\otimes {\Q})^{\vee}$ sitting in the degree $-n$ 
where $n \geq 1$. 

2. {\it The universal Hodge realization functor for the category of mixed Tate motives over a 
number field}. 
We need the following general construction. Let $\{{\cal T}_i\}$, $i \in I$, be a finite 
collection of mixed Tate categories. Then there exists  the universal 
mixed Tate category ${\cal F}_{I}{\cal T}$ equipped with pure functors 
$p_i: {\cal F}_{I}{\cal T} \lra {\cal T}_i$, the coproduct of the mixed Tate categories 
$\{{\cal T}_i\}$. 
 By definition an object $M$ of ${\cal F}_{I}{\cal T}$
is  a collection 
$$
(\psi(M), M_i, \varphi_i) \quad \mbox{where} \quad \psi(M) \in {\rm Vect}_{\bullet}, \quad 
M_i \in  {\cal T}_i, \quad \varphi_i: \psi(M) \stackrel{=}{\lra} \Psi_{{\cal T}_i}(M_i)
$$
where  $\varphi_i$ are isomorphisms 
of graded $K$-vector spaces. The universality property of the category ${\cal F}_{I}{\cal T}$ 
is obvious: if $({\cal M}, K(1)_{{\cal M}})$ is a mixed Tate category and we are 
given pure functors $Q_i: ({\cal M}, K(1)_{{\cal M}}) \lra ({\cal T}_i, K(1)_{i})$ 
then there exists unique functor 
$Q: ({\cal M}, K(1)_{{\cal M}}) \lra ({\cal F}_{I}{\cal T}, K(1))$ which sends 
an object $M$ of ${\cal M}$ to the triple $(\Psi_{{\cal M}}(M), Q_i(M), \varphi_i)$. 
Here the isomorphisms $\varphi_i$ are provided by the fact that for any pure functor 
$G: ({\cal M}, K(1)_{{\cal M}}) \lra ({\cal N}, K(1)_{{\cal N}})$ we have a canonical 
isomorphism $\Psi_{{\cal M}}(M) \lra \Psi_{{\cal N}}(G(M))$.

It follows that  the fundamental Lie algebra 
$L_{\bullet}({\cal F}_{I}{\cal T})$ is the coproduct of $L_{\bullet}({\cal T}_i)$, i.e. a graded 
Lie algebra ${\cal F}_{I}L_{\bullet}({\cal T}_i)$ freely generated by $L_{\bullet}({\cal T}_i)$. 

The realization functors from the triangulated 
category of mixed motives over $\C$ has been constructed by M. Levine ([L]) 
for his category, 
and by A. Huber [Hu] for the Voevodsky triangulated category. 
They include the Hodge and \'etale realizations.

It follows from this that for any embedding $\sigma: F \hookrightarrow \C$ 
there exists the Hodge realization functor, given by 
a pure functor between the mixed Tate categories 
$$
{\rm Hod}_{\sigma}: {\cal M}_T(F) \longrightarrow {\cal H}_T
$$

 Let $\Sigma$ be the set of all complex 
embeddings of the number field $F$ in  $\C$. 
By the universality property the 
 collection of the Hodge realizations ${\rm Hod}_{\sigma}$  provides  pure functor 
$$
{\rm Hod}: {\cal M}_T(F) \longrightarrow {\cal F}_{\Sigma}{\cal H}_T
$$
between the Tate categories 
which we call the {\it universal Hodge realization functor}. 

Let ${\cal F}_{\Sigma}{\cal H}_{\bullet }$ be the coproduct of the Hopf algebras 
${\cal H}_{\bullet}$ indexed by the set $\Sigma$, considered 
in the category of commutative pro-group schemes, 
i.e. ${\rm Hom}(A_1, A_2) := {\rm Hom}_{{\rm Hopf}-{\rm  alg}}(A_2^{\vee}, A_1^{\vee})$. 
Then ${\cal F}_{\Sigma}{\cal H}_{\bullet }$ is the fundamental commutative Hopf algebra 
for the mixed category ${\cal F}_{\Sigma}{\cal H}_T$. Therefore the functor 
${\rm Hod}$
  can be described by a morphism of the corresponding commutative Hopf algebras 
\begin{equation} \label{1.1.01.1}
h_F: {\cal A}_{\bullet}(F) \lra {\cal F}_{\Sigma} {\cal H}_{\bullet }
\end{equation}

\begin{lemma} \label{1.1.01.2} Let $F$ be a number field. Then 
the Hopf algebra map (\ref{1.1.01.1}) is injective. 
\end{lemma}

{\bf Proof}. 
By the Borel theorem [Bo] 
the regulator map on $K_{2n-1}(F)\otimes \Q$ is injective. Thus 
$\oplus_{\sigma \in \Sigma} {\rm Hod}_{\sigma}$ induces an injective map 
\begin{equation} \label{BeBo1}
  Ext_{{\cal M}_T(F)}^1(\Q(0), 
\Q(n))  \hookrightarrow \oplus_{\sigma \in \Sigma} 
Ext_{{\cal H}_T}^1(\Q(0), \Q(n))
\end{equation}
This implies that the natural morphism of the Lie algebras
\begin{equation} \label{BB32}
{\cal F}_{\Sigma} L_{\bullet }({\cal H}_T) \lra L_{\bullet}(F)
\end{equation}
is surjective. Indeed, 
$
H_1{\cal F}_{\Sigma}  L_{\bullet }({\cal H}_T) = \oplus_{\sigma \in \Sigma} H_1L_{\bullet }({\cal H}_T)
$, and  
the map on $H_1$'s induced by morphism (\ref{BB32}) is dual to  map 
(\ref{BeBo1}), and so it is surjective. 
The dual to  ${\cal F}_{\Sigma}{\cal H}_{\bullet }$ is the universal 
enveloping algebra of ${\cal F}_{\Sigma}L_{\bullet}({\cal H}_{T})$. 
So the lemma follows from surjectivity 
of map (\ref{BB32}).


3. {\it The mixed Tate category of lisse $\Q_l$-sheaves on a scheme $X$}. 
It was considered by Beilinson and Deligne in [BD2],  see an exposition of their approach 
in  s. 3.7 of [G4]. 
When $X={\rm Spec} {\cal O}_{F,S}$, where $F$ is a number field,  such a category 
was considered by Hain and Matsumoto in [HM] from a different but equivalent point of view.

Let $X$ be a connected coherent scheme over $\Z_{(l)}$ such that ${\cal O}^*(X)$ does not 
contain  all roots of unity of order equal to a power of $l$. Denote by 
${\cal F}_{\Q_l}(X)$ the Tannakian category of 
lisse $\Q_l$-sheaves on $X$.   
There is  the Tate sheaf $\Q_l(1):= \Q_l(1)_X$. 
Since $\mu_{l^{\infty}} \not \subset {\cal O}^*(X)$ the Tate sheaves $\Q_l(m)_X$ 
 are mutually nonisomorphic. So we can 
apply the general 
construction above and get the mixed Tate category ${\cal T}{\cal F}_{\Q_l}(X)$ 
of lisse $\Q_l$-sheaves on  $X$.

Now let $F$ be a number field, $S$ a finite set of prime ideals in the ring of integers 
containing 
all primes above $l$, and $X_{F,S} = {\rm Spec} {\cal O}_{F,S}$. 
 Then  
${\cal F}_{\Q_l}(X_{F,S})$ is identified with the category 
of finite dimensional $l$-adic representations of ${\rm Gal}(\overline F/F)$ 
unramified outside of ${S}$. 
We denote the corresponding mixed Tate category by 
${\cal T}^{(l)}({\cal O}_{F,S})$, 
and its fundamental Lie algebra by $L_{\bullet}^{(l)}({\cal O}_{F,S})$.  
It follows from a 
theorem of Soule [So], see s. 3.7 of [G4],  that 
the $l$-adic regulator map provides canonical isomorphism 
\begin{equation} \label{8.13.00.1}
H^i_{\cal M}(X_{F,S}, \Q(m))\otimes \Q_l := 
gr^{\gamma}_mK_{2m-i}({\cal O}_{F,S}) \otimes \Q_l 
\stackrel{=}{\lra} 
{\rm Ext}_{{\cal T}^{(l)}({\cal O}_{F,S})}^i(\Q_l(0), \Q_l(m)) 
\end{equation}
 In particular the ${\rm Ext}^i$-groups vanish for $i>1$. 
   This just means that  $L^{(l)}_{\bullet}({\cal O}_{F,S})$ is a free negatively 
graded 
Lie algebra generated by $(K_{2m-1}({\cal O}_{F,S}) \otimes \Q_l)^{\vee}$ in the degree $-m$. 

Now let ${S}':= {S} - \{l\}$. Then the $l$-adic regulator map 
\begin{equation} \label{8.13.00.111}
H^i_{\cal M}(X_{F,S'}, \Q(m))\otimes \Q_l \quad  
\stackrel{}{\lra} \quad
{\rm Ext}_{{\cal T}{\cal F}_{\Q_l}(X_{F,S})}^i(\Q_l(0), \Q_l(m)) 
\end{equation}
is an isomorphism unless $i=1, m=1$. Indeed, in these cases each of the 
groups does not change when we delete 
a closed point from the spectrum. 
If $i=1, m=1$ map (\ref{8.13.00.111}) is injective but not an isomorphism since 
 $$
H^1_{\cal M}({\cal O}_{F,S'}, \Q(1))\otimes \Q_l = {\cal O}_{F,S'}^* \otimes \Q_l 
\hookrightarrow {\rm Ext}_{{\cal T}{\cal F}_{\Q_l}(S)}^1(\Q_l(0), \Q_l(1)) = 
{\cal O}_{F,S} \otimes \Q_l
$$  Thus 
$L^{(l)}_{\bullet}({\cal O}_{F,S'})$ is a quotient of 
$L^{(l)}_{\bullet}({\cal O}_{F,S})$, and the space of 
generators of
 these
 Lie algebras differ only in the degree $-1$.

 The underlying vector space of a Galois 
representation provides a fiber functor $\widetilde \Psi$
  on the category ${\cal T}^{(l)}(X_{F,S})$. 
The image of a Galois group in a continuous  
representation  in a finite dimensional $\Q_l$-vector space  $V$
is an $l$-adic Lie group. Further, any Lie subgroup of a nilpotent $l$-adic Lie 
group is an algebraic group over $\Q_l$. It follows that  for any Galois representation 
$V$ from the category ${\cal T}^{(l)}({\cal O}_{F,S})$ 
the Lie algebra of the image of the Galois group is isomorphic to 
the image of the semidirect product of ${\rm Lie}{\Bbb G}_m$ and 
$L^{(l)}_{\bullet}({\cal O}_{F,S})$ acting 
on $\widetilde \Psi(V)$. 
Let ${\cal G}_V$  be the Lie algebra of the image of 
${\rm Gal}(\overline F/ F(\zeta_{l^{\infty}}))$ in ${\rm Aut}V$.  
Then ${\cal G}_V$ 
is isomorphic to the image of $L^{(l)}_{\bullet}({\cal O}_{F,S})$ 
 in ${\rm Der}\Psi(V)$. 

{\bf Example}.  Let ${\cal P}^{(l)}({\rm G}_m - \mu_N; v_0, v_1)$ be the $l$-adic torsor of path between the tangential base points $v_0, v_1$  at $0$ and $1$ corresponding to the canonical (up to a sign) differential  on
${\rm G}_m$, ([D]). Since ${\rm G}_m - \mu_N$ has a good reduction outside of $N$, this torsor  
 is a pro-object in the category 
${\cal T}^{(l)}(\Z[\zeta_N][\frac{1}{lN}])$. Thus the Lie algebra of the image of the 
${\rm Gal}(\overline \Q/ \Q(\zeta_{l^{\infty}N}))$ (understood as a projective limit 
of Lie algebras of its finite dimensional quotients 
provided by action of the lower central series filtration on 
 $\pi^{(l)}_1({\rm G}_m - \mu_N, v_0))$  coincides with the image 
of $L^{(l)}_{\bullet}(\Z[\zeta_N][\frac{1}{lN}])$.

Similar  results hold for the mixed Tate category of l-adic representations of 
the Galois group of $F$. 
There is the l-adic realization functor 
$$
{\cal M}_T(F) \longrightarrow {\cal T}^{(l)}(F)
$$
provided by the l-adic realization functor on the category ${\cal D}{\cal M}_F$. 
It can be described  via the morphism of the fundamental Lie algebras 
$
L_{\bullet}(F) \lra L^{(l)}_{\bullet}(F)
$. 
Thanks to (\ref{8.13.00.1}) it  leads to an  isomorphism 
$L_{\bullet}(F)\otimes \Q_l \lra L^{(l)}_{\bullet}(F)$. Thus the $l$-adic realization provides 
an equivalence of categories
$$
{\cal M}_T(F)\otimes \Q_l \stackrel{\sim }{\longrightarrow} {\cal T}^{(l)}(F)
$$

{\bf 6. The abelian category of the mixed Tate motives over the ring of integers in 
a number field}. Let ${\cal O}$ be the ring of integers in a number field $F$, 
 $S$ a finite set of prime ideals in ${\cal O}$, and 
${\cal O}_{F, S}$ the  localization of ${\cal O}$ at the primes belonging to  $S$.

\begin{definition} \label{8-24.2/99}
a) The Hopf algebra ${\cal A}_{\bullet}({\cal O}_{F, S})$ is the maximal 
Hopf subalgebra 
of ${\cal A}_{\bullet}(F)$ such that 
$$
{\cal A}_{_{1}}({\cal O}_{F, S}) \quad = \quad {\cal A}_{1}(F) \quad = \quad {\cal O}^*_S \otimes \Q
$$

b) The category of mixed Tate motives over ${\rm Spec}{\cal O}_{F, S}$ 
is   
the category of graded comodules over ${\cal A}_{\bullet}({\cal O}_{F, S})$. 
\end{definition}

Using the canonical equivalence 
$$
\Psi:  {\cal M}_T(F) \quad \lra \quad {\cal A}_{\bullet}(F)-comod
$$
we identify the category ${\cal M}_T({\cal O}_{F, S})$ with a   
subcategory of ${\cal M}_T(F)$: an object of ${\cal M}_T(F)$ belongs to ${\cal M}_T({\cal O}_{F, S})$ if and only if all its matrix coefficients 
(\ref{8-24.5/99}) are in the subalgebra ${\cal A}_{\bullet}({\cal O}_{F, S})$. 

\begin{lemma} \label{5-18.00.2}
The Hopf algebra ${\cal A}_{\bullet}({\cal O}_{F, S})$ can  
be defined inductively: 
\begin{equation} \label{5-28.00.1}
{\cal A}_{n}({\cal O}_{F, S}) = \quad \left\{x \in {\cal A}_{n}(F)\quad | \quad \Delta'(x) \in \quad \oplus_{k=1}^{n-1} {\cal A}_{k}({\cal O}_{F, S}) \otimes 
{\cal A}_{n-k}({\cal O}_{F, S}) \right\}
\end{equation}
\end{lemma}

{\bf Proof}. Follows immediately from the definitions. 

Recall the following construction. Let $A_{\bullet}$ 
be a positively graded coassociative coalgebra with counit and 
the restricted coproduct 
$\Delta': A_{\bullet} \lra A_{>0}\otimes A_{>0}$. Then 
one can  define a  map 
$$
{\Delta'}_{[n]}: A_{n} \lra \otimes^{n} A_{1}, \qquad 
$$
Informally it is obtained by dualization of the product map 
$\otimes^{n} A_{1}^{\vee} \lra A_{n}^{\vee}$. To define the product map we 
have to specify the order of multiplication, so there are $(n-1)!$ possibilities. Thanks to the associativity 
the result is independent of this. 

We can 
define ${\Delta'}_{[n]}$ as the composition
$$
A_{n} \stackrel{\Delta_{n-1,1}}{\lra} A_{n-1}\otimes A_{1} 
\stackrel{\Delta_{n-2,1} \otimes {\rm Id}}{\lra} A_{n-2} \otimes \otimes^2 A_{1} \lra ... \lra \otimes^n A_{1}
$$ 
As before, we can use another chains of arrows to define ${\Delta'}_{[n]}$. 
In fact there are $(n-1)!$ different ways to  define ${\Delta'}_{[n]}$. 
They are described as follows. First, apply the map $\Delta_{k,n-k}$. Then apply either 
$\Delta_{i,k-i} \otimes {\rm Id}$ or ${\rm Id} \otimes \Delta_{j, n-k-j}$, and so on. After the $k$ steps we are getting an element of
$A_{i_1} \otimes ... \otimes  A_{i_k}$, and we can apply any of the operators of the shape ${\rm Id} \otimes ... \otimes \Delta_{j,i_p-j}\otimes ... \otimes {\rm Id}  $ to it.  Thanks to the coassociativity 
they give the same result.  

\begin{lemma} \label{8-24.1/99}
\begin{equation} \label{8-24.10/99}
{\cal A}_{n}({\cal O}_{F, S}) = \{x \in {\cal A}_{n}(F)\quad | \quad
{\Delta'}_{[n]}(x) \in \otimes^n{\cal O}_{F, S}^*\otimes \Q\}
\end{equation}
\end{lemma}

{\bf Proof}. Suppose that $x$ belongs to the right hand side of 
(\ref{8-24.10/99}). Let us show by induction that 
$$
\Delta_{k,n-k}(x) \in 
{\cal A}_{k}({\cal O}_{F, S}) 
\otimes {\cal A}_{n-k}({\cal O}_{F, S}) \quad \mbox{for any 
$ 1 \leq k \leq n-1$ }
$$ 
We can define the map $\Delta'_{[n]}$ as a composition
$$
\Delta'_{[n]} = ({\rm Id}\otimes \Delta'_{[n-k]}) \circ  
(\Delta'_{[k]}\otimes {\rm Id}) \circ \Delta_{k,n-k}
$$
So by the induction assumption 
$$
(\Delta'_{[k]}\otimes {\rm Id}) \circ \Delta_{k,n-k}(x) \in 
\otimes^{k}{\cal O}^*_S \otimes {\cal A}_{n-k}({\cal O}_{F, S})
$$
Applying the induction assumption again  we get the lemma. 

Denote by $H_{(n)}$ the degree $n$ part of $H$. 
\begin{lemma} \label{8-24.3/99} 
Let $n>1$. Then the morphism 
\begin{equation} \label{8-24.4/99}
H_{(n)}^i\Bigl({\cal A}_{\bullet}({\cal O}_{F, S})\Bigr) \lra
H_{(n)}^i\Bigl({\cal A}_{\bullet}(F)\Bigr)
\end{equation}  
induced by the  natural inclusion of the Hopf algebras is an isomorphism. 
\end{lemma}

{\bf Proof}. The cohomology of a positively graded coalgebra $A_{\bullet}$ 
can be computed via the reduced cobar complex
$$
A_{\bullet} \stackrel{\Delta'}{\lra} A_{>0}\otimes A_{>0} 
\stackrel{\Delta' \otimes {\rm Id}   - {\rm Id} 
\otimes \Delta' }{\lra} A_{>0}\otimes A_{>0}\otimes A_{>0} \lra ... 
$$
It follows from the very definitions that (\ref{8-24.4/99}) is an isomorphism for $n>1, i=1$. Let us show that it is an isomorphism for $i=2$, i.e. 
$$
H^2\Bigl({\cal A}_{\bullet}({\cal O}_{F, S})  \lra \otimes^2{\cal A}_{\bullet}({\cal O}_{F, S}) \lra  \otimes^3{\cal A}_{\bullet}({\cal O}_{F, S}) \Bigr) = 0
$$
Since $H^2({\cal A}_{\bullet}(F))=0$, for any cycle in 
$x \in \otimes^2{\cal A}_{\bullet}({\cal O}_{F, S})$ 
there exists $y \in {\cal A}_{\bullet}({\cal O}_{F, S})$ such that i.e. $x = \Delta'y$.
Then $y \in {\cal A}_{\bullet}({\cal O}_{F, S})$ by definition. 
 Since  $H^2$ is zero the higher cohomology vanish: one can prove this using a similar fact for the negatively graded Lie algebras since our Hopf algebra is dual to the universal enveloping of such a Lie algebra. 
The lemma follows.

\begin{corollary} \label{8-24.6/99} 
\begin{equation} \label{8-24.7/99}
{\rm Ext}^i_{{\cal M}_T({\cal O}_{F, S})}\Bigl(\Q(0), \Q(n)\Bigr)  = 
gr^{\gamma}_nK_{2n-i}({\cal O}_{F, S})\otimes \Q
\end{equation}  
\end{corollary}

The right hand side of (\ref{8-24.7/99}) is known to be zero for $i>1$. 

{\bf Proof}. It follows from the previous lemma since 
$$
K_{n}({\cal O}_{F, S})\otimes \Q \quad =  \quad K_{n}(F)\otimes \Q \quad \mbox{for $n>1$}
$$ 

Therefore the only difference between the Ext's in the categories of mixed Tate motives over $F$ and ${\cal O}_{F, S}$ is the group $Ext^1(\Q(0), \Q(1))$. This group is infinite dimensional for the motives over $F$ and finite dimensional for the motives over the $S$-integers. As a result the category of mixed Tate motives over the 
$S$-integers is ``much smaller'' 
then the one over $F$. For instance  all graded components of 
${\cal A}_{\bullet}(F)$ are infinite dimensional $\Q$-vector space, while 
for ${\cal A}_{\bullet}({\cal O}_{F, S})$ they are finite dimensional.

Let $M$ be a mixed Tate motives over $F$, i.e. it is 
an object of the category ${\cal M}_{T}(F)$. We say that $M$ is defined over 
${\cal O}_{F, S}$ if it belongs to the subcategory ${\cal M}_{T}({\cal O}_{F, S})$. 

\begin{definition}
Let $(M, v_0, f_n)$ be an $n$-framing on a mixed Tate motive $M$ over $F$. We say that 
the equivalence class of the $n$-framed mixed Tate motive $(M, v_0, f_n)$ 
is defined over ${\cal O}_{F, S}$ if its image in ${\cal A}_{n}(F)$ belongs to 
 the subspace ${\cal A}_{n}({\cal O}_{F, S})$. 
\end{definition}

It follows from the definitions that $M$ is defined over 
${\cal O}_{F, S}$ if and only if  any framing on $M$ leads to a framed 
motive whose  equivalence class is  
defined over 
${\cal O}_{F, S}$.

\begin{lemma} \label{8-25/99}
Let $(M, v_0, f_n)$ be an $n$-framed motive over 
 ${\cal O}_{F, S}$. Then its minimal representative $\overline M$ is defined over ${\cal O}_{F, S}$.
\end{lemma}

Of course $M$ itself is not necessarily defined over 
 ${\cal O}_{F, S}$. 

{\bf Proof}. Let us assume the opposite. Then a certain matrix 
element of $\overline M$ does not belong to ${\cal A}_{\bullet}({\cal O}_{F, S})$. 
Since the minimal representative is a subquotient of every 
mixed Tate within the  equivalence class of $M$, {\it every} 
motive  from the  equivalence class has such a matrix coefficient (\ref{8-24.5/99}). 
Thus no representative in this  equivalence class 
is defined over ${\cal O}_{F, S}$.  

On the other hand, since $H^2({\cal A}_{\bullet}({\cal O}_{F, S}) ) =0$ there exists an $n$-framed motive 
$(M', v'_0, f'_n)$ over 
 ${\cal O}_{F, S}$ such that 
$$
\Delta((M', v'_0, f'_n)) = \Delta((M, v_0, f_n))
$$
According to (\ref{5-28.00.1}) the class of the $n$-framed object $(M', v'_0, f'_n)-(M, v_0, f_n)$ is defined over 
${\cal O}_{F, S}$. Thus the class of $(M', v'_0, f'_n)$ is also defined over ${\cal O}_{F, S}$. This contradiction proves the lemma.

So to produce a mixed Tate motive over ${\cal O}_{F, S}$ it is enough to produce 
an element of ${\cal A}_{n}({\cal O}_{F, S})$: then we take an $n$-framed mixed Tate motive representing this element and its 
minimal representative.

\begin{theorem} \label{8-25.1/99} Suppose that $S$ contains all primes of ${\cal O}_{F, S}$ above a 
given 
prime integer $l$. 
Then the   $l$-adic realization of a mixed Tate motive $M$ over $F$ 
is unramified over ${\cal O}_{F, S}$ if and only if $M$ is defined over ${\cal O}_{F, S}$. 
In particular there is a functor 
$$
{\cal M}_{T}({\cal O}_{F, S}) \lra {\cal T}^{(l)}({\cal O}_{F, S})
$$
which induces an equivalence of categories $
{\cal M}_{T}({\cal O}_{F, S})\otimes \Q_l 
\stackrel{\sim}{\lra} {\cal T}^{(l)}({\cal O}_{F, S})
$. 

Thus one has $L_{\bullet}^{(l)}({\cal O}_{F,S}) = 
L_{\bullet}({\cal O}_{F,S})\otimes \Q_l$. 
\end{theorem}

{\bf Proof}. The inertia group for a prime ideal ${\cal P}$ over a a prime $p$ 
has  a pro-$p$ subgroup of wild ramification, 
and the quotient is isomorphic to $\widehat \Z(1)$. So it acts on $l$-adic representations via 
its $\Z_l(1)$ quotient, which shifts the weight filtration down by $2$.  
So the statement boils down to the isomorphism $L_{-1}^{(l)}({\cal O}_{F,S}) = 
L_{-1}({\cal O}_{F,S})\otimes \Q_l$, which is clear.
The theorem follows.

\section{ The motivic torsor of path}

Let $F$ be  a number field and $z_i, x,y \in P^1(F)$. 
In this chapter we construct the motivic torsor of path 
${\cal P}_{\cal M}({\Bbb A}^1 - \{z_1, ..., z_{m}\}; x,y)$. 
It is a pro-object of the abelian category 
${\cal M}_T(F)$ of mixed Tate motives over $F$.

A  detailed exposition of different  realizations of the pronilpotent completion 
${\cal P}^{\rm nil}(X; x,y)$ of the topological 
 torsor of path 
   was given by  Deligne in  [D].

We start with a construction of the motivic torsor of path $P(X; x,y)_{\bullet}$ for a 
 regular variety $X$ over an arbitrary field $F$, understood as an object of the 
triangulated category ${\cal D}{\cal M}_F$. The 
construction is 
based on  A.A. Beilinson's unpublished approach to  path torsor 
${\cal P}^{\rm nil}(X; x,y)$.  Its main advantage is that the isomorphism between 
$H_n$ of its  realization and  the corresponding 
version of ${\cal P}^{\rm nil}(X; x,y)$ is very transparent.  
In particular it gives  an  easy 
 way to 
put the mixed Hodge structure on ${\cal P}^{\rm nil}(X; x,y)$. 
We also discuss the   Betty, \'etale,   and  
De Rham  realizations.

A different  approach to the motivic torsor of path was 
developed by Wojtkowiak  [W2],  and 
Shiho [Sh].  

Let $X$ be a variety over a number field $F$ whose motive  $M(X)$ belongs to the  
derived category of mixed Tate motives ${\cal D}_T(F)$ over $F$.  For instance  
one can take  $X =  {\Bbb A}^1 - \{z_1, ..., z_{m}\}$. Then 
using the $t$-structure on 
the category ${\cal D}_T(F)$ we define  
${\cal P}_{{\cal M}}(X; x,y):= H_n^t(P(X; x,y)_{\bullet}$.

{\bf 1. Preliminaries}. 
Let $Y_0, ..., Y_n$ be subspaces of a   a topological space $Y$. 
For an ordered subset $I = \{i_1 < ... < i_k\}$ 
of $\{0, 1, ..., n\}$ set 
\begin{equation} \label{2.11.01.1}
Y_I:= Y_{i_1} \cap ... \cap Y_{i_k}; \qquad Y_{(k)}:= \cup_{|I|=k}Y_I
\end{equation}
Then there is a complex of topological spaces (see s. 4.3)
\begin{equation} \label{2.11.01.2}
Y_{(\bullet)}:= \qquad Y_{(n+1)} \stackrel{d}{\lra} ... \stackrel{d}{\lra} Y_{(2)} 
\stackrel{d}{\lra} Y_{(1)} \stackrel{d}{\lra} Y_{(0)} = Y
\end{equation}
where $Y_k$ is in degree $k$, 
the last map $Y_{(1)} \lra Y$ is provided by the maps $f_i$, 
and for $k >1$ the differential $d: Y_{(k)} \lra Y_{(k-1)}$ 
is defined by 
$$
d := \sum_{|I|=k}\sum_{p=1}^k(-1)^p d_I^p; \qquad d_I^p: Y_{i_1} \cap  ...\cap    Y_{i_k} 
\hookrightarrow Y_{i_1} \cap  ... \widehat Y_{i_p} ... \cap  Y_{i_k}
$$
The Betti homology  $H^B_*(Y_{(\bullet)}; \Z)$ 
coincide with 
$H^B_*(Y, \cup Y_i; \Z)$. Indeed, consider the simplicial space 
whose $(k-1)$-simplices are defined by formula (\ref{2.11.01.1}) where $I$ is now any non empty 
subset $i_1 \leq i_2 \leq ... \leq i_k$. The face and degeneration maps are 
given 
by the natural inclusions and  diagonals. 
The topological realization of this simplicial set maps onto the union 
$\cup Y_i$ with    contractable fibers. 
The complex $S_*(Y_{(\bullet)})$ is isomorphic to the normalization   of 
the simplicial abelian group obtained by  applying the singular complex functor $S_*$ 
to this simplicial space.

{\bf 2. A geometric construction of the pronilpotent completion of the torsor of path 
(after A.A. Beilinson) }. 
Let $X$ be a connected topological manifold and $x, y \in X$. Consider in $X^n$ the  
following ordered sequence of  $n+1$ submanifolds: 
\begin{equation} \label{82}
D^n_0:= \{x = t_1\}; \quad  D^n_i:= \{t_i = t_{i+1}\}; \quad D^n_n := \{t_n = y\}
\end{equation}
Define a complex of manifolds $\widetilde P^n(X; x,y)_{\bullet}$ as 
the  complex (\ref{2.11.01.2})  
for  this ordered collection 
of submanifolds.  For instance  
$\widetilde P^n(X; x,y)_{0} =X^n$ and $\widetilde P^n(X; x,y)_{1}$ is 
the disjoint union of the submanifolds 
(\ref{82}). 

Observe that $\widetilde P^n(X; x,y)_{n+1} $ is the 
intersection of the submanifolds (\ref{82}),  so it is 
empty if $x \not = y$ and it is 
a single point if $x=y$. 
We define the  complex  $P^n(X; x,y)_{\bullet}$ by 
deleting $\widetilde P^n(X; x,y)_{n+1} $. It coincides with $ \widetilde 
P^n(X; x,y)_{\bullet}$ when $x \not = y$. If $x=y$  the natural map 
\begin{equation} \label{2.11.01.3}
H_*(P^n(X; x,x)_{\bullet}; \Z)  \lra H_*(\widetilde P^n(X; x,x)_{\bullet}; \Z)
\end{equation}
is an isomorphism for $* \not = n$, and for $*  = n$ it is surjective 
with the  kernel  isomorphic to $\Z$. 

Similarly one defines a complex of varieties $P^n(X)_{\bullet}$ 
together with the natural projection
\begin{equation} \label{2.16.01.6}
\pi: P^n(X)_{\bullet} \lra X \times X\qquad \mbox{such that} \quad \pi^{-1}(x,y) = 
P^n(X; x,y)_{\bullet}
\end{equation}
Namely, we consider in $X \times X^n \times X = \{x, t_1, ..., t_m, y\}$ the collection of 
subvarieties given by (\ref{82}) and proceed just as we did before.

There is a natural map 
\begin{equation} \label{81112}
\Delta: {\cal P}(X;x,y) \lra H_n(P^n(X; x,y)_{\bullet}; \Z)
\end{equation}
defined as follows. Let $\gamma:[0,1] \lra X$ be a path  from  $x$ to $y$ and 
$$
\Delta_n:= \{(t_1, ..., t_n) \in \R^n| 0 \leq t_1 \leq ... \leq t_n \leq 1\}
$$
 the standard $n$-simplex. Then there is a simplex $
\Delta_{\gamma} := \{\gamma(t_{1}),... , \gamma(t_{n})\} \subset X^n$. 
 For a subset $I \subset \{0, 1, ..., n\}$ let $X_I:= \cap_{p \in I} D^n_p$,  
and $\Delta^I_{\gamma}$ be the corresponding face of the simplex $\Delta_{\gamma}$. 
 The map (\ref{81112})  
assigns to a path $\gamma$ the  cycle  in the singular
bicomplex  of $P^n(X; x,y)_{\bullet}$ whose component in $S_{n-|I|}(X_I)$ is 
 $\Delta^I_{\gamma}$. 
The homology class of this cycle depends only 
on the homotopy class of $\gamma$, so we get the map (\ref{81112}). 
The kernel of  map (\ref{81112})
is  generated by the  
constant path at $x=y$.  

  Denote by ${\cal P}(X; x,y)$ the free abelian group 
generated by the homotopy 
equivalence classes of path from $x$ to $y$ on $X$. It is a left torsor over the 
group algebra $\Z\pi_1(X; x)$ of the fundamental group of $X$ based at $x$, and the right 
torsor over $\Z\pi_1(X; y)$. 
Let ${\cal I}_y$ be the augmentation ideal of $\Z\pi_1(X; y)$, and  ${\cal I}_y^k$ 
is its $k$-th power. 
Denote by ${\cal I}_{x,y}^k$ the image of the natural map
$ \Z\pi_1(X; x,y) \otimes {\cal I}_y^k \lra {\cal P}(X; x,y)$. Set
$$
{\cal P}^n(X; x,y):= {\cal P}(X; x,y)/ {\cal I}_{x,y}^{n+1}; \qquad {\cal P}^{\rm nil}(X; x,y):=
\lim_{\longleftarrow}{\cal P}^n(X; x,y)
$$
The following result was proved long ago by A.A. Beilinson (unpublished).

\begin{theorem} \label{2.2.01.1} One has 
\begin{equation} \label{2.8.01.1}
H_i(P^n(X; x,y)_{\bullet}; \Z) = \left\{ \begin{array}{ll}
{\cal P}^{n}(X; x,y) & i=n\\ 
0 & i< n \end{array}\right.
\end{equation}
\end{theorem}

{\bf Proof}. We will prove the theorem by induction on $n$. 
The exact sequence $$
0 \lra H_1(X;\Z) \lra H_1(X, \{x, y\};\Z) \lra \Z \lra 0
$$
makes clear all the statements of the theorem for $n=1$. 

Consider the collection  of all the submanifolds (\ref{82}) but the first one $x=t_1$. 
It has a natural ordering. Denote  by 
 $T_{y, \bullet}^n$ 
the complex of varieties corresponding to it via 
the general construction (\ref{2.11.01.2}). 
  Let   $\widehat T_y^n$  be the union of all the  
submanifolds (\ref{82}) except the first one. Then  
$H_*(T_{y, \bullet}^n; \Z) = H_*(X^n,  \widehat T_y^n; \Z)$. 

Let   $T^n_{x,y}$ be the union of submanifolds (\ref{82}).   
Observe that the pair $(p^{-1}(t), p^{-1}(t) \cap \widehat T_y^n)$
is identified with the pair $(X^{n-1}, T^{n-1}_{t, y})$. 
So for an arbitrary pair $\{x,y\}$ there is an exact sequence of complexes of varieties
(which should be treated as an exact triangle in ${\cal D}{\cal M}_F$):
$$
0 \lra T_{y, \bullet}^n \lra {P}^n(X; x,y)_{\bullet} \lra {P}^{n-1}(X; x,y)_{\bullet}[-1] \lra 0
$$
  where $H_i(C_{\bullet}[-1]):= H_{i-1}(C_{\bullet})$. It leads to a long exact sequence, 
which for $x \not = y$ is simply the long exact sequence for the relative homology of the 
pair $(X^n,  \widehat T_y^n)$ modulo the fiber of the projection $p$ over $x$, given by the 
pair $(p^{-1}(x), p^{-1}(x) \cap \widehat T_y^n)$:
\begin{equation} \label{881*}
... \lra H_n(P^{n-1}(X; x,y)_{\bullet}; \Z) \stackrel{f}{\lra} 
H_n(X^n,  \widehat T_y^n; \Z) \lra H_n(P^{n}(X; x,y)_{\bullet}; \Z) \lra
\end{equation}
$$
 H_{n-1}(P^{n-1}(X; x,y)_{\bullet}; \Z) \lra 
H_{n-1}(X^n,  \widehat T_y^n; \Z) \lra 
H_{n-1}(P^{n}(X; x,y)_{\bullet}; \Z) \lra ...
$$
By the induction assumption 
one has
\begin{equation} \label{91}
H_{i}(P^{n-1}(X; x,y)_{\bullet}; \Z) = \left\{ \begin{array}{ll}
{\cal P}^{n-1}(X; x,y) & i=n-1\\ 
0 & i< n-1 \end{array}\right.
\end{equation}

\begin{lemma} \label{2.2.01.3}
$H_{i}(X^n,  \widehat T_y^n; \Z) =0$ for $i<n$. One has 
$H_{n}(X^n,  \widehat T_y^n; \Z)/{\rm Im}(f) = {\cal I}_{y}^n/{\cal I}_{y}^{n+1}$. 
\end{lemma}

{\bf Proof}. Consider the projection $p:X^n \lra X$ onto the first factor. 

\begin{center}
\hspace{4.0cm}
\epsffile{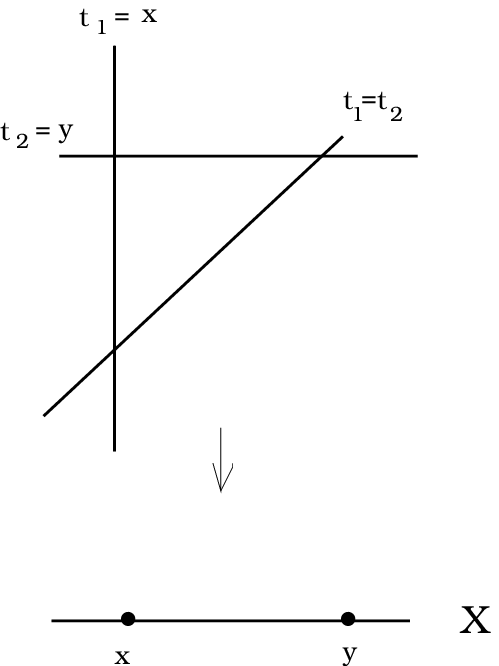}
\end{center}

We  compute $H_{i}(X^n,  \widehat T_y^n; \Z) $ for $i \leq n$ by induction 
on $n$ using the homological Leray spectral sequence for the projection $p$. Its $E^2$-term 
 looks as follows:
$$
E^2_{p,q} = H_p(X, {\cal R}_q) \quad => \quad H_{p+q}(X^n, \widehat T_y^n; \Z)
$$
Here ${\cal R}_q$ is the co sheaf associated with the precosheaf 
$U \lra H_q(p^{-1}(U), p^{-1}(U) \cap \widehat T^n_y; \Z)$. Let us show that 
${\cal R}_q = 0$ for $q<n-1$ and calculate ${\cal R}_{n-1}$. 
Observe that $T^n_{y, n}$ is a single point $(y, ..., y) \subset X^{n}$. Let 
$\tau_{\leq n-1}T^n_{y, \bullet}$ be the stupid  truncation of the  complex of varieties 
$T^n_{y, \bullet}$, i.e. 
\begin{equation} \label{2.13.01.1}
T^n_{y, \bullet} \quad = \quad T^n_{y, n} \stackrel{d}{\lra} \tau_{\leq n-1}T^n_{y, \bullet}
\end{equation}
Let ${\cal L}$ be the local system on $X$ whose fiber ${\cal L}_t$ at $t \in X$ 
is ${\cal P}^{n-1}(X; t,y)$. Then by the induction assumption
$$
R_ip_* (\tau_{\leq n-1}T^n_{y, \bullet}) = \left\{ \begin{array}{ll}
{\cal L} & i=n-1\\ 
0 & i< n-1 \end{array}\right.
$$
It follows from (\ref{2.13.01.1}) that ${\cal R}_q = 0$ for $q<n-1$ and there is an exact 
sequence:


\begin{equation} \label{912}
0  \lra H_1(X, {\cal L}) \lra H_1(X, {\cal R}_{n-1}) \lra 
  \Z \stackrel{g}{\lra} 
H_0(X, {\cal L}) \lra H_0(X, {\cal R}_{n-1})\lra 0
\end{equation}
For any local system ${\cal L}$ on a manifold $X$ with $\Gamma:= \pi_1(X,y)$ one has 
\begin{equation} \label{912*}
H_0(X, {\cal L}) = H_0(\Gamma, {\cal L}_y) = {\cal L}_y/{\cal I}{\cal L}_y; \qquad 
H_1(Y, {\cal L}) = H_1(\Gamma, {\cal L}_y)
\end{equation}
The natural map $g: \Z  \lra H_0(X, {\cal L}) = \Z[\Gamma]/{\cal I} = \Z$ 
is an isomorphism. It follows that $H_0(X, {\cal R}_{n-1})= 0$, and thus 
$H_{n-1}(X^n,  \widehat T_y^n; \Z) =0$. 
Moreover, $H_1(X, {\cal L}) = H_1(X, {\cal R}_{n-1})$ by (\ref{912}). 
Using the exact sequence of $\Gamma$-modules 
$$
0 \lra {\cal I}_y^n \lra \Z[\Gamma] \lra {\cal L}_y \lra  0 
$$
we get $H_1(X, {\cal L}) = {\cal I}_y^n/{\cal I}_y^{n+1}$. Finally, 
the image of the subgroup $H_{0}(X,  {\cal R}_n)$ in $H_n(X^n,  \widehat T_y^n; \Z)$
coincides with ${\rm Im}(f)$. Combining these statements we get the second part of the lemma. 
The lemma is proved. 

The $i<n$ part  of (\ref{2.8.01.1}) 
follows from the first statement of the lemma and (\ref{91}).
 
Observe the exact sequence
\begin{equation} \label{2.13.01.11}
0 \lra \frac{{\cal I}^n_{x,y}}{{\cal I}^{n+1}_{x,y}} \lra {\cal P}^n(X; x,y) 
\lra {\cal P}^{n-1}(X; x,y)
 \lra  0 
\end{equation}
It follows from the lemma that there is an exact sequence
\begin{equation} \label{2.13.01.10}
0 \lra \frac{{\cal I}^n_{y}}{{\cal I}^{n+1}_{y}} \lra  H_{n}(P^n(X; x,y)_{\bullet}; \Z)
\lra H_{n-1}(P^{n-1}(X; x,y)_{\bullet}; \Z)
 \lra  0 
\end{equation}
The map $\Delta$ provides a map of the first exact sequence to the second one. 
The induced map on the groups on the right is given by the isomorphism (\ref{91}). 
Recall  the canonical isomorphism 
${\cal I}^n_{y}/{\cal I}^{n+1}_{y} \lra {\cal I}^n_{x,y}/{\cal I}^{n+1}_{x,y}$ 
given by the right action of the fundamental group on the torsor of path. 
It is easy to see that it also induces the isomorphism of the subgroups on the left. 
So $\Delta$ 
is an isomorphism. The $i<n$ part  of (\ref{2.8.01.1}) is proved. 
The theorem is proved.

{\it The $l$-adic realization}. Suppose that $X$ is a regular
 algebraic  variety over a field $F$ 
and $x, y \in X(F)$. 
Denote by ${\cal P}_{(l)}(X; x,y)$ the $l$-adic torsor of path ([D]) 
with coefficients in $\Q_l$. A point $x \in X(F)$ provides a fiber functor $F_x$ 
on the category of $\Q_l$-local systems on $X$. Elements of 
 ${\cal P}_{(l)}`(X; x,y)$  act as  the natural transformations between the 
fiber functors $F_x$ and $F_y$.   It is a ${\rm Gal}(\overline F/F)$-module. 
Similarly to the classical case there are the finite dimensional 
quotients ${\cal P}^n_{(l)}(X; x,y)$. Set $\overline X:= X\otimes_F \overline F$. 
 The \'etale homology are dual to the \'etale cohomology. 

\begin{theorem} \label{2.4.01.1} Let $X$ be a regular variety over $F$. Then 
there is an isomorphism of ${\rm Gal}(\overline F/F)$-modules
$$
{\cal P}^n_{(l)}(X; x,y) = 
H_n^{et}(P^n_{\cal M}(\overline X; x,y)_{\bullet}; \Q_l)
$$
\end{theorem}

{\bf Proof}. There is the $l$-adic version of the exact sequence (\ref{2.13.01.11}), 
and computing the $l$-adic homology $H_i^{et}(P^n(\overline X; x,y)_{\bullet}; \Q_l)$ 
just the same way we did for the Betty homology during the 
proof of theorem \ref{2.2.01.1} 
we get the $l$-adic version of the exact sequence (\ref{2.13.01.10}). To complete 
the proof we need a natural map of the first to the other. Observe that 
$H_n^{et}(P^n(\overline X; x,y)_{\bullet}; \Q_l)$ is a fiber over $(x,y)$ 
of an $l$-adic  local system on on the \'etale site of 
$X \times X$ given by $L_n\pi_*(P^n(X)_{\bullet})$, see (\ref{2.16.01.6}). 
There is a distinguished element 
$[1] \in H_n^{et}(P^n(\overline X; x,x)_{\bullet}; \Q_l)$ 
corresponding to the identity loop. It is stable by the action of the Galois group. 
Acting on $[1]$ by the elements of ${\cal P}_{(l)}(X; x,y)$ 
 we get a morphism of Galois modules $[1] \otimes {\cal P}_{(l)}(X; x,y) \lra  
H_n^{et}(P^n(\overline X; x,y)_{\bullet}; \Q_l)$. 
The theorem follows. 

{\it The De Rham realization}. Repeating the arguments above in the  De Rham setting we get

\begin{theorem} \label{2.13.01.15} Let $X$ be a regular variety over a field $F$ of 
characteristic zero. Then 
$$
{\cal P}^n_{DR}(X; x,y) = 
H_n^{DR}(P^n_{\cal M}(X; x,y)_{\bullet})
$$
\end{theorem}
Indeed, the only thing one needs to adjust in the proof of theorem \ref{2.2.01.1} is the 
De Rham versions of (\ref{912*}) and  formula 
$H_1(X, {\cal L}) = {\cal I}_y^n/{\cal I}_y^{n+1}$. 

However to compare the De Rham realization with the bar construction it is more natural to use 
a different model for the motivic torsor of path recalled in the subsection  4 below. 

{\it The mixed Hodge structure on ${\cal P}^{\rm nil}(X; x,y)$}. 
Now let $X$ be a smooth 
complex algebraic variety. Taking $H^n$ of the 
Hodge realization of the 
mixed motive $P^n_{\cal M}(X; x,y)_{\bullet}$ we get a 
 mixed Hodge structure on ${\cal P}^n_{\cal M}(X; x,y)$. 
One can show that it is canonically isomorphic to the  mixed Hodge structure
defined in the works of Morgan [M], Hain [H], [H1] and  
Hain-Zucker [HZ]. 

{\bf 3. The mixed motive corresponding to the pronilpotent completion of the torsor of path 
${\cal P}(X; x,y)$}. Let $X$ be a regular algebraic variety over an arbitrary field 
$F$ of characteristic zero and $x,y \in X(F)$. Repeating the construction of the complex 
$P^n(X; x,y)_{\bullet}$ in the category of algebraic varieties we get a complex 
of varieties 
$$
P^n_{\cal M}(X; x,y)_{\bullet} \in {\rm Ob}{\cal D}{\cal M}_F
$$
representing an object of the triangulated category of motives 
${\cal D}{\cal M}_F$. 
The boundary map in (\ref{881*}) provides a morphism $p_n: 
{\cal P}^n_{\cal M}(X; x,y) \lra 
{\cal P}^{n-1}_{\cal M}(X; x,y)$. We define 
$$
{\cal P}_{\cal M}(X; x,y)_{\bullet}:= \quad \lim_{\longleftarrow}
{\cal P}^n_{\cal M}(X; x,y)_{\bullet}
$$
The maps $p_n$ are obviously compatible with their counterparts in all 
the realizations, so the theorems proved in s. 4.2 imply that $H_n$ of the 
Betty, \'etale, De Rham and Hodge realizations 
of ${\cal P}_{\cal M}(X; x,y)_{\bullet}$ coincide with the corresponding 
realizations of ${\cal P}^{\rm nil}(X; x,y)$ considered in [D].

{\bf 4. The torsors of path on punctured $P^1$ as mixed Tate motives}. 
 Now suppose that $F$ is a number field, $z_i \in F$ and 
$z_i \not = z_0, z_{n+1}$ for $1 \leq i \leq n$. 
Denote by $M(X)$ the object of ${\cal D}{\cal M}(F)$ corresponding to a variety $X$. According to theorem 4.1.11 in [V] one has 
$M(P^n) = \oplus_{i=0}^n \Q(i)[2i]$. It follows immediately from this that 
${P}^{m}_{\cal M}({\Bbb A}^1 - \{z_1, ..., z_n\}; z_0, z_{n+1})$ is 
an object of the triangulated category  ${\cal D}_T(F)$
of mixed Tate motives over $F$. Recall  
the $t$-structure $t$ on the triangulated category ${\cal D}_T(F)$ (see  ch.
 5 of [G8]). Its 
abelian heart 
is the category ${\cal M}_T(F)$ of mixed Tate motives over $F$. Set 
$$
{\cal P}^{m}_{\cal M}({\Bbb A}^1 - \{z_1, ..., z_n\}; z_0, z_{n+1}):= \quad 
H^t_{m}\Bigl( {P}^{n}_{\cal M}({\Bbb A}^1 - \{z_1, ..., z_n\}; z_0, z_{n+1})_{\bullet}\Bigr) \in {\cal M}_T(F)
$$ 
Here $H_*^t$ are the homology 
with respect to the $t$-structure $t$ on ${\cal D}_T(F)$.

The boundary map in   (\ref{91}) provides a morphism  
in the category ${\cal M}_T(F)$
$$
{\cal P}^{m}_{\cal M}({\Bbb A}^1 - \{z_1, ..., z_n\}; z_0, z_{n+1}) \lra 
{\cal P}^{m-1}_{\cal M}({\Bbb A}^1 - \{z_1, ..., z_n\}; z_0, z_{n+1})
$$
Denote by ${\cal P}_{\cal M}({\Bbb A}^1 - \{z_1, ..., z_n\}; z_0, z_{n+1})$ the projective limit 
with respect to these maps. 

 Observe that $H^t_1({\Bbb A}^1 - \{z_1, ..., z_n\})$ is a direct sum of $m$ copies of 
$\Q(1)$. 

\begin{lemma} \label{1012}
One has 
\begin{equation} \label{qw4021}
{\rm gr}^W_{-2m}{\cal P}_{\cal M}({\Bbb A}^1 - \{z_1, ..., z_n\}; z_0, z_{n+1}) = 
\otimes^mH^t_1({\Bbb A}^1 - \{z_1, ..., z_n\})
\end{equation}
\end{lemma}

{\bf Proof}. We can show this using the Hodge realization functor ${\rm Hod}$, 
since the corresponding statement for the Hodge realization is well known (and obvious). 
Notice that ${\rm Hod}$ is an exact pure functor between 
the corresponding Tate categories.  
The lemma is proved.

The forms $d\log (t-z_{1}), ..., d\log (t-z_{n})$ 
form a basis in $H_{DR}^1({\Bbb A}^1 - \{z_1, ..., z_n\})$. 
  Thus the forms 
\begin{equation} \label{4021}
 \quad \frac{dt_1}{t_1-z_{i_1}} \otimes  ... \otimes \frac{dt_n}{t_n-z_{i_m}} 
\end{equation}
when $\{z_{i_1}, ..., z_{i_m}\}$  run through all $m$-tuples of the points 
$z_1, ..., z_{n}$ 
  provides a basis in the dual to (\ref{qw4021}). Observe the natural  isomorphism
\begin{equation} \label{40221}
{\rm gr}^W_{0}{\cal P}_{\cal M}({\Bbb A}^1 - \{z_1, ..., z_n\}; z_0, z_{n+1}) = 
\Q(0)
\end{equation}

\begin{definition} \label{16.2.01.1} Let $F$ be a number field, $z_i \in F$, and 
$z_i \not = z_0, z_{n+1}$ for $1 \leq i \leq n$. 
The framed mixed Tate motive ${\rm I}^{\cal M}(z_0; z_{i_1}, ... , z_{i_m}; z_{n+1})$ 
is defined as the mixed Tate motive 
$$
{\cal P}_{\cal M}({\Bbb A}^1 - \{z_1, ..., z_n\}; z_0, z_{n+1})
$$
 with the framing provided  by 
the form (\ref{4021}) and the isomorphism (\ref{40221}). 
\end{definition}

\section {The  multiple polylogarithm variations of Hodge-Tate  structures }

According to Deligne  a variation of mixed Hodge structures corresponding to 
the classical $n$-logarithm 
$Li_n(x)$ is described by the matrix
\begin{equation} \label{4.20.01.1}
\left (\matrix{1&0&0&... &0 \cr 
{\rm Li}_1(x)& 2 \pi i &0 &... & 0 \cr 
{\rm Li}_2(x) & 2 \pi i \log x & (2 \pi i)^2& ... & 0\cr 
    ...  
& ...&    ... & ... &...\cr
{\rm Li}_n(x)& 2 \pi i \frac{\log^{n-1}x}{(n-1)!}&
(2 \pi i)^2 \frac{\log^{n-2}x}{(n-2)!} & ... & (2 \pi i)^n 
\cr}\right ) 
\end{equation}
whose entries are (regularized) iterated integrals along a 
certain fixed path $\gamma$ between $0$ to $x$: 
$$
{\rm Li}_n(x) = \int_0^x\frac{dt}{1-t} \circ \underbrace{\frac{dt}{t}\circ ... \circ 
\frac{dt}{t}}_{n-1 \quad  \mbox{times}}; \qquad 
\frac{\log^{n}x}{n!} =   \int_0^x \underbrace{\frac{dt}{t} \circ  ... 
\circ\frac{dt}{t}}_{n \quad  \mbox{times}} 
$$

The depth $m$ multiple polylogarithms are multivalued analytic functions on 
the configuration space $M_{0,m}$  of  $(m+2)$-tuples  of   
distinct points in $\C$ 
\begin{equation} \label{+CONF}
\{a_0; a_1, ..., a_m; a_{m+1}\}
\end{equation}
considered up to  the action of the translation group. 
The orbits can be 
parametrized by the $(m+1)$-tuples 
\begin{equation} \label{++CONF}
\{0; a_1, ..., a_m; a_{m+1}\}  
\end{equation} 
In fact multiple polylogarithms are invariant under the action of bigger group of all 
affine transformations of the plane, so 
the orbits can be parametrized by the $m$-tuples 
\begin{equation} \label{++CONF+}
\{0; a_1, ..., a_m; 1\}
\end{equation}
However to study the regularized values of the polylogarithms when 
$a_1 \to a_0$ or $a_m \to a_{m+1}$  we 
have to break the symmetry under the 
action of the multiplicative group. 

Multiple polylogarithms are iterated integrals 
on the projective line. One can consider them as functions of one variable, the 
upper limit of integration. 
Their investigation 
as functions on the configuration space (\ref{++CONF}) is a 
more natural problem. 
For the classical polylogarithms these two problems coincide. 

The goal of this section is to  show that the 
multiple polylogarithm function is a period  of variation of framed mixed 
Hodge-Tate structures over this configuration space. This variation in particular  
 provides a neat description 
of the analytic  properties of multiple polylogarithms: 
their monodromy  and differential equations. 
The explicit description of this variation is used in the next section to get 
 an explicit formulas for 
the coproduct of the 
corresponding framed Hodge-Tate structures.

The first three subsections contain mostly a well known material 
(see for instance [D], [H]) and 
serve as a background.

{\bf 1. The Hodge-Tate  structures}. We 
adopt the following version of the  definition of the  Hodge-Tate structures,  
which has been already 
defined in 
s. 3.5. 

A $\Q$-rational Hodge-Tate structure is  the following 
linear algebra data: 

i) A graded finite dimensional $\Q$-vector space $V = \oplus_{n\in \Z} V_{(n)}$.

ii) A rational subspace $V^B \subset V_{\C}:= V\otimes \C$ in the 
complexification of $V$.

iii) The compatibility condition. To spell it we define the weight 
filtration on $V_{\C}$ by 
$$
W_{2n}V_{\C} = W_{2n+1}V_{\C}:= \oplus_{k \leq n} V_{(k)}\otimes \C
$$
It induces  a weight filtration $W_{\bullet}V^B$ on the subspace $V^B$. 
Then 

a) $W_{2n}V_{\C} \cap V^B$ must define a rational structure on 
$W_{2n}V_{\C}$, and 

b) On the quotient ${\rm Gr}^W_{2n}V_{\C}$ the $\Q$-rational 
structure inherited from $V_{2n}$ is $(2\pi i)^n$ times the $\Q$-rational 
structure inherited from ${\rm Gr}^W_{2n}V^{B}$. 

A Hodge-Tate structure is equipped with the Hodge filtration $
F^pV = \oplus_{k \geq p}V_{(k)}$ on $V$. 


Below we will define a Hodge-Tate  structure by exhibiting 
the following data: 

i) Two $\Q$-rational vector spaces $V^B$ and $V^{DR}$. The Betti vector space $V^B$ 
is equipped with a weight filtration $W_{\bullet}$, 
and the de Rham vector space $V^{DR}$ has a 
grading:
\begin{equation} \label{1.1.01.5}
V^{DR} = \oplus_{n \in \Z} V^{DR}_{(n)}
\quad \mbox{ such that} \quad F^pV^{DR} = \oplus_{k \geq p}V^{DR}_{(k)}
\end{equation}  

ii) A $\Q$-linear period map 
$$
P: V^B \lra V^{DR}\otimes_{\Q}\C \qquad \mbox{such that} \quad  P^{\C}: 
V^B\otimes \C \lra V^{DR}\otimes \C \quad \mbox{is an isomorphism}
$$
Then  $(P(V^B), V^{DR})$ provides a Hodge-Tate structure. 

A morphism between the  Hodge-Tate structures defined by period maps 
$V_i^B  \stackrel{P_i}{\lra} V_i^{DR}$, $i=1,2$, is given by  a commutative diagram
$$
\begin{array}{ccc}
V_1^B & \stackrel{P_1}{\lra}& V_1^{DR}\\
&&\\
\varphi^B\downarrow && \downarrow \varphi^{DR}\\
&&\\
V_2^B & \stackrel{P_2}{\lra}& V_2^{DR}
\end{array}
$$
where $\varphi^B$ is a $\Q$-linear map preserving the weight filtration and 
$\varphi^{DR}$ 
preserves the grading.

Let $X$ be a simplicial complex algebraic variety. Then according to Deligne [D4]
there is a mixed Hodge structure in $H_n(X, \Q)$ with   
$$
V^B:= H_n^B(X(\C), \Q); \qquad 
V_{\C}^{DR}:= H_n^{DR}(X, \C):= H^n_{DR}(X, \C)^{\vee}
$$
 The pairing  $H^n_{DR}(X, \C) \otimes H_n^B(X, \Q) \lra \C$ provides the  period map   
$$
P: H_n^B(X, \Q) \lra H_n^{DR}(X, \C)
$$
 There is the weight filtration $W_{\bullet}$ on $H_n^B(X(\C), \Q)$ and the 
Hodge filtration $F^{\bullet}$ on $H_n^{DR}(X, \C)$. If the mixed Hodge structure enjoys
the condition  $h^{p,q} =0$ for $p \not =q$, 
 then we get a grading 
$$
H_n^{DR}(X, \C) = \oplus_p F^pH_n^{DR}(X, \C) \cap W_{2p}H_n^{DR}(X, \C)
$$ 
and   the De Rham $\Q$-rational structure on $H_n^{DR}(X, \C)$ compatible with this grading. 
So we get a Hodge-Tate  structure on $H_n(X(\C), \Q)$. 

{\bf 2. Unipotent variations of Hodge-Tate structures ([HZ], [H])}. 
Let $X$ be a smooth complex algebraic variety and $\overline X$ a 
compactification of $X$ such that $D:= \overline X 
- X$ is a normal crossing divisor. 

A  unipotent 
 variation of  $\Q$-Hodge-Tate structures over 
 $X$ is a local system ${\Bbb V}$ over $\Q$ on $X$ equipped with 
the following additional data:

i) {\it The weight filtration}: an  increasing filtration
$W_{\bullet}$ on the local system ${\Bbb V}$.  

ii) {\it The Hodge filtration}:  a decreasing filtration
${\cal F}^{\bullet}$ by holomorphic sub bundles on the holomorphic vector bundle 
${\cal V}:= {\Bbb V} \otimes_{\Q}{\cal O}_{X, an}$. 

iii) {\it The Griffiths transversality condition}. 
$
\nabla {\cal F}^{p} \subset 
{\cal F}^{p-1} \otimes_{{\cal O}_{X, an}} \Omega_{X, an}^1
$ for every $p$.

iv) For each point $x \in X$ the two filtrations induce the   $\Q$-rational Hodge-Tate
 structure  on the fiber $V_x$ of the local system ${\Bbb V}$ over $x$. 

v) {\it Extension to $\overline X$}.  Let $\overline {\cal V}$ be Deligne's canonical extension of ${\cal V}$ to $\overline X$. 
Then the Hodge bundles ${\cal F}^{p}$ 
can be extended holomorphically to subbundles of 
$\overline {\cal V}$. (A similar condition for the weight bundles 
${\cal W}_p:= W_p \otimes_{\Q} {\cal O}_{X, an}$ is automatically true since 
$W_p$ are flat.)

vi) {\it The Hodge-Tate condition}. ${\rm gr}^W_{m}{\Bbb V} $ is a constant variation 
of a direct sum of the Hodge-Tate structures, i.e. 
${\rm gr}^W_{2k}{\Bbb V}$ is a direct sum of copies of $\Q(-k)$ and 
${\rm gr}^W_{2k+1}{\Bbb V} = 0$. 

vii) {\it Monodromy at infinity}. 
Let $N_{D_i}$ be the logarithm of the local monodromy 
around an irreducible component $D_i$ of the divisor $D$:
\begin{equation} \label{3.10.01.3}
N_{D_i}:=  \log M_{D_i}
\end{equation} 
 Then for any $x \in X$ one has 
$
N_{D_i}(W_kV_x) \subset W_{k-2}V_x
$.

The condition vi) implies 
that the monodromy representation is unipotent. 

The following result of R. Hain ([H], theorem 8.5) 
exhibits an especially simple nature of 
the unipotent variations of Hodge-Tate structures. 

\begin{theorem} Let ${\Bbb V} \to X$ be a unipotent variations of 
Hodge-Tate structures. 
Then the canonical extension ${\overline {\cal V}} \to \overline X$  
of the underlying holomorphic vector bundle is trivial as a holomorphic 
vector bundle, so that there is a complex vector space $V$ and isomorphism 
$$
\begin{array}{ccc}
\overline {\cal V}& \lra & V \times \overline X\\
\downarrow && \downarrow \\
\overline X &= & \overline X
\end{array}
$$
of holomorphic vector bundles. Moreover, there are filtrations $F^{\bullet}$ 
and $W_{\bullet}$ of $V$ such that the Hodge and weight bundles are 
$F^p \times \overline X$ and $W_p \times \overline X$, respectively, 
under the bundle isomorphism.
\end{theorem}

{\bf 3.  The Hodge-Tate structure on the torsor of path in $\C - \{z_1, ..., z_n\}$}.
Let me spell the construction of the Hodge-Tate structure on  
${\cal P}(\C - \{z_1, ..., z_n\}; z_0, z_{n+1})$.

We start with some preliminary material. Let $V$ be a finite dimensional 
$\Q$-vector space. Denote by $L(V)$ the free Lie algebra generated by $V$. 
The universal enveloping algebra of $L(V)$ is identified with
 the tensor algebra $T(V)$ of $V$.
It is a Hopf algebra with a coproduct $\Delta$. 
Let $\widehat L(V)$ be the pronilpotent completion of the Lie algebra $L(V)$. 
Let $\widehat T(V)$ be the completion  of  
$T(V)$ with respect to powers of the augmentation ideal. Then $\widehat L(V)$ 
can be identified with the Lie algebra of primitive elements of $\widehat T(V)$, i.e. 
the elements $X$ such that $\Delta(X) = 
X \otimes 1 + 1 \otimes X$.  The set of all group like elements, i.e. the 
elements $X$ satisfying $\Delta(X) = X \otimes X$,  can be identified with 
${\rm exp} \widehat L(V) \subset \widehat T(V)$. It has a group structure.

Now we specify $V_Z = H_1(X_Z, \Q)$ where $X_Z:= \C - \{z_1, ..., z_n\}$. 
 Let $X_1, ..., X_n$ be the basis of $V_Z$ dual to 
the basis of forms $\{d\log (t- z_i)\} $  in $H_{DR}^1(X_Z)$ so that 
$$
<d\log (t- z_i), X_j> = 2\pi i \cdot \delta_{ij}
$$
Consider the 1-form 
$$
\Omega:= \quad \sum_{i=1}^n  X_i \otimes d\log (t- z_i)
$$
on  $X_{Z}(\C)$ 
with values in the Lie algebra $L(V_Z)$. 
It provides a flat connection $\Delta:= d+ \Omega$ 
in the trivialized bundle over $X_Z$ with the fiber $\widehat T(V_Z)$. The Lie algebra 
acts by multiplication from the right: $X \otimes T \lms T (-X)$.  Its horizontal section 
\begin{equation} \label{2.23.01.1}
P\int_{z_0}^{z_{n+1}} \Omega 
\end{equation}
 along a path 
$\gamma$ from $z_0$ to $z_{n+1}$ which is $1 \in \widehat T(V_Z)$ at $z_0$      
is called   
the ``ordered exponential''. 
  $\widehat T(V_Z)$ is isomorphic to 
the Hopf algebra of formal noncommutative power series $\Q<<X_1, ..., X_n>>$. 
 One has 
the Feynman-Dyson formula 
$$
P\int_{z_0}^{z_{n+1}} \Omega \quad = \quad \sum \Bigl( \int_{z_0}^{z_{n+1}} 
d\log (t- z_{i_1}) \circ ... \circ 
d\log (t- z_{i_p})\Bigr) \cdot X_{i_1} ... X_{i_p}
$$
where the sum is over all basis monomials in $\Q<<X_1, ..., X_n>>$. 
Indeed, the right hand side is a horizontal  section 
of our connection since it, as a function on $t=z_{n+1}$, 
  satisfies the differential equation 
\begin{equation} \label{4.22.01.1}
df(t) = \Bigl( \sum_{i=1}^n \frac{f(t)}{t-z_i}\cdot X_i\Bigr) dt
\end{equation} 
and its value at $z_{n+1} =z_0$ is $1$. 
One also gets  this formula by considering the right hand side  as the limit of the 
Riemann sums for the left integral.  
The shuffle product formula (\ref{schprfo}) for the iterated integrals 
of the forms $d\log (t-z_i)$  is equivalent to the fact that (\ref{2.23.01.1}) 
is a group like element. 

Formula (\ref{BT}) is equivalent to the multiplicativity property of the 
ordered exponential with respect to the composition of path: if a path $\beta$ starts at the end of the path $\alpha$ then 
\begin{equation} \label{4.25.01.1}
\Bigl(P\int_{\alpha} \Omega\Bigr) \cdot \Bigl(P\int_{\beta} \Omega\Bigr)
  = P\int_{\alpha \cdot \beta} \Omega
\end{equation}

We need below a generalization of the product formula 
(\ref{BT}). Let $\alpha_i$ be paths on $X$ such that $\alpha_{i+1}$ starts at 
the end of $\alpha_i$.  Let $\eta_i$ be smooth 1-forms on $X$. Then 
\begin{equation} \label{12.29.00.3} 
\int_{\alpha_1 \cdot ... \cdot \alpha_p} \eta_1 \circ  ... \circ \eta_m = \sum 
\int_{\alpha_1}\eta_{1} \circ  ... \circ \eta_{i_1} \cdot ... \cdot 
\int_{\alpha_p}\eta_{i_{p-1}+1} \circ  ... \circ \eta_{m}
\end{equation}
where the summation is over all $i_1 + ... +  i_{p}=m$ where $i_k$ are nonnegative integers. 
In other words we decompose the set 
of forms $\{\eta_1, ..., \eta_m\}$ into a union of $p$ subsets, some of which might be empty, 
by picking the first $i_1$ forms, then the next $i_2$, after this the next $i_3$, and so on. Then we take the iterated integral of the forms in the $p$-th subset 
along the path $\alpha_p$ and   multiply all these integrals. The factor corresponding to 
the empty subset is $1$. Finally, we  take the 
sum over all such decompositions. 

In particular let $\sigma_i$ be loops based at a fixed point on a manifold $X$. 
Denote by $<\int \eta_{1} \circ   ... \circ 
\eta_{k}, \gamma>$ the evaluation of the iterated integral on an element $\gamma 
\in {\cal P}(X; a,b)$. 
 Then 
\begin{equation} \label{12.29.00.2}
\left<\int \eta_{1} \circ \eta_{2} \circ  ... \circ 
\eta_{k}, (\sigma_1 -1)\cdot ... \cdot  (\sigma_r -1)\right> = \left\{ \begin{array}{llll}
 \prod_{i=1}^k \int_{\sigma_i} \eta_i&  k=r \\ 
0&  k<r \end{array} \right.
 \end{equation}
Indeed, take a particular term of the product formula. It corresponds to a decomposition of the set of forms $\{\eta_1, ..., \eta_k\}$ into $r$ subsets as above. 
If one of the subsets is empty the corresponding factor is zero:
$$
<\int \emptyset, \sigma_i-1> \quad = \quad <\int \emptyset, \sigma_i> - <\int \emptyset, 1> \quad = \quad 1-1 \quad =\quad 0
$$

Let ${\cal C} = \{{i_1}, ..., {i_k}\}$ where $1 \leq i_p \leq n$. Denote by 
$X({\cal C})$ the monomial $X_{i_1} ... X_{i_k}$. 
Let us  define elements 
\begin{equation} \label{12.28.00.1}
\gamma(z_0; {\cal C}; z_{n+1}) \in \quad 
{\cal P}^{k}(X_Z; z_0, z_{n+1})
\end{equation}
 Say that a loop in 
$X_Z$ is a simple loop around point $z_i$ if it is homotopic to zero 
in $\C  - \{Z- z_i\}$ and its index around $z_i$ is $\pm 1$. 
Take a path $\gamma$ between the points  
$z_0$ and $z_{n+1}$ in $X_Z$, choose points 
$x_{i_1}, ..., x_{i_k}$ on $\gamma$ ordered  by  the orientation of the path, denote by $\gamma(x_{i_p}, x_{i_q})$ the segment of path $\gamma$ between 
$x_{i_p}$ and $x_{i_q}$, and choose 
simple loops $\sigma_{{i_p}}$ based 
at the points $x_{i_p}$ going counterclockwise 
around the points $z_{i_p}$. 
Then 

\begin{equation} \label{1.1.00.1} 
\gamma(z_0; {\cal C}; z_{n+1}):= \quad \gamma(z_0, x_{i_1}) \cdot (\sigma_{i_1} -1)
\cdot \gamma(x_{i_1}, x_{i_2})  \cdot (\sigma_{i_k} -1)
\cdot  ...  \cdot \gamma(x_{i_k}, z_{n+1}) 
\end{equation}
Here we compose the paths from the left to the right, and $1$ in $\sigma_{i_k} -1$ 
 stays for the constant path at the point $x_{i_k}$. 
The following picture illustrates the structure of the element (\ref{1.1.00.1}). 
Let me stress that it is  a linear combination of paths.

\begin{center}
\hspace{4.0cm}
\epsffile{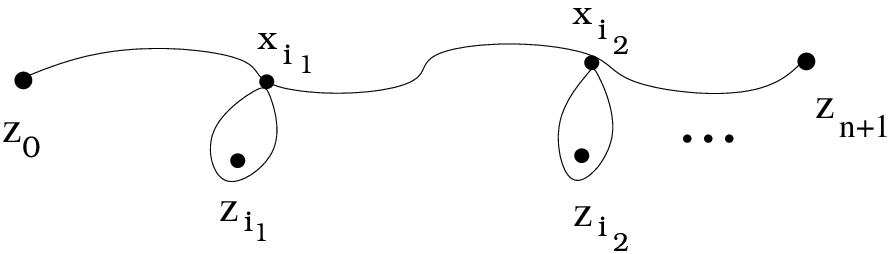}
\end{center}
Observe that the  path on this picture is homotopy equivalent to the following one:
\begin{center}
\hspace{4.0cm}
\epsffile{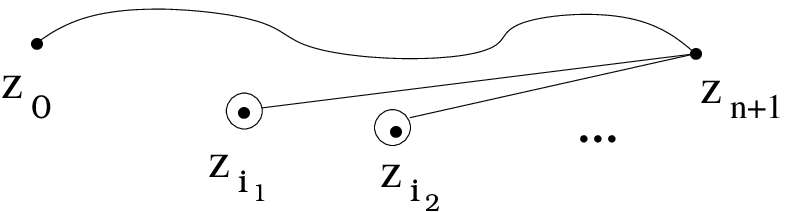}
\end{center}
So the element (\ref{1.1.00.1}) can be written as 
$$
\gamma(z_0; {\cal C}; z_{n+1}) = \quad \gamma \cdot (\widetilde \sigma_{i_1} -1)
\cdot ... \cdot (\widetilde \sigma_{i_k} -1) 
$$
where $\widetilde \sigma_{i_k}$ is a simple loop around $z_{i_k}$ based at $z_{n+1}$.

\begin{lemma} \label{4.13.01.11}
\begin{equation} \label{3.12.01.4}
\int_{\gamma(z_0; {\cal C}; z_{n+1})}\omega_{z_{j_1}} \circ ... \circ 
\omega_{z_{j_r}} =0
\end{equation}
unless ${\cal C}$ is an ordered subset of $\{z_{j_1}, ..., z_{j_r}\}$.  
\end{lemma}

{\bf Proof}. 
Shrinking the loop $\sigma_i$ into a point we see that 
\begin{equation} \label{12.29.00.10}
\mbox{if} \quad i \not \in {\cal J}:= \{j_p, ..., j_q\}  \quad \mbox{ then} \quad 
\int_{\sigma_i} \omega_{a_{j_p}} \circ  ... \circ 
\omega_{a_{j_q}}   = 0
\end{equation}
 Thus  
if $i \in {\cal C}, i \not \in {\cal J}$  then
$<\int \omega_{a_{j_p}} \circ  ... \circ 
\omega_{a_{j_q}}, \sigma_i -1> =0$,  and then we apply (\ref{12.29.00.3}). 
The lemma follows. 
 
When ${\cal C}$ run through all the 
monomials in $X_{i}$ of degree $\leq m$ the elements 
(\ref{12.28.00.1})  provide a basis 
in ${\cal P}^m(X_Z; z_0, z_{n+1})$.

Let $T^{(m)}(V_Z)$ be the quotient of $T(V_Z)$ 
by the $(m+1)$-th power of the augmentation ideal. It 
is a  graded algebra 
such that the degree of $X_i$ is $-1$, and the natural 
projections $T^{(m+1)} \to T^{(m)}$ 
respect this structure. 
The vector space $\widehat T(V_Z)$ is a projective limit of  the
quotients $T^{(m)}(V_Z)$.  

\begin{proposition} \label{2.22.01.1} a) The projection of the  map $P$ onto $T^{(m)}(V_Z) \otimes \C$ 
factorizes through a map 
\begin{equation} \label{2.22.01.3}
P^m: {\cal P}^m(X_Z; z_0, z_{n+1}) \lra T^{(m)}(V_Z)\otimes \C
\end{equation} 

b) This map is an isomorphism after tensoring by $\C$. 

c) The maps $P^m$ are compatible with the natural projections. 
\end{proposition} 

{\bf Proof}. The part a) follows from (\ref{12.29.00.2}).
The part b) is the special case of Chen's theorem ([Ch]) for $X_Z$. 
In fact in our case it 
follows trivially from lemma \ref{4.13.01.11}. Indeed, both vector spaces 
in (\ref{2.22.01.3}) are of the same dimension, and 
the matrix of the period map  written using  this basis
in ${\cal P}^m(X_Z; z_0, z_{n+1})$ and the natural basis formed by  
monomials in $T^{(m)}(V_Z)$ 
is lower triangular by (\ref{12.29.00.3}) and  (\ref{12.29.00.2}). 
Observe that this also gives another proof of the fact that elements (\ref{12.28.00.1})
 form a basis in ${\cal P}(X_Z; z_0, z_{n+1})$.

Finally, c) is straitforward. The proposition is proved.

It follows from this that we get a Hodge-Tate structure 
on ${\cal P}^m(X_Z; z_0, z_{n+1})$ by taking 
$$
V^B = {\cal P}^m(X_Z; z_0, z_{n+1}), \quad V^{DR} = T^{(m)}(V_Z)
$$ 
and defining the period map by the operator $P^m$. The weight filtration on 
${\cal P}^m(X_Z; z_0, z_{n+1})$ is given by the powers of the augmentation ideal. 


{\it The  Hodge-Tate structure on the torsor of path between the tangential base points}.
Choose a tangent vector $v_{a} = \lambda_a \partial_t$ at every point $a \in \C$.  
In s. 2.9 we defined for an 
arbitrary configuration of points 
$a_0, ..., a_{m+1}$ 
in $\C$ and a path $\gamma$ between the tangent vectors $v_{a_0}$ and 
 $v_{a_{m+1}}$ the regularized integral
\begin{equation} \label{12.21.00.1}
{\rm I}_{\{\gamma, v_a\}}(a_0; a_1, ..., a_m; a_{m+1})
\end{equation} 
as the constant term of the  asymptotic expansion as $\varepsilon \to 0$ of 
the integral 
\begin{equation} \label{12.21.00.112}
\int_{a_0 + \varepsilon \lambda_{a_0}}^{a_{m+1} + \varepsilon \lambda_{a_{m+1}}} 
\frac{dt}{t-a_1} \circ ... \circ \frac{dt}{t-a_m}
\end{equation}

Let $X_{0}:= X_{z_0}; X_{n+1}:= X_{z_{n+1}}$  and $\lambda_{0}:= \lambda_{z_0}; 
\lambda_{n+1}:= \lambda_{z_{n+1}}$.  

\begin{lemma} \label{4.22.01.1}
The regularized ordered exponent is given by the formula 
$$
P\int_{v_{z_0}}^{v_{z_{n+1}}}\Omega:= \quad \lim_{\varepsilon_i \to 0} 
(\varepsilon_1 \lambda_1)^{X_{0}}\cdot
\Bigl(P\int_{z_0+ \varepsilon_1 \lambda_{0}}^{z_{n+1} 
+ \varepsilon_2 \lambda_{n+1}} 
\Omega\Bigr)\cdot 
(\varepsilon_2 \lambda_2)^{-X_{{n+1}}}
$$
In particular the limit on  the right exists.  
\end{lemma}

{\bf Proof}.  $f(x,y):= P\int_{x}^{y} \Omega$ as a function on $y$ satisfies 
the differential equation (\ref{4.22.01.1}). So it 
has an asymptotics $(y-z_{n+1})^{X_{n+1}}$ 
near the pont $z_{n+1} \in Z$. Thus $\Bigl(P\int_{x}^{z_{n+1} 
+ \varepsilon} 
\Omega\Bigr)\cdot 
\varepsilon^{-X_{{n+1}}}$ has no singularities as 
$\varepsilon \to 0$. Reversing the path of integration and using (\ref{4.25.01.1}), 
we see that 
 $f(x,y)$ 
has an asymptotics $(x-z_{0})^{-X_{0}}$ near $z_{0} \in Z$. The lemma is proved.

It is straitforward to  generalize all the results of s. 5.3 above by replacing
 everywhere $z_0$ and $z_{n+1}$ by tangential base points 
at $z_0, z_{n+1}$ and the iterated integrals by its regularized values. 
In particular there are basis elements
$$
\gamma_{\{v_a\}}(v_{z_0}; {\cal C}; v_{z_{n+1}}) \in 
{\cal P}^k(X_Z; v_{z_0}, v_{z_{n+1}})
$$
and lemma \ref{4.13.01.11} holds.

{\bf 4. The canonical framed Hodge-Tate  structures related to 
multiple polylogarithms}. Recall the coordinate $t$ 
in $\C P^1 \backslash \{\infty\}$ which identifies it with $\C$. 
Below we set  $v_a:= \partial_t$. 
Let $(a_0; a_1, ..., a_m; a_{m+1})$ be an arbitrary set of 
points in $\C$. Some or all of them may coincide. Let
$$
Z:= \quad \mbox{the union of 
all the points $a_1, ..., a_m$}
$$
\begin{definition} \label{4.12.01.20}
The framed Hodge-Tate structure 
\begin{equation} \label{4.13.01.111}
{\rm I}^{\cal H}_{\{v_a\}}(a_0; a_1, .., a_m; a_{m+1})
 \end{equation}
is the Hodge-Tate structure 
${\cal P}^m(X_Z; v_{a_0}, v_{a_{m+1}})$ with the framing morphisms $$
\Q(0) \lra {\rm Gr}^W_0{\cal P}^m; \qquad 
{\rm Gr}^W_{-2m}{\cal P}^m \lra \Q(m)$$ provided by the  elements
\begin{equation} \label{3.5.01.4}
1 \in {\rm Gr}^W_0{\cal P}^m(X_Z; v_{a_0}, v_{a_{m+1}})
\quad \mbox{and} \quad 
\end{equation}
\begin{equation} \label{3.5.01.5} 
d\log(t-a_1) \otimes ... \otimes d\log(t-a_m) \in \Bigl({\rm Gr}^W_{-2m}{\cal P}^m(X_Z; v_{a_0}, v_{a_{m+1}})\Bigr)^{\vee}
\end{equation}
 
\end{definition}

The period of this framed Hodge-Tate structure is given by the regularized 
integral (\ref{12.21.00.1}).

If $Z'$ is another subset of $\C$ containing $Z$ then 
the Hodge-Tate structure ${\cal P}^m(X_{Z'}; v_{a_0}, v_{a_{m+1}})$ framed 
by (\ref{3.5.01.4}) and (\ref{3.5.01.5}) is equivalent to 
${\rm I}^{\cal H}_{\{v_a\}}(a_0; a_1, .., a_m; a_{m+1})$. Indeed, the inclusion 
$X_{Z'} \subset X_Z$ provides a morphism of the Hodge-Tate structures 
respecting the frames. 

Below we construct explicitly an equivalent framed Hodge-Tate structure 
$\widetilde {\rm I}^{\cal H}_{\{v_a\}}(a_0; a_1, .., a_m; a_{m+1})$
by taking a quotient of  ${\cal P}^m$. This quotient is described in theorem 
\ref{2.22.01.111} below. These Hodge-Tate structures form a 
unipotent variation over an appropriate configuration space, 
see theorem \ref{3.5.01.1}. In particular 
$\widetilde {\rm I}^{\cal H}_{\{v_a\}}(0; x^{-1}, 0, ..., 0; 1)$ is isomorphic to the 
Hodge-Tate structure (\ref{4.20.01.1}).

{\it Construction}. 
Let ${\cal S}:= \{a_1, ..., a_m\}$ and 
 ${\cal O}\{{\cal S}\}$ is the set of  isomorphism classes of ordered subsets of 
${\cal S}$. We define  $V_{\cal S}^{DR}$ a graded $\Q$-vector space 
with a natural basis $h_{{\cal A}}$ 
labelled by the elements ${\cal A} \in {\cal O}\{{\cal S}\}$ 
such that the degree of $h_{{\cal A}}$ is $-|{\cal A}|$. We will also use 
the dual basis  $\{f_{\cal B}\}$ to the basis $\{h_{\cal B}\}$. The Hodge filtration 
$F^{\bullet}$ on $V_{\cal S}^{DR}$ is provided by this grading, 
see (\ref{1.1.01.5}), i.e.  

\vskip 3mm \noindent
$\qquad \qquad F^{-k}V^{DR}_{{\cal S}}$ is spanned by the vectors 
$h_{\cal B}$ for $|{\cal B}| \leq  k$. 
\vskip 3mm \noindent

The vector space $V_{\cal S}^{B}$ is 
a $\Q$-vector space with a weight filtration $W_{\bullet}$ such that 
\vskip 3mm \noindent
$\qquad \qquad W_{-2k} {V}^{B}_{{\cal S}} = W_{-2k+1} {V}^{B}_{{\cal S}} = $ $\Q$-span  of the vectors 
$e_{{\cal A}}$ with $|{\cal A}| \geq k$
\vskip 3mm \noindent
for a certain basis $e_{{\cal A}}$ labelled by the elements of ${\cal O}\{{\cal S}\}$.

Let us construct the period map 
$
P_{\cal S}: {V}^{B}_{{\cal S}} \lra {V}^{DR}_{{\cal S}}\otimes \C
$.  
Let 
$$
{\cal A} = \{a_{i_1}, ..., a_{i_k}\},\qquad 
{\cal B} = \{a_{j_1},  ... ,  a_{j_r}\}
$$

 Denote by $\omega_a$ the 1-form $d\log (t-a)$. We define the period map by setting
\begin{equation} \label{1.1.00.2}
 <P_{\cal S} e_{{\cal A}}, f_{{\cal B}}> := \quad \int_{\gamma_{\{v_a\}}
(a_0; {\cal A}; a_{m+1})} 
\omega_{a_{j_1}} \circ \omega_{a_{j_2}} \circ  ... \circ 
\omega_{a_{j_r}} 
\end{equation}

\begin{proposition} \label{1.1.00.3} 
a) The complexification of the period map $P_{\cal S}$ is an isomorphism. 

b) The image $P_{{\cal S}}(W_{\bullet}{V}^{B}_{{\cal S}})$ of the weight 
filtration  is well defined,
i.e. it does not depend on the choice of path $\gamma$ and loops 
$\sigma_i$.

c) The triple $({V}^{B}_{{\cal S}}, {V}^{DR}_{{\cal S}}, P_{{\cal S}})$ defines a 
Hodge-Tate structure, 
denoted 
\begin{equation} \label{3.5.01.3}
\widetilde {\rm I}_{\{v_a\}}(a_0; {\cal S}; a_{m+1}) \quad = \quad 
\widetilde {\rm I}_{\{v_a\}}(a_0; a_1, ..., a_m; a_{m+1})
\end{equation}
 and equipped 
with a natural framing provided by $(h_{\emptyset}, (2\pi i)^{|{\cal S}|}h_{\cal S})$.
\end{proposition}

{\bf Proof}. a) Consider a partial order $<$ on the set ${\cal O}\{{\cal S}\}$ such that 
${\cal A} < {\cal B}$ if and only if there are ordered subsets ${\cal A}' \subset 
{\cal B}'\subset  {\cal S}$ such that ${\cal A}$ is the isomorphic to ${\cal A}'$ 
and ${\cal B}$  to ${\cal B}'$. We claim that the period matrix (\ref{1.1.00.2}) 
is a lower triangular matrix with 
\begin{equation} \label{1.2.01.2}
<P_{\cal S} e_{{\cal A}}, f_{\cal A}> = (2\pi i)^{|{\cal A}|};\qquad      <P_{\cal S} e_{{\cal A}}, f_{\cal B}> = 0 \quad \mbox{if} \quad {\cal A} \not < {\cal B}
\end{equation}
The  first equality follows from (\ref{12.29.00.2}), and the second from lemma 
\ref{4.13.01.11}. 
This implies a).

b) When we deform  the path $\gamma$ by moving it across one of the points $a_i$ 
the difference between the paths before and after the deformation is an element of  shape 
(\ref{1.1.00.1}) corresponding to an ordered subset $\widetilde {\cal A}$ containing 
${\cal A}$ with $|\widetilde {\cal A}| = |{\cal A}| +1$. Indeed, 
moving the segment $\gamma(x,y)$ across $a_i$ we picking up an extra 
simple loop around the point $a_i$, which is  based on $\gamma(x,y)$. 
If $\widetilde {\cal A} \not < {\cal B}$ then  $<P_{\cal S}e_{\widetilde {\cal A}}, 
f_{{\cal B}}> = 0$. 
If $\widetilde {\cal A}  < {\cal B}$  we changed the vector 
$P_{\cal S}e_{{\cal A}}$ by adding the lower weight vector $P_{\cal S}e_{\widetilde {\cal A}}$ to it. Then c) is clear. 
The proposition is proved. 

{\bf 5. The framed Hodge-Tate structures (\ref{3.5.01.3}) and (\ref{4.13.01.111}) 
are equivalent}. 
{\it The set up}. 
Suppose that $Z = \{z_1, ..., z_n\}$ is a collection of distinct points in $\C$, 
and $a_i \in Z$ for $i=1, ..., m$. 
Let $a_0, a_{m+1}$ be two points which may coincide with 
some of the points $z_i$. 

Recall the ordered subset ${\cal S} = \{a_1, ..., a_m\}$. 
Let ${\cal C} = \{i_1, ..., i_k\}$ where $1 \leq i_p \leq n$, and 
$X({\cal C}):=X_{i_1} ... X_{i_k}$. 
Consider the natural projection
$$
\varphi^{DR}_{{\cal S}}: T^{(m)}(V_Z) \lra V^{DR}_{{\cal S}}; \qquad 
\varphi^{DR}_{{\cal S}}(X({\cal C})) := 
\left\{ \begin{array}{ll}
 h_{\cal C} & \mbox{if} \quad {\cal C} \in {\cal O}\{{\cal S}\}
\\
 0 & \mbox{otherwise} \end{array}\right.
$$
It is a surjective map by its very definition, and it  obviously respects the grading.

Recall the  elements 
$$
\gamma_{\{v_a\}}(v_{a_0}; {\cal C}; v_{a_{m+1}}) \in  {\cal P}^m(X_Z; v_{a_0}, v_{a_{m+1}})
$$
  We define a map 
$$
\varphi^{B}_{{\cal S}}: {\cal P}^m(X_Z; v_{a_0}, v_{a_{m+1}}) \lra V^{B}_{{\cal S}}
$$
between the Betti spaces by the following formula:
$$
\varphi^{B}_{{\cal S}}: \gamma_{\{v_a\}}(v_{a_0}; {\cal C}; v_{a_{m+1}}) \lms \quad 
\left\{ \begin{array}{ll}
 e_{\cal C} & \mbox{if} \quad {\cal C} \in {\cal O}\{{\cal S}\}
\\
 0 & \mbox{otherwise} \end{array}\right.
$$
The map $\varphi^{B}_{{\cal S}}$ clearly preserves the weight filtration. 

The basis elements $\gamma_{\{v_a\}}(v_{a_0}; {\cal C}; v_{a_{m+1}})$ and the period 
matrix $P_{\cal S}$ depend on the choice of the path $\gamma$ and loops $\sigma_i$. 
However the weight filtration is independent of these choices.

\begin{theorem} \label{2.22.01.111} The projections  $\varphi^{DR}_{{\cal S}}$  and 
$\varphi^{B}_{{\cal S}}$ 
provide a surjective morphism of 
Hodge-Tate structures  
\begin{equation} \label{3.5.01.6}
\varphi_{{\cal S}}: {\cal P}^m(X_Z; v_{a_0}, v_{a_{m+1}}) \lra \widetilde 
{\rm I}_{\{v_a\}}(a_0; {\cal S}; a_{m+1})
\end{equation}
 which is compatible with the framing 
(\ref{3.5.01.4})-(\ref{3.5.01.5}) on ${\cal P}^m(X_Z; v_{a_0}, v_{a_{m+1}})$ and the framing 
$(h_{\emptyset}, (2\pi i)^{|{\cal S}|}h_{{\cal S}})$ on (\ref{3.5.01.3})
\end{theorem} 

{\bf Proof}. We need to show that the following diagram is commutative. 
$$
\begin{array}{ccc}
T^{(m)}(V_Z)\otimes \C & \stackrel{\varphi_{{\cal S}}^{DR}\otimes \C}{\lra}& V_{\cal S}^{DR}\otimes \C\\
&&\\
P \uparrow && \uparrow P_{\cal S}\\
&&\\
{\cal P}^m(X_Z; v_{a_0}, v_{a_{m+1}}) & \stackrel{\varphi_{{\cal S}}^B}{\lra}& V_{\cal S}^{B}
\end{array}
$$
Formula (\ref{3.12.01.4}) immediately implies this. The theorem 
is proved.

{\it The unipotent variation of Hodge-Tate structures}. 
Let us say that a configuration of points $(a_0; a_1, ..., a_m; a_{m+1})$ and 
$(a'_0; a'_1, ..., a'_m; a'_{m+1})$ have the same combinatorial type if $a_i = a_j$ if and only if $a'_i = a'_j$. Denote by $M_{\cal S}$ the moduli space of configurations 
 of the same combinatorial type as the configuration 
$(a_0; {\cal S}; a_{m+1})$.

Consider the trivialized complex holomorphic bundle ${\cal V}^{DR}_{\cal S}$ over 
$M_{\cal S}$ with 
 fiber ${V}^{DR}_{{\cal S}} \otimes \C$. It comes
 equipped with a Hodge filtration ${\cal F}^{\bullet}$ provided by the Hodge 
filtration in
 the fiber.

The $\Q$-local system ${\Bbb V}^B_{{\cal S}}$ over $M_{\cal S}$ 
is a subsheaf of ${\cal V}^{DR}_{{\cal S}}$, and  
its fiber over a configuration of points 
$(a_0; {\cal S}; a_{m+1})$ is spanned over $\Q$ by the vectors  
$P_{\cal S} e_{\cal A}$. It is 
equipped with a weight filtration given by $P(W_{\bullet}V^B_{{\cal S}})$. 
The connection on this local system is the Gauss-Manin connection; it is given   
 by continuous deformation of  the integration paths (\ref{12.28.00.1}) 
in the complex plane. 
Proposition  \ref{1.1.00.3} insures that 
we get  a well defined $\Q$-local system filtered by the weight filtration.

\begin{theorem} \label{3.5.01.1} The canonical Hodge-Tate structures 
$\widetilde {\rm I}_{\{v_a\}}(a_0; a_1,  ..., a_m; a_{m+1})$ form  a unipotent 
variation of Hodge-Tate structures of the space of all configurations 
points $(a_0; a_1, .., a_m; a_{m+1})$ of given combinatorial type. 

In particular the multiple polylogarithm Hodge-Tate structures 
$\widetilde {\rm I}_{n_1, ..., n_m}(a_1:  ... : a_{m+1})$ form a unipotent 
variation over the configuration space of distinct points $M_{0,m}$. 
\end{theorem}

{\bf Proof}. It follows immediately from theorem \ref{2.22.01.111} that 
the map $\varphi^{DR}_{{\cal S}}$ gives rise  to a quotient of the canonical 
 unipotent variation of Hodge-Tate structures $ {\cal P}^m(X_Z; v_{a_0}, v_{a_{m+1}})$.

{\bf Warning}. These  Hodge-Tate structures  
 {\it do not} 
provide a variation over 
the configuration space of all points $(a_0; a_1, .., a_m; a_{m+1})$ since their 
dimensions can jump.

{\it Specialization}. 
Consider the variation $\widetilde I(z_0, a_1, ..., a_m; z_{m+1})$ over the 
space 
\begin{equation} \label{3.5.01.10}
(z_0, z_{m+1}) \in 
\Bigl(\C - \{a_1, ..., a_m\}\Bigr) \times \Bigl(\C - \{a_1, ..., a_m\}\Bigr) 
\end{equation}
The specialization  functor transforms 
 it to a variation of Hodge-Tate structures 
on the tangent space to the point $(a_0, a_{m+1})$ punctured at zero. 
\begin{lemma} \label{5.6}
The fiber at the vector $(v_{a_0}, v_{a_{m+1}})$ 
of the specialization  functor applied to 
the variation of Hodge-Tate structures 
$\widetilde I(z_0, a_1, ..., a_m; z_{m+1})$ over the space (\ref{3.5.01.10})  is 
canonically isomorphic to the Hodge-Tate structure 
$\widetilde I_{v_a}(a_0, a_1, ..., a_m; a_{m+1})$.
\end{lemma}

{\bf Proof}. The definition of the specialization of a 
unipotent variation of 
Hodge-Tate structures convinient  for our purposes is spelled  
by R. Hain in ch. 7 of [H]. We assume the reader is familiar with it. 
Its comparizon  with the canonical regularization
 procedure used in s. 2.9 and 5.3 makes the lemma obvious. 
Let me add some comments which might help the reader. 
As explained in [H],   
the only problem is how to define the Betty subspace of the specialization, 
and it is solved as follows. For simplicity consider the specialization when $z_0 \to a_0$. 
Let $\varepsilon:= z_0-a_0$. We consider the restriction of our variation into a 
little neighborhood of the point $z_0$, and view it as a variation on the punctured 
disc with canonical coordinate 
$\varepsilon$. Denote by 
$s_1, .., s_p$ a basis of the local system at the point $\varepsilon =1$.  Let $N$ be 
the logarithm of the monodromy matrix in this basis, see (\ref{3.10.01.3}). The sections 
$\widetilde s_i(\varepsilon):= s_i(\varepsilon) \varepsilon^{-N}$ 
have no monodromy around the point $\varepsilon =0$, and extend to this point.  Then 
$\{\widetilde s_i(0)\}$ form a basis of the Betty space of the specialization. 
 In our situation the entries of the columns $s_i(\varepsilon)$ are iterated integrals 
which admit an asymptotic expansion near $\varepsilon =0$ of form 
$$
f(\varepsilon) = f_0(\varepsilon)  + f_1(\varepsilon) \log(\varepsilon) + ... + 
f_n(\varepsilon) \log^n(\varepsilon)
$$ 
where $f_i(\varepsilon)$ are continuous at $\varepsilon =0$. 
Such an asymptotic expansion is unique if exists. Thus 
the entries of the column 
$\widetilde s_i(\varepsilon)$ are the regularized values 
of the iterated integrals. 
The lemma is proved.  

{\bf Remark}. In this and the next chapters we use 
the tangential base points $\partial_t$ at every point $a \in \C$. 
On the other hand for the canonical regularization discussed in chapter 2 we used 
the tangential base point $\partial_t$ at $a_0$ and $-\partial_t$ at $a_{m+1}$. 
However the element of the Hopf algebra ${\cal H}_n$ obtained 
using the latter tangential base points coincides modulo $2$-torsion 
with the one $\widetilde I_{\{v_a\}}(a_0;{\cal S} ;a_{m+1})$.

{\bf 6. An explicit computation of the multiple polylogarithm  Hodge-Tate structures}. 
Let us define another  period map $P_{\cal S}': V_{{\cal S}}^B \lra 
V_{{\cal S}}^{DR}$. Set 
$<{P}'_{\cal S}e_{\cal A}, f_{\cal B}> =0$ if ${\cal A} \not < {\cal B}$. 
To define the rest of the matrix coefficients we need some notations. 
Let $ {\cal A} \subset {\cal B}$ be two ordered subsets of ${\cal S}$. 
We present the ordered set ${\cal B}- {\cal A}$ as a union of ordered subsets 
${\cal B}\{i_p, i_{p+1}\}$, where $p=0, ..., k$ such  that  
$$
{\cal B} \quad = \quad {\cal B}\{0  \to  i_{1}\}\cup  \{a_{i_1}\} \cup {\cal B}\{i_1\to  i_{2}\} \cup 
\{a_{i_2}\} \cup ... \cup \{a_{i_r}\} \cup {\cal B}\{i_{r} \to  m+1\} 
$$
as ordered sets.  Then   
\begin{equation} \label{1.5.00.1}
<P'_{\cal S}e_{\cal A}, f_{\cal B}>:= \quad (2\pi i)^{|{\cal A}|} 
\sum_{{\cal A}' \subset {\cal B}} \prod_{p=0}^r 
{\rm I}_{\{\gamma, v_a\}}(a_{i_p}; {\cal B}\{i_p \to i_{p+1}\}; a_{i_{p+1}})
\end{equation}
where the sum is over all ordered subsets  ${\cal A}' \subset {\cal B}$ which are 
isomorphic to ${\cal A}$ as ordered sets. 

\begin{proposition} \label{3.1.01.1} 
The Hodge-Tate structures defined by the 
period maps  $P_{\cal S}'$ and $P$ are canonically isomorphic.   
\end{proposition}

{\bf Proof}. To calculate integral (\ref{1.1.00.2}) we deform the element $\gamma(a_0; {\cal A}; a_{m+1})$ by shrinking the 
loops $\sigma_{i_k}$ around the points $a_{i_k}$, and moving the segments 
$\gamma(x_{i_k}, x_{i_{k+1}})$ towards  paths $\overline \gamma(a_{i_k}, a_{i_{k+1}})$ 
between the tangential base points 
at $a_{i_k}$ and $a_{i_{k+1}}$. Denote 
the limiting element by 
$\overline \gamma(a_{0}; {\cal A}; a_{m+1})$. 
Its shape  is illustrated  on the picture. 

\begin{center}
\hspace{4.0cm}
\epsffile{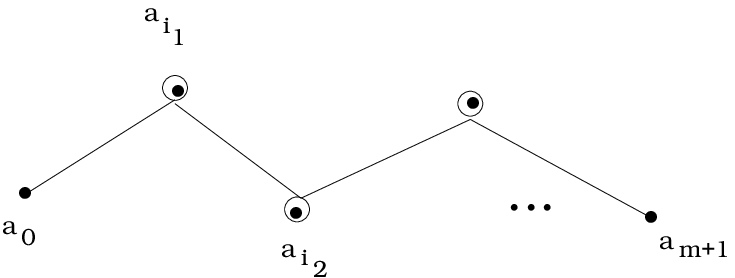}
\end{center}

It follows from the product formula (\ref{12.29.00.3}) and lemma \ref{hov3} that 
\begin{equation} \label{1.6.01.1}
P\int_{\overline \gamma(a_{0}; {\cal A}; a_{m+1})}\Omega \quad = \quad 
P\int_{\overline \gamma(a_{0}, a_{i_{1}})}\Omega \cdot \prod_{p=1}^k 
\Bigl((e^{2\pi i X_{i_k}} - 1)  P\int_{\overline \gamma(a_{i_k}, a_{i_{k+1}})}\Omega\Bigr)
\end{equation}
To check this consider the path $\gamma_{\varepsilon}$ near the limiting path 
$\overline \gamma(a_{0}; {\cal A}; a_{m+1})$ such that the loops $\sigma_i$ became the radius 
$\varepsilon$ circles around $a_i$. Then thanks to lemma \ref{hov3} 
$\lim_{\varepsilon \to 0}\int_{\sigma^{\varepsilon}_i} \omega_{i_p} \circ ... \circ\omega_{i_q} =0$ 
of one of the forms in integrand is different from $\omega_{i}$. Further, 
$$
\int_{\sigma^{\varepsilon}_i} d\log(t-a_i) \circ ... \circ d\log(t-a_i) = 
\frac{(2\pi i)^n }{n!}
$$
This proves formula (\ref{1.6.01.1}). 

We use this  to compute  integral (\ref{1.1.00.2}) over 
$\overline \gamma(a_{0}; {\cal A}; a_{m+1})$. 
Observe that if 
 $a = a_{i_k} = a_{i_{k+1}}$ we choose  $\overline \gamma(a_{i_k}, a_{i_{k+1}})$ to
be the identity  path at $v_a$. Assuming this, if 
 ${\cal B}\{i_k \to i_{k+1}\} \not = \emptyset$ then 
${\rm I}_{\{\gamma, v_a\}}(a_{i_k};{\cal B}\{i_k \to i_{k+1}\} ; a_{i_{k+1}}) =0$, and the corresponding matrix element $<P'_{\cal S}e_{{\cal A}}, f_{{\cal B}}> =0$. 

The period matrices $P_{\cal S}$ and $P'_{\cal S}$ differ by  
multiplication on a unipotent operator in $V_{\cal S}^B$ acting as identity on 
${\rm Gr}^WV_{\cal S}^B$. For the  generic configurations of points $\{a_i\}$ 
these period matrices coincide. 
In general the contribution to the integral coming only from 
the linear terms in $e^{2\pi i X_{i_k}} - 1$ matches the terms in (\ref{1.5.00.1}). 
The contribution of the higher degree terms amounts to adding the lower weight columns 
to the one $P'_{\cal S}e_{\cal A}$. The proposition follows. 

{\bf 7.  An example:  variations of framed Hodge-Tate structures related to  the 
multiple logarithms}. 
{\it The multiple logarithm Hodge-Tate structures}. 
Let  
${\cal S}_m$ be the set of all subsets of $\{1, 2, ..., m\}$, 
so $|{\cal S}_m|=2^m$. 
Suppose that $a_0, a_1, ... , a_{m+1}$ are distinct points on the complex plane. Let 
$$
{\cal I} = \{0 < i_1 < ... < i_k < m+1\}, \quad 
{\cal J} = \{0 < j_1 < ... < j_r < m+1\}
$$ 
Then ${\cal I}, {\cal J} \in {\cal S}_m$ and the period map 
$P_{{\cal S}_m}: V^B_{{\cal S}_m} \lra V_{{\cal S}_m}^{DR} \otimes \C$ is given by by
\begin{equation} \label{12.29.00.4}
 <P_{{\cal S}_m} e_{\cal I}, f_{\cal J}> := \quad \int_{\gamma(a_0; a_{i_1}, ..., a_{i_k}; a_{m+1})} 
\omega_{a_{j_1}} \circ \omega_{a_{j_2}} \circ  ... \circ 
\omega_{a_{j_r}} 
\end{equation}
The integration path is defined as in (\ref{1.1.00.1}).

{\it Explicit computation}. 
Below we choose a certain system of path $\{\gamma\}$ connecting pairs of the points 
$a_{\alpha}, a_{\beta}$, $\alpha < \beta$. Let us 
define a period map 
${P}'_{{\cal S}_m}$. First, 
$<{P}'_{{\cal S}_m}{\{\gamma\}}e_{\cal I}, f_{\cal J}>
 =0$ if ${\cal I} \not \subset {\cal J}$. 
The other matrix coefficients are defined as follows. 
Suppose that ${\cal I} \subset {\cal J}$. Let 
\begin{equation} \label{4.17.01.1}
{\cal I} = \{0 < i_1 < ... < i_k < m+1\}\qquad i_0 := 0, \quad i_{k+1}:= m+1
 \end{equation}
The complement to $I$ in $J$ is the disjoint union 
$$
{\cal J} - {\cal I} \quad = \quad {\cal J}{\{0\ra 1\}} \cup {\cal J}{\{1\ra 2\}} \cup  ... \cup  
{\cal J}{\{k\ra k+1\}}  
$$
where ${\cal J}{\{p \ra p+1\}}$ is the subset of  $ {\cal J}\backslash {\cal I}$ 
between $i_p $ and $i_{p+1}$. 
We  use a shorthand
$$
{\rm I}_{\gamma}(a_{i_p}; a_{J{\{p \ra p+1 \}}}; a_{i_{p+1}}):=  \quad 
{\rm I}_{\gamma}(a_{i_p}; a_{j_{p}^1}, ..., a_{j_{p}^{k_p}}; a_{i_{p+1}})
$$ 
By definition it is  $1$  if 
$J{\{p \ra p+1\}}$ is the empty set. Set
\begin{equation} \label{yCO0}
<{P}'_{{\cal S}_m} e_{{\cal I}}, f_{\cal J}> := \quad (2\pi i)^{k} \prod_{p=0}^k
{\rm I}_{\gamma(a_{i_p}, a_{i_{p+1}})}(a_{i_p}; a_{J{\{p \ra p+1}\}}; a_{i_{p+1}})
\end{equation}
The collection of path ${\{\gamma\}}$ is chosen as follows. 
For a given subset $I$ we choose  paths 
$$
\gamma(a_{i_{0}}, a_{i_{1}}), \quad ... \quad \gamma(a_{i_{k-1}}, a_{i_{k}}), 
\quad \gamma(a_{i_k}, a_{i_{m+1}})
$$
where the path $\gamma(a,b)$ connects the points $a$ and $b$, and its interior 
part is in 
$\C - \{a_0, ..., a_{m+1}\}$. 
The paths for different $I$'s may be unrelated.

 \begin{corollary} \label{CPa2} Multiple logarithm Hodge-Tate structures form 
 a variation of Hodge-Tate structures over the configuration space $M_{0,m}$.  
\end{corollary}

{\bf Proof}. These is a special case of theorem \ref{3.5.01.1}.

{\bf Example}. {\it The variation of Hodge-Tate structures for the double logarithm}. The matrix of the operator $P'_{{\cal S}_2} (a_0; a_1, a_2; a_3)$ is 
$$ 
\left (\matrix{1&0&0&0 \cr {{\rm I}}_{\gamma(0,3)}(a_0; a_1; a_3) & 2 \pi i &0&0 \cr {{\rm I}}_{\gamma(0,3)}(a_0; a_2; a_3) &0&
    2 \pi i&0\cr {\rm I}_{\gamma(0,3)}(a_0; a_1, a_2; a_3) 
& 
2 \pi i \cdot {{\rm I}}_{\gamma(1,3)}(a_1; a_2; a_3) &    
 2 \pi i \cdot {{\rm I}}_{\gamma(0,2)}(a_0; a_1; a_2)&(2 \pi i)^2 \cr}\right ) 
$$ 
The columns are the vectors $P'_{{\cal S}_2}e_{\{\emptyset\}}, 
P'_{{\cal S}_2}e_{\{1\}}, P'_{{\cal S}_2}e_{\{2\}}, P'_{{\cal S}_2}e_{\{1,2\}}$ 
written in the basis $h_{J}$. 
The conditions that the columns form a horizontal basis of a 
local system are the following:
$$
\nabla h_{\{1,2\}} = 0; \quad \nabla h_{\{2\}}  = - d{\rm I}(a_0; a_1; a_2) 
h_{\{1,2\}}; \quad 
\nabla h_{\{1\}} = - d{\rm I}(a_1; a_2; a_3) h_{\{1,2\}}; \quad 
$$
$$
\nabla \Bigl(h_{\{\emptyset\}} + {\rm I}_{\gamma(0,3)}(a_0; a_1; a_3) h_{\{1\}} + 
{\rm I}_{\gamma(0,3)}(a_0; a_2; a_3) h_{\{2\}} + 
{\rm I}_{\gamma(0,3)}(a_0; a_1, a_2;a_3) 
h_{\{1,2\}}\Bigr) = 0
$$
Notice that $d{\rm I}(a; b; c)$ is independent 
of the choice of path $\gamma$, and we used the same path 
$\gamma(0,3)$ in the other three 
terms. 

Solving the last equation we see that the Griffiths transversality condition is equivalent to
$$
\nabla h_{\{\emptyset\}} = - d{\rm I}(a_0; a_1; a_3) h_{\{1\}} - d{\rm I}(a_0; a_2; a_3) h_{\{2\}}
$$
and  the differential equation 
$$
d{\rm I}_{\gamma(0,3)}(a_0; a_1, a_2;a_3) \quad = \quad 
{\rm I}_{\gamma(0,3)}(a_0; a_1; a_3) d{\rm I}(a_1; a_2; a_3) + 
{\rm I}_{\gamma(0,3)}(a_0; a_2; a_3) d{\rm I}(a_0; a_1; a_2)
$$
 which is valid thanks to (\ref{300}).

Let us  show that the Griffiths transversality condition for $P'_{{\cal S}_m}$  
is equivalent to the differential 
equations for multiple logarithm function.  
 Define a connection 
$\nabla$ as  follows. Let
$$
0=j_0 < j_1 < ... < j_r < j_{r+1} = m+1 \quad \mbox{and} \quad  
{\cal J} = \{j_1, ..., j_r \}
$$ 
Then  set 
$$
\nabla h_{\cal J} = - \sum_{\alpha =0}^{r} \sum_{j_{\alpha} < s < j_{\alpha+1}}
d{\rm I}(j_{\alpha}; s; j_{\alpha+1})\cdot h_{\{{\cal J} \cup s\}}
$$
The Griffiths transversality condition is obviously valid 
for this
 connection. However we have to prove that this connection coincides with 
the one defined above, i.e  
$\nabla(P'_{{\cal S}_m} e_{\cal J}) = 0$. This would also imply 
 that it is 
an integrable connection. 
Let us check first that the section $P'_{{\cal S}_m}e_{\emptyset}$ 
is covariantly constant. 
We use a shorthand ${\rm I}(a_{\cal J}):= {\rm I}(a_0; a_{j_1}, ..., a_{j_r}; a_{m+1})$. 
One has
$$
\nabla (\sum {\rm I}(a_{\cal J}) h_{\cal J}) \quad = \quad 
\sum_{{\cal J}' = {\cal J} \cup s} \Bigl(d{\rm I}(a_{{\cal J}'}) - {\rm I}(a_{\cal J}) 
d{\rm I}(j_{\alpha}; s; j_{\alpha+1})\Bigr)
h_{{\cal J}'}
$$
The right hand side of this formula is precisely the differential equation 
(\ref{300}), and so it is zero. The computation for the general 
column $P'_{{\cal S}_m}e_{\cal I}$ is similar: the restriction of the connection to the 
subspace with the coordinates $a_i$, $i \in 
{\cal J}{\{p \to p+1\}}$ is computed just  as 
 above.

\section{The multiple polylogarithm Hopf algebra}

Recall the (commutative) 
Hopf 
algebra ${\cal H}_{\bullet}$ of framed Hodge-Tate structures over $\Q$. 
In this section we calculate explicitly the coproduct for the multiple polylogarithm framed  
Hodge-Tate structures. This immediately implies that 
they form a Hopf algebra, which we call the multiple polylogarithm Hopf algebra and denote 
${\cal Z}^{\cal H}_{\bullet}(\C^*)$. 
We conjecture that it is isomorphic to  (the dual to) the 
fundamental Hopf algebra of the (still hypothetical) 
abelian category of mixed Tate motives over $\C$. 

The framed Hodge-Tate structures
related to multiple polylogarithms whose arguments 
belong to a subgroup $G \subset \C^*$  provide a 
graded Hopf subalgebra ${\cal Z}^{\cal H}_{\bullet}(G)$ of ${\cal H}_{\bullet}$. 

In particular when $G = \mu_N$ the  spectrum 
of the Hopf algebra ${\cal Z}^{\cal H}_{\bullet}(\mu_N)$ 
is a prounipotent algebraic group over $\Q$. 
The grading of the Hopf algebra is provided  
by an action of ${\Bbb G}_m$ on this group. 
The semidirect product of ${\Bbb G}_m$ and this 
group is isomorphic to 
the image of the  fundamental group of the category of Hodge-Tate structures 
acting on the Hodge realization ${\cal P}^{{\cal H}}({\Bbb G}_m - \mu_N; v_0, v_{1})$
 of the  
torsor of path.  
The corresponding Lie algebra ${C}_{\bullet}(\mu_N)$ is called the 
 cyclotomic Lie algebra. 

There is a  natural depth filtration on the  Hopf algebra  
${\cal Z}^{\cal H}_{\bullet}(\C^*)$ which induces the depth filtration on the 
Hopf algebras ${\cal Z}^{\cal H}_{\bullet}(G)$, and on the corresponding Lie algebras. 

The cobracket 
for the  multiple polylogarithm Lie coalgebra 
was defined in [G0]. The formulas for the coproduct in the Hopf algebra 
in the depth 2 case see in [G2], theorem 4.5. A good exercise is to show that 
they are equivalent to the ones given below.  The formulas for the coproduct 
 in the multiple logarithm 
case see in theorem 3.6 in [G5].

{\bf 1. The coproduct in the case of multiple logarithms.} 
Recall the framed 
Hodge-Tate structure
${\widetilde {\rm I}}(a_0; a_1, a_2, ..., a_m; a_{m+1})$  related to the multiple logarithm 
at a point of $M_{0,m}$. 

\begin{theorem} \label{CP1} Let us assume that $a_i \not = a_j$ for $i \not = j$. 
Then the coproduct of the framed Hodge-Tate structure corresponding to the multiple 
logarithm is computed as follows: 
\begin{equation} \label{CP2}
\Delta {\widetilde {\rm I}}(a_0; a_1, a_2, ..., a_m; a_{m+1})= 
\end{equation}
$$
\sum_{0 = i_0 < i_1 < ... < i_k < i_{k+1} = m} {\widetilde {\rm I}}(a_0; a_{i_1}, ..., 
a_{i_k}; a_{m+1}) \otimes \prod_{p =0 }^k
{\widetilde {\rm I}}(a_{i_{p}}; a_{i_{p}+1}, ..., a_{i_{p+1}-1}; a_{i_{p+1}})
$$
\end{theorem}

{\bf Proof}. It follows immediately from the definitions and 
description of the Hodge-Tate 
structure corresponding to ${\rm I}(a_0; a_1, a_2, ..., a_m; a_{m+1})$ given in s. 5.7. 
Let me spell the details. 

 One has 
$$
Gr^W_{-2k}{\widetilde {\rm I}}(a_0; a_1, a_2, ..., a_m; a_{m+1}) \quad = \quad 
\Q(k) \boxtimes V^{DR}_{-2k}
$$
A basis  in $V^{DR}_{-2k}$ is given by  $\{h_{\cal I}\}$, $|{\cal I}|=k$, where 
${\cal I}$ is as in  (\ref{4.17.01.1}).

We claim that there is  a projection of Hodge-Tate  structures 
$$
\varphi_{\cal I}: \quad {\widetilde {\rm I}}(a_0; a_1, a_2, ..., a_m; a_{m+1}) \quad \lra \quad 
{\widetilde {\rm I}}(a_0; a_{i_1},  ..., a_{i_k}; a_{m+1})
$$
given as follows. For any subset ${\cal A} \subset \{1, ..., m\}$ one has:
$$
\varphi_{\cal I}^B(e_{\cal A}) = \left\{ \begin{array}{llll}
 0&  {\cal A} \not \subset {\cal I} \\ 
e_{\cal A}&  {\cal A} \subset {\cal I} \end{array} \right.; \qquad \varphi_{\cal I}^{DR}(h_{\cal A}) = \left\{ \begin{array}{llll}
 0&  {\cal A} \not \subset {\cal I} \\ 
h_{\cal A}&  {\cal A} \subset {\cal I} \end{array} \right. 
$$
To show that we get a morphism of Hodge-Tate structures one needs to 
check commutativity of the diagram
$$
\begin{array}{ccc}
V^B_{{\cal S}_m}& \stackrel{\varphi^B_{\cal I}}{\lra}& V^B_{\cal I}`\\
&&\\
P_{{\cal S}_m}\downarrow && \downarrow P_{\cal I}\\
&&\\
V^{DR}_{{\cal S}_m}\otimes \C&\stackrel{\varphi^{DR}_{\cal I}}{\lra}&
V^{DR}_{\cal I}\otimes \C
\end{array}
$$
Indeed, 
$$
P_{{\cal S}_m}e_{\cal A}= \sum_{{\cal A} \subset {\cal B}} 
(P_{{\cal S}_m}e_{\cal A}, f_{\cal B})h_{\cal B};
\qquad P_{\cal I}e_{\cal A}= \sum_{{\cal A} \subset {\cal B} \subset {\cal I}} 
(P_{\cal I}e_{\cal A}, f_{\cal B})h_{\cal B}
$$
and 
$$
(\varphi^{DR}_{\cal I} \circ P_{{\cal S}_m}) e_{\cal A} = 
\sum_{{\cal A} \subset {\cal B} \subset {\cal I}} (P_{{\cal S}_m} e_{\cal A}, f_{\cal B})h_{\cal B};
\qquad (P_{\cal I} \circ \varphi^{B}_I) e_{\cal A} = 
\sum_{{\cal A} \subset {\cal B} \subset {\cal I}} (P_{I} e_{\cal A}, f_{\cal B})h_{\cal B}
$$
Observe that if ${\cal A} \subset {\cal B} \subset {\cal I}$ then $
(P_{{\cal S}_m} e_{\cal A}, f_{\cal B}) = (P_{\cal I}e_{\cal A}, f_{\cal B})$ 
by the very definitions.  
The commutativity of the diagram follows.

The morphism $\varphi_I$ obviously preserves the frames, thus 
provides an equivalence 
of the framed Hodge-Tate 
structures 
$$
({\widetilde {\rm I}}(a_0; a_1, a_2, ..., a_m; a_{m+1}); h_{\emptyset}, (2\pi i)^{|{\cal I}|}h_{\cal I})  
\quad \sim \quad ({\widetilde {\rm I}}(a_0; a_{i_1},  ..., a_{i_k}; a_{m+1}); h_{\emptyset}, (2\pi i)^{|{\cal I}|}h_{\cal I}) 
$$ 
Notice  that the framing  on the right is the standard one.  

Now consider the subspaces $$
V_{({\cal I})}^B \subset V_{{\cal S}_m}^B, \qquad 
V_{({\cal I})}^{DR} \subset V_{{\cal S}_m}^{DR}
$$ 
generated by all the basis vectors $e_{\cal B}$ and $h_{\cal B}$ such that ${\cal I} \subset {\cal B}$. 
Thanks to (\ref{1.2.01.2}) 
the period map induces the one $P_{({\cal I})}: V_{({\cal I})}^B \lra 
V_{({\cal I})}^{DR}\otimes \C$, and we get a Hodge-Tate structure 
$H_{({\cal I})}$. 
So the natural inclusion $H_{({\cal I})} \subset 
{\widetilde {\rm I}}(a_0; {\cal S}_m; a_{m+1})$ is a map of Hodge-Tate structures. 
It induces an equivalence
$$
[H_{({\cal I})}; (2\pi i)^{|{\cal I}|}h_{\cal I}), (2\pi i)^{m}h_{{\cal S}_m}] 
\quad\sim \quad [{\widetilde {\rm I}}(a_0; {\cal S}_m; a_{m+1}); 
(2\pi i)^{|{\cal I}|}h_{\cal I}), (2\pi i)^{m}h_{{\cal S}_m}]
$$
The description of the Hodge-Tate 
structure  $\widetilde {\rm I}(a_0; a_1, ..., a_m; a_{m+1})$ given in s. 5.7 shows 
that there is an isomorphic of mixed Hodge structures
\begin{equation} \label{1.2.01.3}
H_{({\cal I})}\quad = \quad \otimes_{p=0}^k
\widetilde {\rm I}(a_{i_p}; a_{J{\{p \ra p+1}\}}; a_{i_{p+1}})
\end{equation}
where the frame of $H_{({\cal I})}$ 
is provided  by  $h_{\cal I}$ and $h_{\{1, ..., m\}}$, 
and the Hodge-Tate structure on the right is framed as a tensor product of 
framed Hodge-Tate structures. 
The isomorphism obviously respects the frames. Thus there is an equivalence 
of the framed Hodge-Tate structures
$$
\Bigl({\widetilde {\rm I}}(a_0; a_1, a_2, ..., a_m; a_{m+1});  (2\pi i)^{|{\cal I}|}h_{\cal I}, 
(2\pi i)^{m}h_{\{1, ..., m\}}\Bigr)(-|{\cal I}|) = 
$$
$$ 
\otimes_{p=0}^k
\widetilde {\rm I}(a_{i_p}; a_{J{\{p \ra p+1}\}}; a_{i_{p+1}})
$$ 
 The theorem is proved.

{\bf Examples}. 1. {\it The double logarithm} (compare with s. 2.2 in [G2]). 
The coproduct of the framed Hodge-Tate structure 
${\widetilde {\rm I}}(a_0; a_1, a_2; a_3)$ 
is given by the formula
$$
\Delta {\widetilde {\rm I}}(a_0; a_1, a_2; a_3) \quad = \quad 1 \otimes {\widetilde {\rm I}}(a_0; a_1, a_2; a_3) +
$$
$$
{\widetilde {\rm I}}(a_0; a_1; a_3) \otimes {\widetilde {\rm I}}(a_1; a_2; a_3) +  
{\widetilde {\rm I}}(a_0; a_2; a_3) \otimes {\widetilde {\rm I}}(a_0; a_1; a_2) +
{\widetilde {\rm I}}(a_0; a_1, a_2; a_3) \otimes 1 
$$

2. {\it The triple logarithm}. The coproduct of the framed Hodge-Tate structure ${\widetilde {\rm I}}(a_0; a_1, a_2, a_3;a_4)$ 
is given by 
$$
\Delta {\widetilde {\rm I}}(a_0; a_1, a_2, a_3; a_4) \quad = \quad 1 
\otimes {\widetilde {\rm I}}(a_0; a_1, a_2, a_3; a_4) +
$$
$$
{\widetilde {\rm I}}(a_0; a_1; a_4) \otimes {\widetilde {\rm I}}(a_1; a_2, a_3; a_4) +  
{\widetilde {\rm I}}(a_0; a_2; a_4) \otimes {\widetilde {\rm I}}(a_0; a_1; a_2) \cdot {\widetilde {\rm I}}(a_2; a_3; a_4) +
$$
$$
{\widetilde {\rm I}}(a_0; a_3; a_4) \otimes {\widetilde {\rm I}}(a_0; a_1, a_2; a_3) +  
{\widetilde {\rm I}}(a_0; a_1, a_2; a_4) \otimes {\widetilde {\rm I}}(a_2; a_3; a_4) +
$$
$$
{\widetilde {\rm I}}(a_0; a_1, a_3; a_4) \otimes {\widetilde {\rm I}}(a_1; a_2; a_3) +
{\widetilde {\rm I}}(a_0; a_2, a_3; a_4) \otimes {\widetilde {\rm I}}(a_0; a_1; a_2) +
{\widetilde {\rm I}}(a_0; a_1, a_2, a_3; a_4) \otimes 1 
$$

{\bf 2. The coproduct in the general case}. Below we introduce another framed Hodge-Tate 
structure 
$\widehat {\rm I}_{\{\gamma, v_a\}}(a_0; a_1, ..., a_m; a_{m+1})$ corresponding to 
an {\it arbitrary} configuration 
of points $(a_0; a_1, ..., a_m; a_{m+1})$. It coincides with the canonical one 
if the points $a_i$ are distinct, but in 
general it is bigger then the  canonical mixed Hodge  structure 
$\widetilde {\rm I}_{\{v_a\}}$
However we show that they are  equivalent as 
{\it framed} Hodge-Tate structures. 
The second one 
is  handy for the computation of the coproduct.

Choose for each pair of indices $\alpha < \beta$ a path between the tangential base points 
$v_{a_{\alpha}}$ and $v_{a_{\beta}}$. Denote by  $\{\gamma\}$ this collection of path. 
Then define the period map  
\begin{equation} \label{12.30.00.1}
\widehat P_{\{\gamma, v_a\}}(a_0; a_1, ..., a_m; a_{m+1}): 
\quad V_{{\cal S}_m}^B \lra V_{{\cal S}_m}^{DR}\otimes \C
\end{equation} 
by copying 
 the definitions from s. 5.7, replacing everywhere 
the integrals ${\rm I}_{\gamma}$ there by their regularized values
${\rm I}_{\{\gamma, v_a\}}$.

\begin{proposition} \label{12.30.00.2}
For any 
configuration of points $(a_0; a_1, .., a_m; a_{m+1})$ the period map  (\ref{12.30.00.1}) 
provides a well defined 
framed Hodge-Tate  structure 
$
\widehat {\rm I}_{\{v_a\}}(a_0; a_1, .., a_m; a_{m+1})
$
\end{proposition}

{\bf Proof}. One needs to check  that 
a different choice of the collection of 
path $\{\gamma\}$ between the tangential base points does not change 
the weight filtration on the Betti space $V_{{\cal S}_m}^B$. This is done the same 
way as the proof of proposition \ref{1.1.00.3}. The proposition is proved. 

\begin{proposition} \label{1.1.01.1} 
The framed Hodge-Tate structures 
$$
\widehat {\rm I}_{\{v_a\}}(a_0; a_1, ..., a_m; a_{m+1})\quad \mbox{and} \quad  
\widetilde {\rm I}_{\{v_a\}}(a_0; a_1, ..., a_m; a_{m+1})
$$ 
are equivalent. 
\end{proposition} 

{\bf Proof}. Thanks to proposition \ref{3.1.01.1} we can work with the Hodge-Tate structure 
${\widetilde {\rm I}}'_{\{v_a\}}(a_0; {\cal S}; a_{m+1})$ 
defined by the map $P'_{{\cal S}}$. The map of sets 
$$
\varphi: \mbox{subsets of  $\{1, ..., m\}$} \quad \lra 
\quad {\cal O}\{S\}; \qquad \{i_1, ..., i_k\} \lms 
\{a_{i_1}, ..., a_{i_k}\}
$$
provides linear maps
$
\varphi^B: V^B_{{\cal S}} \lra V^B_{2^m}, \quad \varphi^{DR}: V^{DR}_{{\cal S}} \lra V^{DR}_{2^m}
$ 
given by 
$$
\varphi^B(e_{{\cal A}}):= \sum_{\varphi(A) ={\cal A} }e_{{ A}}; \qquad 
\varphi^{DR}(h_{{\cal B}}):= \sum_{\varphi(B) ={\cal B} }h_{{ B}}
$$
These maps induce a morphism of mixed Hodge structures 
${\widetilde {\rm I}}'_{\{v_a\}}(a_0; {\cal S}; a_{m+1}) \lra 
{\widehat {\rm I}}(a_0; {\cal S}; a_{m+1})$. Indeed, 
$(\varphi^{DR})^* f_{B} = f_{{\varphi(B)}}$, so comparing (\ref{1.5.00.1}) and (\ref{yCO0}) we see that 
the maps $\varphi^B,\varphi^{DR}$ commute with the period maps. 
This morphism obviously preserves the frames and thus induces an 
equivalence of the corresponding framed 
Hodge-Tate structures. The proposition is proved.

\begin{theorem} \label{CP11} Let $(a_0; a_1, a_2, ..., a_m; a_{m+1})$ be an 
arbitrary configuration of points in $\C$. Then the coproduct  
of the canonical framed Hodge-Tate structure 
${\widetilde {\rm I}}_{\{v_a\}}(a_0; a_1, ..., a_m; a_{m+1})$ is given by
\begin{equation} \label{CP2*}
\Delta {\widetilde {\rm I}}_{\{v_a\}}(a_0; a_1,  ..., a_m; a_{m+1})= 
\end{equation}
$$
\sum_{0 = i_0 < i_1 < ... < i_k < i_{k+1} = m} 
{\widetilde {\rm I}}_{\{v_a\}}(a_0; a_{i_1}, ..., a_{i_k}; a_{m+1}) \otimes \prod_{p =0 }^k
{\widetilde {\rm I}}_{\{v_a\}}(a_{i_{p}}; a_{i_{p}+1}, ..., a_{i_{p+1}-1}; a_{i_{p+1}})
$$
\end{theorem}

In other words 
we can drop the assumption 
$a_i \not = a_j$ in theorem \ref{CP1}. 
  
{\bf Proof}.  Thanks to proposition \ref{1.1.01.1} we can do the calculation 
 with the Hodge-Tate structure $\widehat {\rm I}_{\{v_a\}}(a_0; a_1, ..., a_m; a_{m+1})$.
The shape of the matrix defining this Hodge-Tate structure 
does not change when we restrict to degenerate configurations of points, 
so the proof is identical with the one of theorem \ref{CP1}.

{\bf 3. Calculation of the coproduct for the multiple polylogarithms}. Set
\begin{equation}  \label{12.18.00.2}
\widetilde {\rm I}_{n_0, n_1+1, ..., n_m+1}(a_0; a_1, .., a_m; a_{m+1}):= 
\end{equation} 
$$
\quad \widetilde 
{\rm I}_{\{v_a\}}(a_0; \underbrace{0, ..., 0,}_{\mbox{$n_0$ times}} a_1, 
\underbrace{0, ..., 0,}_{\mbox{$n_1$ times}} a_2,\quad ...  \quad 
\underbrace{0, ..., 0}_{\mbox{$n_m$ times}}; a_{m+1})
 $$
This framed Hodge-Tate structure  is clearly invariant under the translations. Since it 
depends on the 
choice of the tangent vectors $v_a$, in general  it is 
no longer 
invariant under the 
action of ${\Bbb G}_m$ given by 
 $a_i \lms \lambda a_i$. However, if the iterated integral 
${\rm I}_{n_0, ..., n_m}(a_0; a_1, ... , a_m ; a_{m+1})$ is convergent, then 
the corresponding framed Hodge-Tate structure is invariant under the 
action of the affine group, see lemma \ref{9.30.00.2} below.

Let $V$ be a vector space. 
Denote by $V[[t_1, ..., t_m]]$ the vector space of the formal power 
series in $t_i$ whose coefficients are vectors of $V$. 
Let us make the generating series
\begin{equation} \label{CO1}
\widetilde {\rm I}(a_0; a_1, ... , a_m ; a_{m+1}| t_0; t_1; ... ; t_m) := 
\end{equation}
$$
\sum_{n_i \geq 0}\widetilde {\rm I}_{n_0, n_1+1, ..., n_m+1}(a_0; a_1, ... , a_m ; a_{m+1}) t_0^{n_0} ... t_m^{n_m} 
 \in \quad {\cal H}_{\bullet}[[ t_0, ..., t_m ]]
$$
To visualize them  
consider a line segment with the following additional data, called 
{\it coloring}:

i) The beginning of the segment is labelled by $a_0$, the end by $a_{m+1}$.

ii) There are $m$ points {\it inside} of the segment labelled by $a_1, ..., a_m$ from 
the left to the right. 

iii) These points cut the segment on $m+1$ arcs labelled  by $t_0, t_1, ..., t_m$.

\begin{center}
\hspace{4.0cm}
\epsffile{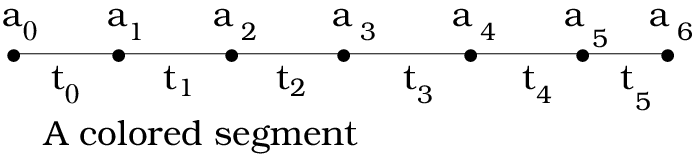}
\end{center}

{\bf Remark}. The way $t$'s sit between the $a$'s reflects  the shape of the iterated integral 
${\rm I}_{n_0, ..., n_m}(a_0; a_1, ..., a_m; a_{m+1})$. 

As we  show in theorem \ref{ur8-4,3} below, the 
terms of the coproduct of the element (\ref{CO1}) correspond to 
the colored segments equipped with the following 
additional data, called marking:

a) Mark (by making them boldface on the picture) 
points $a_0; a_{i_1}, ..., a_{i_{k}}; 
a_{m+1}$  
so that  
\begin{equation}  \label{ur11}
0 = i_0 < i_1 < ... < i_k < i_{k+1} = m+1
\end{equation} 

b) Mark (by cross) segments $t_{j_0}, ..., t_{j_k}$ 
such that there is just one marked segment between any two 
neighboring marked points.

\begin{center}
\hspace{4.0cm}
\epsffile{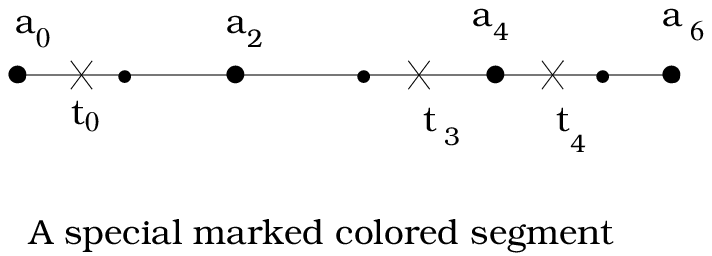}
\end{center}

The conditions on the crosses for a marked colored segment just mean that 
\begin{equation}  \label{ur10}
i_{\alpha } \leq j_{\alpha } < i_{\alpha +1} \quad \mbox{ for any 
$0\leq \alpha \leq k$}
\end{equation}
The marks provide a new colored segment:
\begin{equation}  \label{ur106}
(a_0|t_{j_0} |a_{i_0} | t_{j_1} |
a_{i_1} |\quad ... \quad | a_{i_k} |t_{j_k}|a_{m+1})
\end{equation}

\begin{theorem} \label{ur8-4,3} Let us suppose that $a_i \not = 0$ for $i=1, ..., n$. Then 
\begin{equation}  \label{ur100}
\Delta \widetilde {I}(a_0; a_1, ..., a_m; a_{m+1}| t_0; ...; t_{m}) \quad = \quad \sum \widetilde {\rm I}(a_0; a_{i_1}, ..., a_{i_{k}}; a_{m+1}| t_{j_0}; t_{j_1}; ... ;
  t_{j_{k}}) \otimes
\end{equation}
$$
\prod_{\alpha =0}^k \Bigl(
\widetilde {\rm I}(a_{i_{\alpha}}; a_{i_{\alpha}+1}, ..., a_{j_{\alpha}};0| 
t_{i_{\alpha}}; ... ; t_{j_{\alpha}}) \cdot 
  \widetilde {\rm I}(0; a_{j_{\alpha}+1}, ..., a_{i_{\alpha+1}-1}; a_{i_{\alpha+1}}|
t_{j_{\alpha}}; t_{j_{\alpha}+1}; ... ; t_{i_{\alpha+1}-1} )\Bigr)
$$
where the sum is over all marked colored segments, i.e. over all 
sequences $\{i_{\alpha}\}$ and $\{j_{\alpha}\}$ 
satisfying inequality (\ref{ur10}). 
\end{theorem}

{\bf Proof}. One immediately sees that 
\begin{equation} \label{12.18.00.1}
\widetilde {\rm I}(0; a_1, ..., a_m; 0) =0
\end{equation}
Indeed, it follows from the very definition that the left column of the matrix 
$A_{\gamma}(0; a_1, ..., a_m; 0)$ defining 
the framed Hodge structure $\widetilde {\rm I}(0; a_1, ..., a_m; 0)$ is 
$(1, 0, ..., 0)$, which means that $\Q(0)$ 
is its direct summand. Therefore  we get  (\ref{12.18.00.1}). 

We apply the coproduct formula from theorem \ref{CP1} to the framed mixed Hodge-Tate structure 
(\ref{12.18.00.2}) 
and then keep track of the non zero terms taking into account (\ref{12.18.00.1}).

The left hand side factors of the non zero terms in the formula for the coproduct 
correspond to certain  subsets
$$
A \in \{a_0; \underbrace{0, ..., 0,}_{\mbox{$n_0$ times}} a_1, 
\underbrace{0, ..., 0,}_{\mbox{$n_1$ times}} \quad ...  \quad , a_m,
 \underbrace{0, ..., 0}_{\mbox{$n_m$ times}}; a_{m+1}   \}
$$ 
containing $a_0$ and $a_{m+1}$, and called the admissible subsets. Such a 
subset $A$ determines the subset 
$
I = \{i_1 < ... < i_k\}\quad 
$
where $a_0, a_{i_1}, ..., a_{i_k}, a_{m+1}$ are precisely the set of all  $a_i$'s 
containing in $A$. 
A subset $A$ is called  admissible if 
 it satisfies the following properties:

i) $A$  contains $a_0$ and $a_{m+1}$.
 
ii) The set of $0$'s in $A$ located between 
$a_{i_{\alpha}}$ and $a_{i_{\alpha+1}}$ must be a string of {\it consecutive} $0$'s located 
between $a_{j_{\alpha}}$ and $a_{j_{\alpha}+1}$ for some $i_{\alpha} \leq j_{\alpha} 
<  i_{{\alpha}+1}$.

In other words the factors in the coproduct are parametrized by 

{\it a marked colored segment} and 
{\it a connected string of $0$'s in each of the crossed arcs}. 

The connected string of $0$'s in some  of the crossed arcs might be empty.

 An admissible subset $A$ provides  the following framed Hodge-Tate structure, which is the left hand side of the corresponding term in the coproduct: 
$$
\widetilde {\rm I}(a_0; \widehat a_{j_{0}},  
\underbrace{0, ..., 0,}_{\mbox{$s_{j_0}$ times}}   \widehat a_{j_{0}+1}, 
a_{i_1}, \widehat a_{j_{1}},  \underbrace{0, ..., 0,}_{\mbox{$s_{j_1}$ times}}  
\widehat a_{j_{1}+1}, \quad ... \quad ; a_{m+1})
$$
  So the string of zero's between $a_{i_{\alpha}}$ and $a_{i_{\alpha}+1}$ satisfying ii) 
looks as follows: 
\begin{equation}  \label{ur8-4,5}
\{a_{i_{\alpha}}, \widehat a_{j_{\alpha}}, \underbrace{\widehat 0, ..., \widehat 0, }_{\mbox{$p_{j_{\alpha}}$}} \underbrace{0, ..., 0, }_{\mbox{$s_{j_{\alpha}}$}} \underbrace{\widehat 0, ..., \widehat 0, }_{\mbox{$q_{j_{\alpha}}$}} \widehat a_{j_{\alpha}+1}, a_{i_{\alpha}} \}
\end{equation}
where  $p_{j_{\alpha}} + q_{j_{\alpha}} + s_{j_{\alpha}} = n_{j_{\alpha}}$.
  This notation emphasizes that all $0$'s located between $a_{i_{\alpha}}$ and 
$a_{i_{\alpha+1}}$ are in fact  located between $a_{j_{\alpha}}$ and 
$a_{j_{\alpha}+1}$, and form a connected segment of length $s_{j_{\alpha}}$. 

\begin{center}
\hspace{4.0cm}
\epsffile{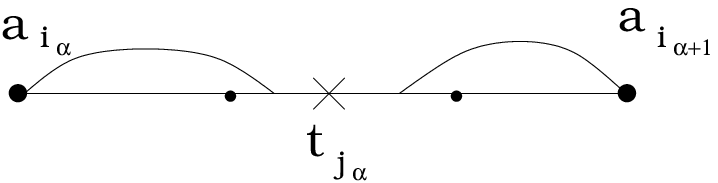}
\end{center}

The right hand side of the term in the coproduct corresponding to the subset $A$ 
is a product over $0 \leq \alpha \leq k$ of framed Hodge-Tate structures of the following shape:
$$
\widetilde {\rm I}(a_{i_{\alpha}}; \underbrace{0, ..., 0, }_{\mbox{$n_{i_{\alpha}}$}} 
a_{i_{\alpha+1}}, \underbrace{0, ..., 0, }_{\mbox{$n_{i_{\alpha}+1}$}} \quad ... \quad a_{j_{\alpha}}, \underbrace{0, ..., 0}_{\mbox{$p_{j_{\alpha}}$}};0) \cdot \underbrace{\widetilde {\rm I}(0;0) \cdot ... \cdot \widetilde {\rm I}(0;0)}_{s_{j_{\alpha}}}
$$
$$
\widetilde {\rm I}(0; \underbrace{0, ..., 0, }_{\mbox{$q_{j_{\alpha}}$}} 
a_{j_{\alpha+1}}, \underbrace{0, ..., 0, }_{\mbox{$n_{j_{\alpha}+1}$}}  \quad ... \quad a_{i_{\alpha+1}-1}, \underbrace{0, ..., 0}_{\mbox{$n_{i_{\alpha+1}-1}$}};
a_{i_{\alpha+1}})
$$
where  the middle factor 
$\widetilde {\rm I}(0;0) \cdot ... \cdot \widetilde {\rm I}(0;0)$ is equal to $1$ since $\widetilde {\rm I}(0;0)= 1$ according to (\ref{12.18.00.1}). 
Translating this into the language of the generating series we get the promised 
formula for the coproduct.  
 The theorem is proved.

{\it A geometric interpretation of formula 
(\ref{ur100}) }. It is surprisingly similar to the 
one  for multiple logarithms. 
Recall that the expression (\ref{CO1}) 
is encoded by a colored segment 
\begin{equation}  \label{ur101} 
(a_0| t_0| a_1| t_1| \quad ... \quad |a_m|t_m|a_{m+1}) 
\end{equation}
The terms of the coproduct are in the bijection with the marked colored 
segments $S$ obtained from  a given colored segment (\ref{ur101}). 
Denote by $L_S \otimes R_S$ the term in the coproduct corresponding to $S$. 
The left factor $L_S$ 
is encoded by the colored segment (\ref{ur106}) obtained from the marked points and arcs. For example 
for the marked colored segment on the picture  we get
$
L_S = \widetilde {\rm I}(a_{0}|t_0|a_2|t_3|a_4|t_4|a_6)
$. 

The marks (which consist of $k+1$ crosses and $k+2$ boldface points) 
determine a decomposition of the segment (\ref{ur101}) 
into $2(k+1)$ 
little colored segments in the following way. Cutting 
the initial segment in all the marked points and crosses we get 
$2(k+1)$ little segments. For instance, the very right one is 
the segment between the last cross and  point $a_m$, and so on.  
Each of these segments either starts
 from a marked point and ends by a cross, or starts from a 
cross and ends by a marked point. 

There is a natural way to make 
 a {\it colored} segment out of each of these  
little segments: mark
 the ``cross endpoint'' of the little segment by the point $0$, and for the arc which is just near to this marked point use the 
letter originally attached to the arc containing it. 
For example the marked colored segment on the picture above produces 
the following sequence of little colored segments:
$$
(a_0|t_0|0), \quad (0|t_0|a_1|t_1|a_2) \quad (a_2|t_2|a_3|t_3|0) \quad (0|t_3| 
a_4) \quad (a_4| t_4|0) \quad (0|t_4|a_5|t_5|a_6)
$$
Then the factor $R_S$ is the product of the generating series for framed Hodge-Tate 
structures corresponding to these little colored segments. 
For example for the marked colored segment on the picture we get 
$$
R_S \quad = \quad \widetilde {\rm I}(a_0;  0|t_0) \cdot \widetilde {\rm I}(0; a_1; a_2|t_0; t_1) 
\cdot 
\widetilde {\rm I}(a_2; a_3; 0|t_2;t_3) \cdot \widetilde {\rm I}(0; a_4|t_3)  \cdot 
\widetilde {\rm I}(a_4; 0|t_4) \cdot \widetilde {\rm I}(0; a_5; a_6|t_4; t_5)
$$ 

To check that formula for the coproduct of the multiple logarithms fits 
into this description we use formula (\ref{ur102}) from the following 
\begin{lemma} \label{9.24.00.1} One has the shuffle product formula 
on the level of the framed Hodge-Tate structures
\begin{equation}  \label{ur102q} 
\widetilde {\rm I}(x; a_{1}, a_{2},  ...  , a_{p}; y)
\widetilde {\rm I}(x; a_{p+1}, a_{2},  ...  , a_{p+q}; y) = 
\end{equation}
$$
\quad 
\sum_{\sigma \in \Sigma_{p,q}}\widetilde 
{\rm I}(x; a_{\sigma(1)}, a_{\sigma(2)},  ...  , a_{\sigma(p+q)}; y)
$$
and 
\begin{equation}  \label{ur102} 
\widetilde {\rm I}(a_{p}; a_{p+1}, a_{p+2},  ...  ; a_{p+q}) = \quad 
\widetilde {\rm I}(0; a_{p+1}, a_{p+2},  ...  ; a_{p+q}) + 
 \end{equation}
$$
\sum^{q-2}_{k=1} \widetilde {\rm I}(a_{p}; a_{{p}+1},  ... a_{p+k}; 0) \cdot
\widetilde {\rm I}(0; a_{p+k+1},  ...,  a_{p+q-1}; a_{p+q})  + 
\widetilde {\rm I}(a_{p}; a_{p+1},  ... a_{p+q-1}; 0)
$$
\end{lemma}
together with  the fact that 
each term  on  the right hand side of this formula correspond to a marked 
colored segment shown on the picture:
\begin{center}
\hspace{4.0cm}
\epsffile{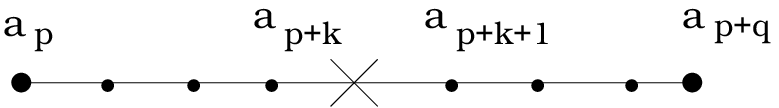}
\end{center}

{\bf Proof of lemma \ref{9.24.00.1}}. Let us prove the shuffle product formula. 
We will assume that $a_i \not = x, a_i \not = y$. The general case follows 
from this by using the specialization. Indeed,  by lemma \ref{5.6} the regularized 
Hodge-Tate structures  $\widetilde {\rm I}(x_0; a_{1}, a_{2},  ...  , a_{p}; y_0)$
can be defined by using the specialization when $x \to x_0, y \to y_0$. 

Let $Z$ be the union of points $a_1, ..., a_{p+q}$. Recall that 
$\widetilde {\rm I}(x; a_{1}, a_{2},  ...  , a_{p}; y)$ can be defined as the 
Hodge-Tate structure 
$$
H_p:= H^p((\C - Z)^p; {\cal D}_p \cap (\C - Z)^p )
$$
where ${\cal D}_p$ is the algebraic simplex given by the union of hyperplanes 
$x =t_1, t_1 = t_2, ..., t_{k-1} = t_k, t_k = y$.  The framing is given by the form 
$d \log (t_1 - a_1) \otimes ... \otimes d \log (t_p - a_p)$, and by a relative cycle 
$\gamma_{{\cal D}_p}$ which provids a generator of $H_p(\C^p; {\cal D}_p)$. 
We construct such a cycle as follows. Let $\gamma: [0,1] \to \C - Z$ 
be a path from $x$ to $y$. Recall the standard simplex 
$$
\Delta_p = \{0 \leq t_1 \leq ... \leq t_p \leq 1\} 
$$ 
Then $\gamma_{{\cal D}_p}:= \gamma(\Delta_p)$. 
Similar interpretation have the other framed Hodge-Tate structures involved in the formula 
(\ref{ur102q}). 
There is a well known  decomposition of the product of simplices
$$
\Delta_p \times \Delta_q \quad = \quad \cup_{\sigma \in \Sigma_{p,q}}\Delta_{\sigma}
$$ 
into the union of simplices $\Delta_{\sigma}$
parametrized to the shuffles $\sigma$. It provides 
a decomposition of ${\cal D}_p \times {\cal D}_q \subset \C^{p+q}$ into the  algebraic simplices denoted 
${\cal D}_{\sigma}$. Let 
$ {\cal D}$ be the union of the hyperplane faces of all the simplices ${\cal D}_{\sigma}$. 
Consider the following 
Hodge-Tate structures 
$$
H:= H^{p+q}((\C - Z)^{p+q};  {\cal D} \cap (\C - Z)^{p+q}); 
$$
$$
 H_{\sigma}:= H^{p+q}((\C - Z)^{p+q}; 
 {\cal D}_{\sigma} \cap (\C - Z)^{p+q} ); 
$$
They have  natural framings. Namely, one component of these framings
 is  provided by the form 
$d \log (t_1 - a_1) \otimes ... \otimes d \log (t_{p+q} - a_{p+q})$. 
The other is given by  
relative cycle $\gamma(\Delta_p \times \Delta_q)$ for $H$ and by 
$\gamma(\Delta_{\sigma})$ for $H_{\sigma}$.  We need to show that 
$H_p \otimes H_q$ is equivalent to $\oplus_{\sigma \in \Sigma_{p,q}} H_{\sigma}$. 

 There are natural morphism 
of Hodge-Tate structures $H \lra  H_{\sigma}$ as well as $H \lra  H_p$ and $H \lra  H_q$.  
They provide the morphisms
$$
H \lra \oplus_{\sigma \in \Sigma_{p,q}}  H_{\sigma}\qquad 
H \lra H_{p} \otimes H_{q}
$$
Since $\gamma(\Delta_{p} \times \Delta_{q}) = \cup \gamma(\Delta_{\sigma})$ these maps 
respect the frames and thus induce the desired equivalence. 
The shuffle product formula is proved. 

The proof of formula (\ref{ur102}) is completely similar. 
Assume first that $a_i \not \in \{ a_p, a_q, 0\}$ for $i \in \{p+1, ..., p+q-1\}$. 
The general case  follows by the specialization. 
Take $Z =\{a_{p+1}, ..., a_{p+q-1}\}$ and consider the composition  $\gamma := 
\alpha \cdot \beta$ of path $\alpha$ from $a_p$ to $0$ and $\beta$ from $0$ to $a_{q}$. Then
$$
\gamma(\Delta_{q-1}) = \cup_{i=0}^{q-1} \alpha(\Delta_{q-i-1}) \times \beta(\Delta_{i})
$$
Consider the Hodge-Tate structure 
\begin{equation}  \label{3.10.01.4}
H^{q-1}((\C - Z)^{q-1}, {\cal D} \cap (\C - Z)^{q-1})
\end{equation}
equipped with several different framings. One component of all these 
framings is the same for all of them and it is provided by 
$$
\frac{dt}{t-a_{p+1}} \otimes ... \otimes \frac{dt}{t-a_{p+q-1}} \in \quad 
H^{q-1}_{DR}((\C - Z)^{q-1})
$$
The second components of these framings are given by relative homology cycles 
$$
\alpha(\Delta_{q-1-i}) \times \beta(\Delta_{i}) \quad \mbox{and} \quad 
\gamma(\Delta_{q-1})
$$
We call $H_i$ and $H$  the corresponding framed Hodge-Tate structures. 
The identity morphism of (\ref{3.10.01.4}) provides a morphism of framed Hodge-Tate 
structures $H \lra \oplus H_i$ which proves the identity. 
The lemma is proved.

\begin{lemma} \label{9.30.00.2} One has 
$
\widetilde {\rm I}(a; 0|t) = a^{-t}$, $ \widetilde {\rm I}(0; a|t) = a^{t}
$, 
and more generally  
\begin{equation}  \label{3.7.01.19}
\widetilde {\rm I}(0; a_1, ..., a_m; a_{m+1}| t_0; t_1; ... ;t_m) \quad = \quad
a_{m+1}^{t_0} \widetilde {\rm I}(a_1: ... : a_m: a_{m+1}| t_1-t_0, ..., t_m-t_0)
\end{equation}
and 
$$
\widetilde {\rm I}(a_1; a_2, ..., a_m; 0| t_1; ...; t_m) \quad = \quad 
(-1)^{m-1}\widetilde {\rm I}(0; a_m, ..., a_2; a_1| -t_m; -t_{m-1}; 
  ...;  - t_1) \quad = 
$$
\begin{equation}  \label{ur8-9}
(-1)^{m-1}a_1^{-t_m} \widetilde {\rm I}(a_m: ... : a_2: a_1| t_m-t_{m-1}, 
 t_m-t_{m-2}, ..., t_m - t_1)
\end{equation}
\end{lemma}

{\bf Proof}. The formula (\ref{3.7.01.19}) is a Hodge version of proposition \ref{GENP}. 
Its proof uses 
the shuffle product formula and specialization and copies 
 the second proof of proposition \ref{GENP}.  

 The second  equality in (\ref{ur8-9})
 follows from the first one and (\ref{3.7.01.19}). The proof of the first 
equality in (\ref{ur8-9}) is completly similar 
to the proof of lemma 
\ref{9.24.00.1}. Namely, it follows by a specialization argument from the equality 
$$
\widetilde {\rm I}(a_1; a_2, ..., a_m; 0) = (-1)^{m-1} 
\widetilde {\rm I}(0; a_m, ..., a_2; a_1)
$$
where we can assume that $a_2, ..., a_m \not = 0, a_1$. To prove this equality 
consider two different framings of  the Hodge-Tate structure 
$$
H^{m-1}((\C - Z)^{m-1}, {\cal D}_{m-1} \cap (\C - Z)^{m-1})
$$
where $Z := \{a_2, ..., a_m\}$ and ${\cal D}_{m-1}$ is the usual algebraic simplex. 
One  is provided by a form $$
d \log(t_2 -a_2) \wedge ... \wedge  d \log(t_m -a_m)\quad \mbox{and a relative cycle} \quad \{\gamma(t_2), ..., \gamma(t_m)\}
$$
where  $\gamma: [0, 1] \to \C - Z$ is a path from $a_1$ to $0$, and 
$0 \leq t_2 \leq ... \leq t_m\leq 1 $. 
The other is given 
by 
$$
d \log(t_m -a_m) \wedge ... \wedge  d \log(t_2 -a_2)
\quad \mbox{and a relative cycle} \quad
\{\gamma^0(t_m), ..., \gamma^0(t_2)\}
$$ where 
$\gamma^0: [1,0] \lra \C -Z$ is the path  opposite to $\gamma$. The components of these frames differ by signes, and the total sign difference is $(-1)^{m-1}$. 
The lemma is proved. 

Recall the notation
$$
\widetilde {\rm I}(a_1: ... : a_{m+1}| t_0: ... :t_{m}):= \quad 
\widetilde {\rm I}(0; a_1, ... , a_m; a_{m+1}| t_0; ... ;t_{m})
$$

A marked colored segment is {\it special} if the first cross is marking 
the segment $t_0$. 
A colored segment with such a  data is called a {\it special marked colored segment}. 
See an example on the picture above. 
The conditions on the crosses for a marked colored segment just mean that 
\begin{equation}  \label{ur10*}
i_{\alpha } \leq j_{\alpha } < i_{\alpha +1} \quad \mbox{ for any 
$0\leq \alpha \leq k$}, \qquad j_0 = i_0 = 0
\end{equation}

\begin{proposition} \label{12.11.00.11} One has
$$
\Delta \widetilde {\rm I}(a_1: ... : a_{m+1}| t_0: ... :t_{m}) \quad = \quad 
\sum \widetilde {\rm I}( a_{i_1}: ...: a_{i_{k}}: a_{m+1}| t_{j_0}: t_{j_1}: ... :
  t_{j_{k}}) \otimes
$$
\begin{equation}  \label{ur100100}
\prod_{\alpha =0}^k \Bigl(
(-1)^{j_{\alpha}-i_{\alpha}}\widetilde {\rm I}(a_{j_{\alpha}}: a_{j_{\alpha}-1}: ... : 
 a_{i_{\alpha}}| 
-t_{j_{\alpha}}: -t_{j_{\alpha}-1}: ... : - t_{i_{\alpha}}) \cdot 
\end{equation}
$$
  \widetilde {\rm I}(a_{j_{\alpha}+1}: ... : a_{i_{\alpha+1}-1}: a_{i_{\alpha+1}}|
t_{j_{\alpha}}: t_{j_{\alpha}+1}: ... : t_{i_{\alpha+1}-1} )\Bigr)
$$
where the sum is over all special marked colored segment, i.e. over all 
sequences $\{i_{\alpha}\}$ and $\{j_{\alpha}\}$ 
satisfying inequality (\ref{ur10*}). 
\end{proposition}

{\bf Proof}. Follows immediately from theorem \ref{ur8-4,3} and   lemma \ref{9.30.00.2}.

Let us define several other generating series 
for framed Hodge-Tate structures related to  multiple polylogarithms. 
Consider the following two pairs of  sets of variables:
$$
i) \quad (x_0, ... , x_m) \quad \mbox{such that $x_0 ...  x_m=1$}; \qquad ii) 
\quad (a_1: ... : a_{m+1})
$$
$$
iii) \quad (t_0: ... : t_m); \qquad iv) 
\quad (u_1, ..., u_{m+1}) \quad \mbox{such that $u_1 +  ... + u_{m+1} = 0$}
$$
The relationship between them is given by 

$$
x_i := \frac{a_{i+1}}{a_i}, \quad i = 1, ..., m+1; x_0 = \frac{a_1}{a_{m+1}}
\qquad  u_i := t_i -t_{i-1}, \quad u_{m+1} := t_0 - t_m
$$
the indices are modulo $m+1$, and is illustrated on the picture below

\begin{center}
\hspace{4.0cm}
\epsffile{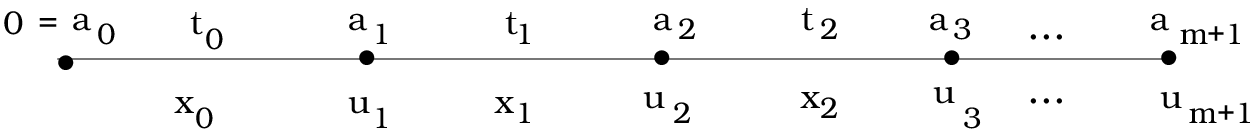}
\end{center}
Observe that $x_i, a_i$ are the multiplicative variables, 
and $t_i, u_i$ are additive variables. 

We introduce the $\widetilde {\rm Li}$-generating series
\begin{equation} \label{11.16.00.1}
{\widetilde {\rm Li}}(*, x_1, ..., x_m| t_0: ... : t_m):= \quad 
\end{equation}
$$
\widetilde {\rm Li}(x_0, ..., x_m| t_0 : ... : t_m) := \quad (-1)^m 
\widetilde {\rm I}( (x_1 ... x_m)^{-1}: 
(x_2 ... x_m)^{-1}: ... : x^{-1}_m :1 | t_0 : ... : t_m) 
$$
$$
=: (-1)^m 
\widetilde {\rm I}( a_1: a_2: ... : a_{m+1}| t_0 : ... : t_m)  \quad =: \widetilde {\rm I}( a_1 : a_2 : ... : a_{m+1}| u_1,  ... u_{m+1})
$$
{\bf Remark}. The $(,)$-notation is used for the variables which sum to zero 
(under the appropriate group structure), and the $(:)$-notation is used those 
sets of variables which are essentially homogeneous with respect to 
 the multiplication by a common factor, 
see lemma \ref{9.30.00.2}. 

{\it The coproduct in terms of the $\widetilde {\rm Li}$-generating series}. 
To state the formula we need the following notations. Let $x_i$ be elements of a group $G$. Set 
$$
X_{a\to b}:= \quad \prod_{s=a}^{b-1} x_{s}
$$

\begin{proposition} \label{ur8-4,3**} Let us suppose that $x_i \not = 0$. Then 
$$
\Delta {\widetilde {\rm Li}}(x_0, x_1, ..., x_m | t_0: t_1: ... :t_{m}) \quad = \quad 
$$
\begin{equation}  \label{ur10000}
\sum 
{\widetilde {\rm Li}}(X_{i_0 \to i_{1}}, X_{i_1 \to i_{2}}, ... , X_{i_k \to m}| 
t_{j_0}: t_{j_1}:  ... :
  t_{j_{k}}) \otimes
\end{equation}
\begin{equation}  \label{12-15.00.1}
\prod_{p = 0}^k \Bigl( (-1)^{j_p-i_p}X_{i_p \to i_{p+1}}^{t_{j_p}}
{\widetilde {\rm Li}} (*, x^{-1}_{{j_p}-1}, x^{-1}_{{j_p}-2}, ..., x^{-1}_{i_p}|-t_{j_{p}}:  
-t_{j_{p}-1}:  ... : -t_{i_{p}}) \cdot 
\end{equation}
\begin{equation}  \label{12-15.00.2}
{\widetilde {\rm Li}} (*, x_{j_{p}+1}, x_{j_{p}+2}, ..., x_{i_{p+1}-1}| 
t_{j_{p}}: t_{j_{p}+1}:  ... : t_{i_{p+1}-1} )\Bigr)  
\end{equation}
Here the sum is over special marked colored segments, i.e.
sequences $\{i_{p}\}$, $\{j_{p}\}$ 
satisfying (\ref{ur10*}). 
\end{proposition}

The $p$-th factor in the product on the right is encoded by the data 
on the $p$-th segment:

\begin{center}
\hspace{4.0cm}
\epsffile{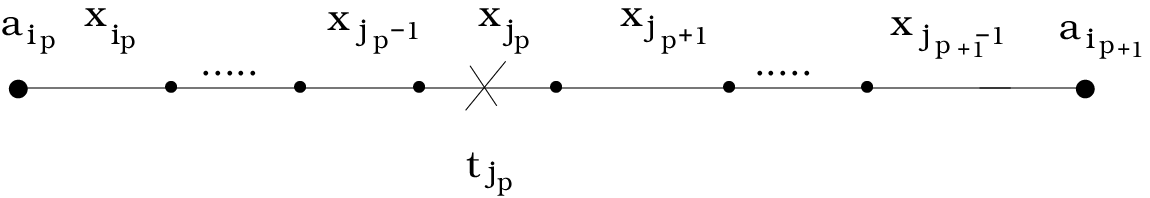}
\end{center}
 Namely, $X_{i_p \to i_{p+1}}$ is the product of all $x_i$ on this segment. 
The first term (\ref{12-15.00.1}) 
is encoded by the segment  between $t_{j_p}$ and $t_{i_p}$, which we read from 
the right to the 
left. The factor (\ref{12-15.00.2})  is encoded by 
the segment  between $t_{j_p}$ and 
$t_{i_{p+1}}$ which we read from the  left to the right.

{\bf Proof}. Follows from lemma \ref{9.30.00.2}, 
proposition \ref{12.11.00.11} and (\ref{11.16.00.1}).

{\bf 4. An example: the coproduct in the classical polylogarithm case}. 
Recall that ${\cal H}_{\bullet}$ is the commutative 
Hopf algebra of the framed $\Q$-Hodge-Tate structures with the coproduct $\Delta$ and the product $\ast$. 

One has $\Delta(1) = 1 \otimes 1, \quad \Delta(x) = 
x \otimes 1 +  1 \otimes x$. 
It is easy to check that 
\begin{equation} \label{d0f}
\Delta: x^t \longmapsto  x^t\otimes x^t
\end{equation}

Recall 
the restricted coproduct: $\Delta'(X): = \Delta(X)
- (X\otimes 1 + 1\otimes X)$. Notice that $\Delta$ is a
homomorphism of algebras and $\Delta'$ is not.

\begin{proposition}  \label{9.30.00.1}
\begin{equation} \label{d1f}
\Delta: {\widetilde {\rm Li}}(x|t) \longmapsto {\widetilde {\rm Li}}(x|t
)\otimes x^{t} + 1 \otimes \widetilde   Li(x|t)
\end{equation}
\end{proposition}

{\bf Proof}. 
This formula just means that 
$$
\Delta' {\widetilde {\rm Li}}_n(x) = {\widetilde {\rm Li}}_{n-1}(x) \otimes \widetilde \log x + {\widetilde {\rm Li}}_{n-2}(x) \otimes \frac{\widetilde \log^2 x}{2} + ... + {\widetilde {\rm Li}}_{1}(x) \otimes \frac{\widetilde \log^{n-1} x}{(n-1)!}
$$

{\bf 5. An example: the coproduct for the depth two
 multiple polylogarithms}. 
We will use both types of $I$ notations for multiple polylogarithms and the corresponding 
Hodge-Tate structures, so for instance 
$$
\widetilde {\rm I}(a_1:a_2:1|t_1,t_2) \quad = \quad 
\widetilde {\rm I}(0; a_1, a_2; 1|t_1;t_2)
$$ 
Set $\widetilde \zeta(t_1, ..., t_m) := \widetilde {\rm I}(1: ... : 1 | t_1, ..., t_m)$.

 \begin{proposition} \label{ur107} a) One has 
$$
\Delta \widetilde {\rm I}(a_1:a_2:1|t_1,t_2) \quad = 1 \otimes 
\widetilde {\rm I}(a_1:a_2:1|t_1,t_2) 
$$
$$
\widetilde {\rm I}(a_1:a_2:1|t_1,t_2) \otimes a_1^{-t_1}\ast a_2^{t_1-t_2} \quad +
 \quad \widetilde {\rm I}(a_1:1|t_1) \otimes a_1^{-t_1}\ast  \widetilde {\rm I}(a_2:1|t_2-t_1)
$$
$$
-\widetilde {\rm I}(a_1:1|t_2) \otimes a_1^{-t_2} \ast \widetilde {\rm I}(a_2:a_1|t_2-t_1)
 \quad + \quad \widetilde {\rm I}(a_2:1|t_2) \otimes \widetilde 
{\rm I}(a_1:a_2|t_1)\ast   a_2^{-t_2}
$$

b) Let us suppose that $a_1^N = a_2^N = 1$. Then modulo the $N$-torsion 
one has 
$$
\Delta' \widetilde {\rm I}(a_1:a_2:1|t_1,t_2) \quad = \quad 
\widetilde {\rm I}(a_1:1|t_1) \otimes   \widetilde {\rm I}(a_2:1|t_2-t_1)
$$
$$
-\widetilde {\rm I}(a_1:1|t_2) \otimes \widetilde {\rm I}(a_2:a_1|t_2-t_1)
 \quad + \quad \widetilde {\rm I}(a_2:1|t_2) \otimes \widetilde 
{\rm I}(a_1:a_2|t_1)
$$
In particular
$$
\Delta' \widetilde \zeta(t_1,t_2) \quad = \quad \widetilde \zeta(t_1) \otimes 
\widetilde \zeta(t_2-t_1) -\widetilde \zeta(t_2) \otimes \widetilde \zeta(t_2-t_1) + \widetilde \zeta(t_2) \otimes \widetilde 
\zeta(t_1)
$$
\end{proposition} 

{\bf Proof}. a) Since in our case $a_0=0$ and  $t_0=0$ the nonzero contribution can be obtained only from those marked colored segments where $t_0$-arc is not marked. Let us call the pictures where $a_0=0$, $a_{m+1} =1$,  $t_0=0$ and 
the $t_0$-arc is not marked by {\it special} marked colored segments. 

The five  terms in the formula above  correspond to the five special 
marked colored segments presented on the picture. 
\begin{center}
\hspace{4.0cm}
\epsffile{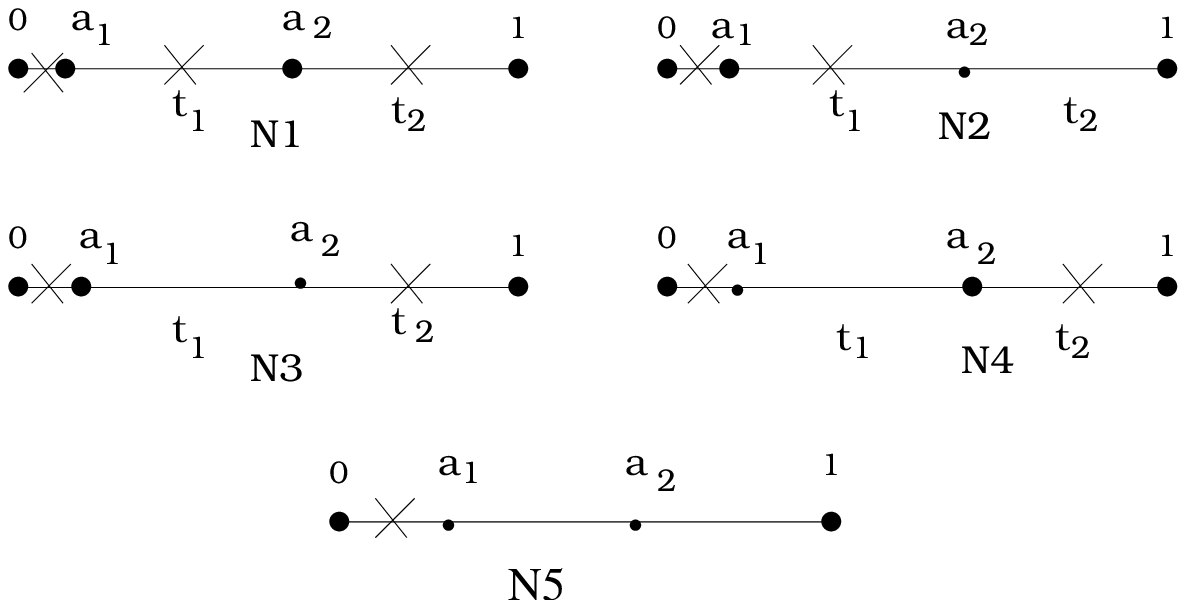}
\end{center}

Using the formulas from lemma \ref{9.30.00.2}
we get the following four terms corresponding to the terms N1-N4 on the 
picture: 
$$
\widetilde {\rm I}(a_1:a_2:1|t_1,t_2) \otimes \widetilde {\rm I}(a_1;0|t_1) \cdot \widetilde {\rm I}(0;a_2|t_1) \cdot 
\widetilde {\rm I}(a_2;0|t_2)\cdot \widetilde {\rm I}(0;1|t_2) = 
$$
$$
\widetilde {\rm I}(a_1:a_2:1|t_1,t_2) \otimes a_1^{-t_1} \cdot a_2^{t_1-t_2}
$$

$$
\widetilde {\rm I}(a_1:1|t_1) \otimes \widetilde {\rm I}(a_1;0|t_1) \cdot  \widetilde {\rm I}(0; a_2; 1|t_1;t_2) \quad = \quad 
\widetilde {\rm I}(a_1:1|t_1)\otimes a_1^{-t_1} \cdot \widetilde {\rm I}(a_2:1|t_2-t_1)
$$

$$
I\widetilde (a_1:1|t_2) \otimes \widetilde {\rm I}(a_1;a_2;0|t_1;t_2) \cdot  \widetilde {\rm I}(0;  1|t_2) = 
$$
$$
- \widetilde {\rm I}(a_1:1|t_2) \otimes \widetilde {\rm I}(0;a_2;a_1|-t_2;-t_1) \quad = \quad - \widetilde {\rm I}(a_1:1|t_2) \otimes a_1^{-t_2}\cdot {\rm I}(a_2: a_1|t_2-t_1) 
$$

$$
\widetilde {\rm I}(0;a_2;1|0;t_2) \otimes \widetilde {\rm I}(0; a_1;a_2|0;t_1) \cdot  {\rm I}(a_2; 0|t_2) 
\cdot  {\rm I}(0;1|t_2) =
$$
$$
\widetilde {\rm I}(a_2: 1|t_2) \otimes \widetilde {\rm I}(a_1 : a_2|t_1) \cdot a_2^{-t_2}
$$

The part b) follows from a) if we notice that $a^N =1$ 
provides $a^t =1$ modulo the $N$-torsion. 
The proposition is proved. 

{\bf Remark}. In theorem 4.5 of [G2] the reader can find a different way to write the 
formulas for the coproduct in the depth $2$ case. It is easy to see that the formulas 
given there are equivalent to the formulas above.


{\bf Example}.  Below we give explicit formulas for  the coproduct in the Hopf algebra for double polylogarithms   up to the  weight three. 
$$
\Delta': {\widetilde {\rm Li}}_{2,1}(x, y) \quad  \longmapsto \quad{\widetilde {\rm Li}}_{1,1}(x, y)\otimes x +
{\widetilde {\rm Li}}_{1}(y) \otimes {\widetilde {\rm Li}}_{2}(x)
 +{\widetilde {\rm Li}}_{2}(x y) \otimes {\widetilde {\rm Li}}_{1}(y)
$$
$$
- {\widetilde {\rm Li}}_{1}(x y) \otimes \Bigl({\widetilde {\rm Li}}_{2}(x) + {\widetilde {\rm Li}}_{2}(y) - {\widetilde {\rm Li}}_{1}(y)\cdot xy + \frac{x^2}{2}\Bigr) 
$$

$$
\Delta': {\widetilde {\rm Li}}_{1,2}(x, y)  \quad   \longmapsto \quad {\widetilde {\rm Li}}_{1,1}(x, y)\otimes y -{\widetilde {\rm Li}}_{2}(xy)\otimes x - {\widetilde {\rm Li}}_{1}(xy)\otimes x \cdot y
$$
$$
+ {\widetilde {\rm Li}}_{2}(y) \otimes {\widetilde {\rm Li}}_{1}(x) + {\widetilde {\rm Li}}_{1}(y) \otimes {\widetilde {\rm Li}}_{1}(x)\cdot y 
+ {\widetilde {\rm Li}}_{1}(x y) \otimes  {\widetilde {\rm Li}}_{2}(y)  
$$
$$
- {\widetilde {\rm Li}}_{2}(x y) \otimes {\widetilde {\rm Li}}_{1}(x) - {\widetilde {\rm Li}}_{1}(x y) \otimes  {\widetilde {\rm Li}}_{1}(x)\cdot xy - {\widetilde {\rm Li}}_{1}(x y) \otimes  {\widetilde {\rm Li}}_{2}(x)   
$$

{\bf Example}. {\it The Lie coalgebra structure in the depth two (see [G0], p. 13)}. 
Recall that 
the space of the indecomposables
$$
{\cal L}({\cal H}_{>0}):= \quad \frac{{\cal H}_{>0}}
{{\cal H}_{>0} \cdot {\cal H}_{>0}}
$$ 
inherits  a structure of the graded Lie coalgebra with the cobracket $\delta$.  
I will describe
$$
\sum _{n_1, ... , n_k > 0} \delta \widetilde {I}_{n_1, ... n_k}(a_1, ... ,a_k) \cdot t_1^{n_1 - 1 } ...
t_k^{n_k - 1 } \quad \in \quad \wedge^2{\cal H}_{\bullet} \widehat \otimes {\Z} [[t_1, t_2, ... ]]
$$
In the formulas below $\delta$ acts  on the first factor in 
${\cal H}_{\bullet} \widehat \otimes {\Z} [[t_1, t_2, ... ]]$.

$$
 \delta \bigr(\sum_{m> 0,n> 0}  {\widetilde {\rm I}}_{m,n}(a:b:c)\cdot t_{1}^{m- 1}t_{2}^{n-1}\bigl) 
 \quad = \quad
$$
$$
\sum_{m> 0, n > 0} \biggr({\widetilde {\rm I}}_{m,n}(a:b:c)\cdot t_{1}^{m- 1}t_{2}^{n-1})
\wedge (\frac{b}{a}\cdot t_{1} + \frac{c}{b}\cdot t_{2}) \quad - \quad  
{\widetilde {\rm I}}_{m}(a:b) \cdot t_{1}^{m- 1} 
\wedge 
{\widetilde {\rm I}}_{n}(b:c) 
\cdot t_{2}^{n-1} 
$$
$$
 + \quad  {\widetilde {\rm I}}_{m}(a:c) \cdot t_{1}^{m- 1} \wedge {\widetilde {\rm I}}_{n}(b:c) 
\cdot (t_{2} - t_{1})^{n-1} 
\quad - \quad {\widetilde {\rm I}}_{n}(a:c) \cdot t_{2}^{n- 1} \wedge 
{\widetilde {\rm I}}_{m}(b:a) 
\cdot (t_{2} - t_{1})^{m - 1} \biggl) 
$$

Set  ${\rm I}_{0,n} = {\rm I}_{n,0} =0$.  Here is a more concrete formula for  $\delta$:  
$$
\delta  {\widetilde {\rm I}}_{m,n}(a,b)
\quad = 
$$
$$
\quad {\widetilde {\rm I}}_{m-1,n}(a,b) \wedge \frac{b}{a} \quad + \quad
{\widetilde {\rm I}}_{m,n-1}(a,b) \wedge \frac{1}{b} \quad 
- \quad {\widetilde {\rm I}}_{m}(\frac{a}{b}) \wedge {\widetilde {\rm I}}_{n}(b) \quad +
$$
$$
\sum_{i=0}^{m-1} (-1)^{i}  { n+ i - 1 \choose i }  {\widetilde {\rm I}}_{m - i}(a) \wedge {\widetilde {\rm I}}_{n + i}(b) 
$$
$$
 - (-1)^{m+j-1} \sum_{j=0}^{n-1} (-1)^{j}  { m+ j - 1 \choose j }  {\widetilde {\rm I}}_{n - j}(a) \wedge {\widetilde {\rm I}}_{m + j}(\frac{b}{a})
$$

{\bf 6. The Hopf algebra of framed Hodge-Tate structures related to multiple polylogarithms}. 
Let $G$ be a subgroup of $\C^*$. 
Denote by ${\cal Z}^{\cal H}_{w}(G) \subset {\cal H}_{w}$ the $\Q$-vector subspace 
generated by the $w$-framed Hodge-Tate 
structures 
\begin{equation} \label{1.11.01.1}
\widetilde {\rm I}_{n_1, ..., n_m} (a_1, ..., a_m), \quad a_i \in G \subset \C^*, 
\quad w = n_1 + ... + n_m
\end{equation}
Set ${\cal Z}^{\cal H}_{\bullet}(G):= \oplus_{ w \geq 1} {\cal Z}^{\cal H}_{w}(G)$.  It is equipped with a depth filtration ${\cal F}_{ \bullet}^{{\cal D}}$. Namely, recall 
the logarithm Hodge-Tate structure $\widetilde \log x \in 
{\rm Ext}_{{\rm MHS}/\Q}^1(\Q(0), \Q(1))$. The depth filtration is defined as follows: 
${\cal F}_{ 0}^{{\cal D}}{\cal Z}^
{{\cal H}}_{\bullet}(G) $ is spanned by  products of $\widetilde \log (a)$, $a \in G$, and 
$
{\cal F}_{ k}^{{\cal D}}{\cal Z}^{{\cal H}}_{\bullet}(G) 
$ for $k \geq 1$ by the elements  (\ref{1.11.01.1}) with 
$m \leq k$.

\begin{theorem} \label{9.30.00.3}
Let $G$ be any subgroup of  $\C^*$. Then 
${\cal Z}^{\cal H}_{\bullet}(G)$ is a graded 
Hopf subalgebra of ${\cal H}_{\bullet}$. 
The depth provides a filtration on this Hopf algebra.
\end{theorem}

{\bf Proof}. 
The graded vector space ${\cal Z}^{\cal H}_{\bullet}(G)$ is closed under the coproduct by  
proposition \ref{12.11.00.11}. The statement about the depth filtration is evident from the formula for coproduct given in proposition \ref{ur8-4,3**}. 
It is a graded algebra  by the shuffle product formula from lemma \ref{9.24.00.1}. 
The theorem  is proved. 

{\bf Remark}. The depth filtration is not defined by a grading of the algebra
${\cal Z}^{\cal H}_{\bullet}(G)$ because of the  
relations like 
$$
\widetilde {\rm Li}_{n}(x) \cdot \widetilde {\rm Li}_{m}(y) \quad = \quad 
\widetilde {\rm Li}_{n,m}(x,y) + 
\widetilde {\rm Li}_{m,n}(y,x) + \widetilde {\rm Li}_{n+m}(xy)
$$

\begin{conjecture} The Hopf algebra ${\cal Z}^{\cal H}_{\bullet}(\C^*)$ is isomorphic to the 
motivic Hopf algebra of $\C$. 
\end{conjecture}
We expect a similar story for any subfield $F \subset \C$.

{\bf 7. The cyclotomic Lie algebras}. Recall the definition
$$
 {\cal C}^{{\cal H}}_{\bullet}(\mu_N):= \quad 
\frac{{\cal Z}^{{\cal H}}_{>0}(\mu_N)}{
{\cal Z}^{{\cal H}}_{>0}(\mu_N) \cdot {\cal Z}^{{\cal M}}_{>0}(\mu_N)}
$$
 \begin{corollary} \label{ur8-11} 
a) ${\cal C}^{{\cal H}}_{\bullet}(\mu_N)$ is a graded Lie coalgebra. 

b). 
${\cal C}^{{\cal H}}_{1}(\mu_N) \quad = \quad {\cal Z}^{{\cal H}}_{1}(\mu_N) \quad \stackrel{\sim}{=} \quad
\Bigl(\mbox{the group of cyclotomic units in $\Z[\zeta_N][1/N]\Bigr) \otimes \Q$}$.
\end{corollary} 

{\bf Proof}. a). Clear. 

b) 
Since $$
\widetilde {\rm I}_1(a) = \widetilde \log (1-a)\quad \mbox{and} \quad  
\widetilde \log (1-a) - \widetilde \log (1-a^{-1}) = \widetilde \log (a)
$$
 the  weight $1$ component ${\cal Z}^{{\cal H}}_{1}(G)$
is generated by $\widetilde \log (1-a)$ and $\widetilde \log (a)$. Notice that if $a^N=1$ then 
$N \cdot \widetilde \log (a) = 0$. This proves b). The corollary is proved.

We call 
${\cal C}^{{\cal H}}_{\bullet}(\mu_N)$ the {\it cyclotomic Lie coalgebra}.   Its dual 
${C}^{{\cal H}}_{\bullet}(\mu_N)$ is the {\it cyclotomic Lie algebra}. 
The dual to the 
universal enveloping algebra  of the cyclotomic Lie algebra is 
isomorphic to ${\cal Z}^{{\cal H}}_{\bullet}(\mu_N)$.

Let ${\cal C}_{m}^{\Delta}(\mu_N)$ be the $\Q$- subspace of  
${\cal C}_{m}(\mu_N)$ generated by 
$
\widetilde {\rm I}(0; a_1, ..., a_m; 1)$ where $ a_i^N =1$. 

\begin{corollary} \label{ur8-11/1} ${\cal C}_{\bullet}^{\Delta}(\mu_N):= \oplus_{m\geq 1}{\cal C}_{m}^{\Delta}(\mu_N)$ 
is a graded Lie coalgebra. 
\end{corollary}

We call it the  
 {\it diagonal cyclotomic Lie coalgebra}. 

{\bf Proof}. Clear from the  proposition \ref{12.11.00.11}. 

{\bf 8. Expressing the multiple polylogarithm Hodge-Tate structures  of weight
  $3$ via the  classical polylogarithms}. Lemma \ref{9.24.00.1} allows to express 
$\widetilde {\rm I}_{1,1}(a_1, a_2)$, and hence $\widetilde {\rm Li}_{1,1}(x,y)$,  
via the dilogarithm and products of logarithms. 
Here are the similar results for the depth three case. 
The key thing is 
the following formula
$$
\widetilde {\rm Li}_{1,2}(x,y) \quad =\quad  \widetilde {\rm Li}_3^*(\frac{x-xy}{1-xy}) + \widetilde {\rm Li}_3^*(xy) -  \widetilde {\rm Li}_3(\frac{x-xy}{1-xy}) +\widetilde {\rm Li}_3(y) -
\widetilde {\rm Li}_3(xy) +
$$
$$
  \widetilde {\rm Li}_1(xy)(\widetilde {\rm Li}_2(x) 
+ \widetilde {\rm Li}_2(y)) - \frac{1}{2}\widetilde\log^2(\frac{1-x}{1-xy}) \widetilde\log(\frac{1-y}{1-xy})
$$
where
$$
\widetilde {\rm Li}_3^*(x) := \quad \widetilde {\rm Li}_3(1) - \widetilde {\rm Li}_3(1-x)  -
\widetilde {\rm Li}_2(1)\widetilde {\rm Li}_1(x) - - \frac{1}{2}\widetilde \log(x) \widetilde {\rm Li}_1^2(1-x)
$$ 
To prove the formula we show by  a tedious calculation  
 that applying $\Delta$ to both parts of the formula we get zero. 
Thus the difference is a variation of Hodge-Tate structures over the $(x,y)$ space 
representing an element in ${\rm Ext}^1(\Q(0), \Q(3))$. 
Therefore, thanks to the Griffiths transversality condition, it  is a constant variation. 
So it remains to check that it is zero at one point. 
Consider the specialization at the $y=0$ divisor. Then the right hand side is clearly zero, 
and it is easy to check the left hand side is also zero. This proves the formula. 
  
Now lemma \ref{9.24.00.1} shows how to express $\widetilde {\rm Li}_{1,1,1}$ via 
$\widetilde {\rm Li}_{1,2}$, $\widetilde {\rm Li}_{2,1}$ and products of logarithms 
and dilogarithms. Further, $\widetilde {\rm Li}_{2,1}$ is reduced to 
$\widetilde {\rm Li}_{1,2}$ using one of the shuffle relations.

Substituting the expressions above to  formula
$$
\sum _{\sigma \in S_{3}} {{\widetilde {\rm Li}}}_{1,1,1}(x_{\sigma(1)},x_{\sigma(2)},
x_{\sigma(3)}) = 2 {{\widetilde {\rm Li}}}_{3}(x_1 x_2 x_3) 
$$
and applying to the framed Hodge-Tate structures the canonical Lie-period map (see [BD1]; it is 
the one which kills the products)
we get the following functional equation for the single-valued version 
$$
{\cal L}_3(z):= Re\Bigl( Li_3(z) -
 Li_3(z)log|z| + \frac{1}{3} Li_3(z)log^2(|z|)
$$
of the classical trilogarithm:
\begin{equation} \label{3.5*}
Sym\Biggr(  {\cal L}_3 (x) -   {\cal L}_3 (xy) - \frac{1}{3}  {\cal L}_3 (xyz) + 2  {\cal L}_3 (\frac{1-x}{1-y^{-1}}) +   {\cal L}_3 \Bigl(\frac{1-xyz}{1-z}\Bigr) - 
\end{equation}
$$
  {\cal L}_3 \Bigl(
\frac{(1-x)(1-z)}{1-y^{-1})(1-xyz)}\Bigr) +   {\cal L}_3 \Bigr(\frac
{(1-x)yz}{1-xyz}\Bigr) - {\cal L}_3(1) \Biggr) \quad = \quad 0
$$
Here $Sym$ is the symmetrization of variables $x,y,z$. It is
interesting that after symmetrization we get formula where all
the terms appear with the same coefficient $2$.

Similar considerations in the weight two give the classical 
 5-term functional equation for the dilogarithm written in symmetric form.

\section{Some applications and conjectures}

We explain in s. 7.2 how the theorewms from the introduction follow from our results. 
Assuming the motivic formalism we formulate in s. 7.3 
 conjectures relating multiple polylogarithm motives to the structure of
 the motivic Lie algebras of an 
arbitrary field $F$.

{\bf 1. The framed mixed Tate  motives  related  to multiple polylogarithms}. Let $F$ be  a number field,   $a_i \in F$ 
 and $a_0, a_{m+1}  \not = a_i$ for $i=1, ..., m$. Recall 
the $\Q(m)$-framed mixed Tate motive $I^{\cal M}(a_0; a_1, ..., a_m; a_{m+1})$ 
from definition \ref{16.2.01.1}. 
\begin{theorem} \label{14.1.01.1}   Let $F$ be  a number field,   $a_i \in F$ 
 and $a_0, a_{m+1}  \not = a_i$. Then 
for any embedding $\sigma: F \hookrightarrow \C$ the  Hodge realization of 
$I^{\cal M}(a_0; a_1, ..., a_m; a_{m+1})$ corresponding to $\sigma$
is isomorphic to 
$\widetilde I(\sigma(a_0); \sigma(a_1), ..., \sigma(a_m); \sigma(a_{m+1})$.
\end{theorem}

{\bf Proof}. Follows from the very definitions and the results of chapter 3. 

This together with  lemma \ref{1.1.01.2} immediately imply the following useful fact:

\begin{corollary} \label{3.12.01.1}
Assume the conditions of theorem \ref{14.1.01.1}. Then the coproduct of the 
framed mixed Tate motives $I^{\cal M}(a_0; a_1, ..., a_m; a_{m+1})$ 
is given  by the same formulas as in chapter 6 for their Hodge counterparts. 
\end{corollary}

\begin{corollary} \label{3.12.01.2} Assume that  and $a_i^N=1$ and $a_i \not = 1$ for $1 \leq i \leq m$. 
Then 
$I^{\cal M}(0; a_1, ..., a_m; 1)$ is a mixed Tate motive over the scheme $S_N$. 
\end{corollary}

{\bf Proof}. Follows from corollary \ref{3.12.01.1} 
and the definitions of chapter 3. 


{\bf 2. The theorems from the introduction.} Theorems \ref{3.6.01.3} and \ref{3.6.01.1} 
are special 
cases of theorems \ref{3.6.01.33} and \ref{th1.221}. The second part of theorem 
\ref{3.6.01.33} follows 
from proposition \ref{3.1.01.1} and formula (143) combined with (68). 

The proof of the first part of theorem \ref{3.6.01.33} follows from propositions 
\ref{12.30.00.2},  \ref{1.1.01.1}
following  the lines  of proof of theorem \ref{CP1}. 

{\bf Proof of theorem \ref{th1.221}}. The homomorphism 
${\cal Z}_{\bullet}^{\cal H}(\mu_N) \lra {\rm Gr}^W_{\bullet}\widetilde {\cal Z}(\mu_N)$ 
 is given by setting
$$
Li^{\cal H}_{n_1, ..., n_m}(\zeta_N^{\alpha_1}, ..., \zeta_N^{\alpha_m}) \lms 
(2\pi i)^{-w}Li_{n_1, ..., n_m}(\zeta_N^{\alpha_1}, ..., \zeta_N^{\alpha_m})
$$
Observe that the number on the right is well defined modulo the lower weight 
values of multiple polylogarithms at $N$-th roots of unity. 
This follows from our explicit description of the period matrix of the 
Hodge-Tate structure on the left which is established in proposition 5.7. 
More precisely, this proposition shows the following. Let $P$ be a period matrix 
describing the Hodge-Tate structure $Li_{n_1, ..., n_m}(\zeta_N^{\alpha_1}, ..., \zeta_N^{\alpha_m})$. Let $I$ be the matrix of the operator acting on the De Rham basis 
by $h_{\cal A} \lms 
(2\pi i)^{-|{\cal A}|}h_{\cal A}$. Then by proposition 5.7 combined with formula (59) 
the  entries of the matrix 
$PI$ belong to $W_{2w}\widetilde {\cal Z}(\mu_N)$, and applying to 
the entries the canonical projection 
$W_{2w}\widetilde {\cal Z}(\mu_N) \lra {\rm Gr}^W_{2w}\widetilde {\cal Z}(\mu_N)$ we 
get a matrix which has only one non zero entry, that is its maximal period 
$(2\pi i)^{-w}Li_{n_1, ..., n_m}(\zeta_N^{\alpha_1}, ..., \zeta_N^{\alpha_m})$. Indeed, 
all the entries of the 
matrix $P$ except $Li_{n_1, ..., n_m}(\zeta_N^{\alpha_1}, ..., \zeta_N^{\alpha_m})$ are 
 multiple polylogarithms at $N$-th roots of unity of the weight $w' <w$, 
multiplied by an appropriate  power $2\pi i$. 

Now let $H:= H_1 \oplus ... \oplus H_k$ and  $H':= H'_1 \oplus ... \oplus H'_l$,
 where each of the summands is a Hodge-Tate structure of  type 
$Li^{\cal H}_{n_1, ..., n_m}(\zeta_N^{\alpha_1}, ..., \zeta_N^{\alpha_m})$. 
Suppose that we have an equivalence  $H \sim H'$. By 
 lemma \ref{2.3.01.2}  this means that there is a minimal subquotient $\overline H$ of both 
$H$ and $H'$  which is  equivalent to $H$ and $H'$. The period matrix of 
 $\overline H$, being reduced to 
${\rm Gr}^W_{2w}\widetilde {\cal Z}(\mu_N)$,  has only one non zero entry.  
Therefore this entry 
is the period $H$ as well of  $H'$ with respect to the framings.   
So these periods coincide. 
This means that
 the map above is well defined on equivalence classes.  The theorem follows. 

{\bf 3. The depth filtration on the motivic Lie algebra of a field $F$ 
and multiple polylogarithms}. {\it The set up}. 
In this subsection we will 
assume the existence of the {\it abelian} category ${\cal M}_T(F)$ of mixed Tate motives over
 an arbitrary field $F$ with all the standard properties. In particular we assume 
that the category ${\cal M}_T(F)$ is a mixed Tate category. Denote by 
 $L(F)_{\bullet}$  the fundamental Lie algebra of the category ${\cal M}_T(F)$, by 
${\cal L}(F)_{\bullet}$ its dual, called the motivic Lie coalgebra of $F$, and by 
${\cal U}(F)_{\bullet}$ the dual to the universal enveloping algebra of $L(F)_{\bullet}$. 

Moreover we assume that for any  $a_1, ..., a_m \in F^*$ there exists an 
object ${\rm I}^{{\cal M}}_{n_1, ..., n_m}(a_1, ..., a_m)$ 
of the category ${\cal M}_T(F)$ framed by $\Q(0)$ and $\Q(w)$, called the motivic 
multiple polylogarithm. If $F$ is a number field the results of this paper  
already give  all this.

In particular we assume that for any field $F$ the coproduct of the 
framed mixed Tate motives ${\rm I}^{{\cal M}}_{n_1, ..., n_m}(a_1, ..., a_m)$ is given by the 
same formula as in theorem \ref{ur8-4,3}. 

{\it The motivic logarithm}. We have an  isomorphism ${\cal U}(F)_{1} \stackrel{\sim}{=}
 F^{*}\otimes \Q$. It is materialized by the motivic logarithm providing an isomorphism 
$$
\log^{{\cal M}}:  F^* \otimes \Q \lra  {\rm Ext}_{{\cal M}_T(F)}^1(\Q(0), \Q(1)) 
 \stackrel{\sim}{=}  {\cal U}(F)_{1}; \qquad a \lms 
\log^{{\cal M}}(a)
$$  
Its l-adic realization is given by the Kummer extension, and if $F \subset \C$ its Hodge 
realization is  $\widetilde \log (a)$.  

{\it The depth filtration(s)}. It follows from this that motivic multiple polylogarithms 
of depth $\leq m$ provide a  filtration ${\cal F}^{\cal D}$ 
by coideals  in  ${\cal U}(F)_{\bullet}$ indexed by integers $m \geq 0$. Set
$$
{\cal F}_{0}^{\cal D}{\cal U}(F)_{\bullet} : = 
\quad \mbox{the subspace spanned over $\Q$ by the products of $\log^{{\cal M}} (a)$},  
$$
\begin{equation} \label{1.10.01.2}
{\cal F}_{m}^{\cal D}{\cal U}(F)_{\bullet}:= \quad \mbox{the subspace spanned by 
${\rm I}^{{\cal M}}_{n_1, ..., n_k}(a_1, ..., a_k)$, $a_i \in F^*$,  $1 \leq k \leq m$}
\end{equation}
Since 
$\widetilde \log (a) = \widetilde {\rm I}_{\{v_a\}}(0; 0; a)$ is 
the depth zero multiple polylogarithm, this agrees with (\ref{1.10.01.2}). 

The universality conjecture  can be reformulated as follows: 
\begin{conjecture} \label{1.10.01.3} Let $F$ be an arbitrary field. Then  
$
\cup_m {\cal F}_{m}^{\cal D}{\cal U}(F)_{\bullet} \quad = \quad {\cal U}(F)_{\bullet}
$. So  every framed mixed Tate motive over $F$ is equivalent to a subquotient of 
the motivic torsor of path ${\cal P}^{{\cal M}}({\rm Spec} F(x); v_0, v_1)$ between tangential base points at $0$ and $1$. 
\end{conjecture}

{\it The universality conjecture}. 
It tells us 
(see conjecture 17a) in [G1]) that 
every 
framed mixed Tate motive over $F$ is equivalent to a $\Q$-linear combination of 
the framed motives ${\rm I}^{{\cal M}}_{n_1, ..., n_m}(a_1, ..., a_m)$ with $a_i \in F^*$.

The following result, which is a motivic version theorem \ref{th2.7}, 
 shows  that all  n-framed mixed Tate motives but perhaps a countable set are 
given by multiple polylogarithm motives. This is a strong support for the universality conjecture.

\begin{theorem} \label{3.2} Let us assume the motivic formalism. 
Let ${\Bbb V}$ be a variation of n-framed mixed Tate motives over a connected rational variety $Y$. 
Then for any two points $y_1,y_2 \in Y$ the difference $V_{y_1} - V_{y_2}$
is a sum of multiple polylogarithms motives.
\end{theorem}

{\bf Proof}. One can suppose without the loss of generality that there is a rational curve $X$ passing through $y_1$ and $y_2$. Let ${\cal P}^{{\cal M}}(X; y_1,y_2)$ be the torsor of motivic paths from $y_1$ to $y_2$ on $X$. The crucial fact is that any its subquotient is equivalent as framed mixed Tate motive to a hyperlogarithmic one. There is a morphism of mixed Tate motives $p_{y_1,y_2}^X :V_{y_1}\otimes
 {\cal P}^{{\cal M}}(X;y_1,y_2)
\rightarrow V_{y_2}$ (parallel transport along paths from $y_1$ to $y_2$
on $X$). Let ${\cal A}$ be the kernel of the action of ${L}(F)_{\bullet}$
on   ${\cal P}^{{\cal M}}(P^1; v_0, v_1)$. Then ${\cal P}^{{\cal M}}(X; y_1,y_2)$
  is a trivial  ${\cal A}$ - module. Therefore any $p \in  {\cal P}^{{\cal M}}(X; y_1,y_2)$ defines an isomorphism of ${\cal A}$ - modules $V_{y_1}\otimes p 
\rightarrow V_{y_2}$. This is a reformulation of the statement of the theorem.

{\it The main conjecture}. 
The depth filtration on ${\cal U}(F)_{\bullet}$ induces the depth filtration on the corresponding 
Lie coalgebra 
${\cal L}(F)_{\bullet}$. In particular ${\cal F}_{0}^{\cal D}{\cal L}(F)_{\bullet} = 
{\cal L}(F)_{1} = F^* \otimes \Q$. 

Consider the ideal of the Lie algebra $L(F)_{\bullet}$ given by 
$$
 I(F)_{\bullet}:= \oplus^{\infty}_{n=2}  L(F)_{-n}
$$
We define the {\it depth filtration} ${\cal F}^D$ 
on the Lie algebra $L(F)_{\bullet}$ as an increasing  filtration  
indexed by integers $m \leq 0$ and given by the powers of the ideal $I(F)_{\bullet}$
\begin{equation} \label{1.10.01.1}
{\cal F}^{D}_0L(F)_{\bullet} = L(F)_{\bullet}; \quad 
{\cal F}^{D}_{-1}L(F)_{\bullet} = I(F)_{\bullet} 
; \quad {\cal F}^{D}_{-m-1}L(F)_{\bullet} = [I(F)_{\bullet}, {\cal F}^{D}_{-m}L(F)_{\bullet}]
\end{equation}
Thus there are two filtrations on the Lie coalgebra ${\cal L}(F)_{\bullet}$:  the 
dual to  depth 
filtration (\ref{1.10.01.1}) and the 
filtration  by the depth of multiple polylogarithms which is induced  by (\ref{1.10.01.2}). 
We conjecture that they coincide:
 
\begin{conjecture} \label{13.1.01.1}
The dual to the depth filtration (\ref{1.10.01.1}) coincides with the 
filtration induced by (\ref{1.10.01.2}). 
\end{conjecture}
This conjecture, of course, implies conjecture \ref{1.10.01.3}. 
It also implies, when $F$ is a number field, Zagier's conjecture. 


{\it The role of classical polylogarithms}. Let  ${\cal B}_{n}(F)$ be 
the $\Q$-vector space in 
${\cal L}(F)_{n}$ spanned by the classical $n$-logarithm framed motives 
${\rm Li}^{{\cal M}}_n(a)$, $a \in F^*$. (This definition differs 
from the one given in [G6-7], also the corresponding groups  expected to be isomorphic.) 
Denote by  $H^{i}_{(n)}I(F)_{\bullet}$  the degree $n$ part of $H^{i}I(F)_{\bullet}$.
The following conjecture was stated in [G1]:

\begin{conjecture} \label{7.3}
a)  $H^{1}_{(n)}I(F)_{\bullet} \cong {\cal
B}_{n}(F)$ for $n\geq 2$, i.e. $ I(F)_{\bullet}$ 
is generated as a graded Lie algebra by the spaces ${\cal
B}_{n}(F)^{\vee}$ sitting in degree $-n$.

b) $I(F)_{\bullet}$  is a free graded (pro) -
Lie algebra.
\end{conjecture}

Conjecture \ref{13.1.01.1} obviously implies the part a) of conjecture \ref{7.3}.

\vskip 3mm \noindent
{\bf REFERENCES}
\begin{itemize}
\item[{[BD1]}] Beilinson A.A., Deligne P.: {\it Interpr\'etation motivique de la conjecture de Zagier reliant polylogarithmes et r\'egulateurs}. Motives (Seattle, WA, 1991), 97--121, 
Proc. Sympos. Pure Math., 55, Part 2,
AMS, Providence, RI, 1994.
\item[{[BD2]}] Beilinson A.A., Deligne P.: 
{\it Motivic polylogarithms and Zagier's conjecture} 
Manuscript, version of 1992. 
\item[{[Be1]}] Beilinson A.A.: {\it Higher regulators and values of $L$-functions}. (Russian) Current problems in mathematics, Vol. 24, 181--238, Itogi
Nauki i Tekhniki, Akad. Nauk SSSR, Vsesoyuz. Inst. Nauchn. i Tekhn. Inform., Moscow, 1984. 
\item[{[BMS]}] Beilinson A.A., MacPherson R.D. Schechtman V.V: {\it Notes on motivic cohomology}. Duke Math. J., 1987 vol 55 p. 679-710
\item[{[BGSV]}] Beilinson A.A., Goncharov A.A., Schechtman V.V., Varchenko A.N.: {\it Aomoto dilogarithms, mixed Hodge structures and motivic cohomology of a pair of triangles in the plane}, the Grothendieck Feschtrift, Birkhauser, vol 86, 1990, p. 135-171.
\item[{[Bo]}] Borel A.: {\it Cohomologie de ${\rm SL}\sb{n}$ et valeurs de fonctions zeta aux points entiers} Ann. Scuola Norm. Sup. Pisa Cl.
Sci. (4) 4 (1977), no. 4, 613--636.
\item[{[Br]}] Broadhurst D.J., {\it On the enumeration of irreducible $k$-fold sums and their role in knot theory and field theory} Preprint hep-th/9604128. 
\item[{[Ch]}] Chen K.-T.: {\it Iterated path integrals}. 
Bull. Amer. Math. Soc. 83 (1977), 323-338. 
\item[{[D]}] Deligne P.: {\it Le group fondamental de la droite projective moine trois points}, In: Galois groups over $\Q$. Publ. MSRI, no. 16 (1989) 79-298.   
\item[{[D2]}] P. Deligne, {\it A letter to D. Broadhurst}, June 1997.
\item[{[D3]}] P. Deligne, {\it Letter to the author}. July 25, 2000.
\item[{[D4]}] P. Deligne, {\it Th\'eorie de Hodge II, III}. Publ. IHES No. 40 (1971), 5--57;
No. 44 (1974), 5--77.
\item[{[Dr]}] Drinfeld V.G.: {\it On quasi-triangular quasi-Hopf algebras and some group related to clousely associated with Gal$(\overline{\Q}/\Q)$}. Leningrad Math. Journal, 1991. (In Russian).
\item[{[E]}] L. Euler: "Opera Omnia," Ser. 1, Vol XV, Teubner,
Berlin 1917, 217-267.  
\item[{[ES]}] Etingof P., Schiffmann O.: {\it Lectures on Quantum groups}.
 International press, 1998. 
\item[{[G0]}] Goncharov A.B.: {\it Multiple $\zeta$-numbers,
    hyperlogarithms and mixed Tate motives}, Preprint MSRI 058-93, June 1993.
\item[{[G1]}] Goncharov A.B.: {\it Polylogarithms in arithmetic and geometry}, 
Proc. of the International Congress of Mathematicians, Vol. 1, 2
(Zurich, 1994), 374--387, Birkhauser, Basel, 1995.
\item[{[G2]}] Goncharov A.B.: {\it The double logarithm and Manin's
complex for modular curves}. Math. Res. Letters,
 vol. 4. N 5 (1997), pp. 617-636. 
\item[{[G3]}] Goncharov A.B.: {\it Multiple polylogarithms, cyclotomy and modular complexes
},     Math. Res. Letters, 
 vol. 5. (1998), pp. 497-516. K-theory e-print archive, www.math.uiuc.edu/K-theory/ N 297.
\item[{[G4]}] Goncharov A.B.: {\it The dihedral Lie algebras and galois symmetries of 
$\pi_1^{(l)}({\Bbb P}^1 \backslash 0,1, \infty)$}. Duke Math. J.  (2001). 
Alg-geom e-print archive, math.AG/0009121.
\item[{[G5]}] Goncharov A.B.: {\it Multiple $\zeta$-values, Galois groups and geometry of 
modular varieties} Proc. of the Third European Congress of mathematicians, Barcelona (2000). Alg-geom e-print archive, math.AG/0005069.
\item[{[G6]}] Goncharov A.B.: {\it Polylogarithms and motivic Galois groups}, Motives (Seattle, WA, 1991), 43--96, Proc. Sympos. Pure Math., 55, Part 2,
Amer. Math. Soc., Providence, RI, 1994.
\item[{[G7]}] Goncharov A.B.: {\it Geometry of configurations, Polylogarithms and motivic cohomology}.  Adv. in Math.  114 (1995), no. 2, 197--318.
\item[{[G8]}] Goncharov A.B.: {\it Volumes of hyperbolic manifolds and mixed Tate motives} J. Amer. Math. Soc. 12 (1999) N2, 569-618. math.alg-geom/9601021
\item[{[G9]}] Goncharov A.B.: {\it Mixed elliptic motives}. Galois representations in arithmetic algebraic geometry (Durham, 1996), 147--221, London
Math. Soc. Lecture Note Ser., 254, Cambridge Univ. Press, Cambridge, 1998. 
K-theory e-print archive, www.math.uiuc.edu/K-theory/ N 228.
\item[{[G10]}] Goncharov A.B.: {\it Multiple polylogarithms, higher
cylotomy, and geometry of modular varieties}. To appear. 
\item[{[Gr]}] Grothendieck A.:{\it Esquisse d'un programme.} 
Mimeographed note (1984).
\item[{[GS]}] Gelfand I.M., Shilov G. E.: {\it Generalized functions}. Vol. 1. Properties and operations. 
Academic Press, New York-London, 1964 [1977].
\item[{[Ha]}] Hanamura M.: {\it Mixed motives and algebraic cycles}. I. 
Math. Res. Lett. 2 (1995), no. 6, 811--821. 
\item[{[H]}] Hain R., {\it Classical polylogarithms}. Proceedings of AMS Research Summer Conference ``Motives'' 1994, vol. 2, 1-42.
\item[{[H1]}] Hain R., 
{\it The geometry of mixed Hodge structure on $\pi_1$}.  
 Algebraic geometry, Bowdoin, 1985
(Brunswick, Maine, 1985), 247--282, Proc. Sympos. Pure Math., 46, Part 2, Amer. Math. Soc., Providence, RI, 1987. 
\item[{[HM]}] Hain R., Matsumoto M..: {\it Weighted Completion of Galois Groups and Some Conjectures of Deligne}. math.AG/0006158. 
\item[{[HZ]}] Hain R., Zucker S.: {\it Unipotent variations of mixed Hodge structure}. Inventiones Math., 88, (1987) 83-124. 
\item[{[Hof]}] Hoffman M.: Multiple harmonic series, Pacific J. Math., 
152 (1992), 275-290.
\item[{[Hu]}] Huber, A.: {\it Realizations of Voevodsky's motives}. J. Algebraic Geom. 9 (2000), no. 4, 755--799. 
\item[{[Ih]}] Ihara Y.: {\it Braids, Galois groups, and some arithmetic functions}. 
Proc. of the Int. Congress of Mathematicians, Vol. I, II
(Kyoto, 1990), 99--120, Math. Soc. Japan, Tokyo, 1991.
\item[{[Ih1]}] Ihara Y.: {\it Some arithmetical aspects of Galois 
action on the pro-$p$ fundamental group of 
$\widehat \pi_1({\Bbb P}^1 - \{0,1,\infty\})$}, Preprint RIMS-1229, 1999.
\item[{[K1]}] M. Kontsevich, {\it Private communication} (Fall 1991). 
\item[{[K2]}] M. Kontsevich, {\it Operads and motives in deformation quantization},
Lett. Math. Phys., 48 (1999) N1, 35.
\item[{[Kr]}] Kreimer D. {\it Renotmalisation and knot theory} 
 J. Knot Theory Ramifications 6 (1997), no. 4, 479--581.
\item[{[Ku]}] Kummer E.E.: {\it Uber die Transcendeten, welche aus wiederholten Integrationen rationaler Formeln entstehen}, J. reine ang. Math. 21 (1840), pp. 74--90, 193--225, 328--371 ( = C.P., vol.11, N 15).
\item[{[L]}] Levine, M: {\it Mixed motives}. Mathematical Surveys and Monographs, 57. American Mathematical Society, Providence, RI, 1998
\item[{[L1]}] Levine, M: {\it Tate motives and the vanishing conjectures for algebraic 
$K$-theory}. Algebraic $K$-theory and algebraic topology (Lake
Louise, AB, 1991), 167--188, NATO Adv. Sci. Inst. Ser. C Math. Phys. Sci., 407, 
Kluwer Acad. Publ., Dordrecht, 1993.
\item[{[LD]}] Lappo-Danilevsky: {\it Resolution algorithmique des probl\`em es reguliers de Poincare et de Riemann}. J.Soc. Phisico-Mathematique de St.Petersbourg. vol 3, 1911.
\item[{[M]}] Morgan J.: {\it The algebraic topology of smooth algebraic varieties}. 
Publ. IHES. No 48 (1978), 137-204.
\item[{[R]}] Racinet G.: {\it S\'eries g\'en\'eratrices non commutatives de polyzetas et associateurs de Drinfel'd}. Th\'ese, (2000). 
\item[{[Sh]}] Shiho, A.:. {\it Crystalline fundamental groups and $p$-adic Hodge theory}. 
In Arithmetic and Geometry of Algebraic Cycles, 
381--398, CRM Proc. Lecture Notes, 24, Amer. Math. Soc., Providence, RI, 2000. 
\item[{[Z]}] Zagier D.: {\it Values of zeta functions and their applications}. First European Congress of Mathematics, Vol. II (Paris, 1992), 497--512, Progr.
Math., 120, Birkhauser, Basel, 1994. 
\item[{[Z1]}] Zagier D.:{\it Periods of modular forms, traces of Hecke operators, and multiple $\zeta$-values} Research into automorphic 
forms and $L$-functions (Japanese) (Kyoto 1992), Surikaisekikenkyusho Kokyuroku No. 843 (1993), 162-170.
\item[{[V]}] Voevodsky V.:  {\it Triangulated category of motives over a field} in 
Cycles, transfers, and motivic homology theories. Annals of Mathematics Studies,
143. Princeton University Press, Princeton, NJ, 2000. 
\item[{[V2]}] Voevodsky V.:  {\it Stable homotopy theory} To appear.
\item[{[Vol]}] Vologodsky V.:  
{\it The mixed Hodge structure on the fundamental group and 
p-adic integration}. Thesis, Harvard University, 2001.
\item[{[We]}] Weil A.: {\it Elliptic functions according to Eisenstein and Kronecker}. Springer Verlag, 1976. 
\item[{[W1]}] Woitkowiyak Z.: 
{\it The basic structure of polylogarithmic equations}. Ch.10 in 
Structural properties of polylogarithms. 
Edited by L. Lewin. Mathematical Surveys and Monographs, 37. 
AMS, Providence, RI, 1991.
\item[{[W2]}] Woitkowiyak Z.: {\it Cosimplicial objects in algebraic geometry}. 
Algebraic $K$-theory and algebraic topology (Lake Louise, AB, 1991),
287--327, NATO Adv. Sci. Inst. Ser. C Math. Phys. Sci., 407, 
Kluwer Acad. Publ., 
Dordrecht, 1993. 

\end{itemize}

\end{document}